\documentclass[11pt]{amsart}

\usepackage{fullpage, graphicx, enumitem, mathtools,stmaryrd}
\SetSymbolFont{stmry}{bold}{U}{stmry}{m}{n}
\usepackage{url}
\usepackage{verbatim}
\usepackage{bbm}
\usepackage{amssymb,amsmath}
\usepackage{amsfonts,savesym,bm,amsthm}
\usepackage{xcolor}
\usepackage[mathscr]{eucal}
\usepackage[all]{xy}
\usepackage{tikz}
\usepackage{etoolbox}

\usepackage{float}

\definecolor{hot}{RGB}{65,105,225}
\usepackage[pagebackref=true,colorlinks=true, linkcolor=hot ,  citecolor=hot, urlcolor=hot]{hyperref}

\newtheorem{theorem}{Theorem}[section]
\newtheorem{lemma}[theorem]{Lemma}

\newtheorem{theorem-definition}[theorem]{Theorem-Definition}
\newtheorem{corollary}[theorem]{Corollary}
\newtheorem{proposition}[theorem]{Proposition}

\newtheorem{definition-theorem}[theorem]{Definition-Theorem}
\newtheorem{algorithm}[theorem]{Algorithm}

\newtheorem{theorem-defintion}[theorem]{Theorem-Definition}
\newtheorem{corollary-definition}[theorem]{Corollary-Definition}
\newtheorem{definition-proposition}[theorem]{Definition-Proposition}

\theoremstyle{definition}
\newtheorem{example}[theorem]{Example}
\newtheorem{definition}[theorem]{Definition}
\newtheorem{defn}[theorem]{Definition}

\newtheorem{notation}[theorem]{Notation}
\newtheorem{remark}[theorem]{Remark}

\newtheorem{convention}[theorem]{Convention}
\newtheorem{subs}[theorem]{ }

\numberwithin{equation}{section}

\def\be{\begin{equation}}
	\def\ee{\end{equation}}

\def\bZ{\mathbb{Z}}

\def\bQ{\mathbb{Q}}

\def\ord{\text{ord}}

\def\divi{\mathrm{div}}

\def\bs{{\bm s}}

\def\Zmot{Z^{\mathrm{mot}}}

\def\codim{\mathrm{codim}\;}

\def\brT{\{\!\{ T_1^{\pm},\ldots,T_r^{\pm}\}\!\}}
\def\td{\mathrm{td}}
\def\cR{\mathcal R}
\def\cD{\mathcal D}

\newcommand{\C}{\mathbb C}
\newcommand{\Q}{\mathbb Q}
\newcommand{\R}{\mathbb R}
\newcommand{\Z}{\mathbb Z}
\DeclareMathOperator{\Div}{div}

\newcommand{\red}{\mathrm{red}}

\title{A remarkable subset of polar hyperplanes of the multivariable motivic zeta function}

\author{Marko Čmrlec}
\address{Department of Mathematics, KU Leuven, Celestijnenlaan 200B, 3001 Leuven, Belgium;  BCAM, Mazarredo 14, 48009 Bilbao, Spain.}
\email{marko.cmrlec@kuleuven.be}

\author{Shokhzod Kurokboev}
\address{Department of Mathematics, KU Leuven, Celestijnenlaan 200B, 3001 Leuven, Belgium.}
\email{shokhzod.kurokboev@kuleuven.be}

\author{Quan Shi}
\address{Department of Mathematical Sciences, Tsinghua University, Beijing, 100084, P. R. China.}
\email{shiq24@mails.tsinghua.edu.cn / thusq20@gmail.com}

\author{Huaiqing Zuo}
\address{Department of Mathematical Sciences, Tsinghua University, Beijing, 100084, P. R. China.}
\email{hqzuo@mail.tsinghua.edu.cn}

\begin{document}
	
	\begin{abstract}
		For a tuple of complex multivariate polynomials $F=(f_1,\ldots, f_r)$ we find a remarkable subset of polar hyperplanes
		of the multivariable motivic zeta function. They are combinatorially determined by any log resolution of the product $f$ of the polynomials, admit an intrinsic interpretation in terms of arc contact loci, can provide non-trivial information beyond the remarkable numbers of $f$, and admit a generalization to tuples of rational functions. 
	
	The results of this paper, done without AI,  are complemented by an Appendix by a subset of the co-authors using AI and proving the strong monodromy conjecture for remarkable hyperplanes. It gives an alternative proof of the conjecture for the remarkable numbers of $f$.

	\end{abstract}
	\maketitle
	\raggedbottom
	
	\tableofcontents

	\section{Introduction}
	
	Let $X$ be an $n$-dimensional smooth irreducible complex variety and $F = (f_1,\ldots,f_r) : X \to \mathbb A^r$ be a regular map such that $f = f_1\cdots f_r$ is nonzero and non-invertible. The multivariable monodromy conjecture, cf. \cite{Bud15}, predicts that the polar locus of the multivariable motivic zeta function $Z_F^{\mathrm{mot}}(T_1,\ldots,T_r)$, a set of hyperplanes in $\mathbb R^r$, gives rise to a subset of the monodromy support and the zero locus of the Bernstein-Sato ideal. These invariants, as generalizations of classical ones from the case when $r=1$, capture the motivic and complex nature of $f_1,\ldots,f_r$ simultaneously. %The conjecture thus relates the two isolated aspects. 

	Any log resolution of $f$ gives a set of candidate polar hyperplanes. However, most of them are not actual polar hyperplanes, and the cancellations remain mysterious. Consequently, any lower bound for the polar locus is useful. In \cite{BPSZ}, when $r=1$, a remarkable subset of poles of the motivic zeta function was found, determined by the combinatorial data of any log resolution of $f$. 
	
	In this paper, a remarkable subset of polar hyperplanes is found for a tuple of not-necessarily-effective divisors. This includes the case of regular maps $F$ as above, as well as the case of a single rational function $f$, not necessarily regular. 
	
	An interesting by-product of our analysis is one of the original motivations behind this paper. Namely, the remarkable hyperplanes of a splitting $F$ of a regular function $f = f_1\cdots f_r$ can highlight combinatorially interesting exceptional divisors over the zero locus of $f$ which, by some combinatorial mishap, do not compute remarkable numbers of $f$. For example, all rupture components of a multibranch plane curve singularity are detectable via the remarkable hyperplanes, whereas, in contrast with the unibranch case, there are examples such that not every rupture component gives rise to a remarkable number.

	To state our main result, we introduce the multivariable motivic zeta function of an ordered tuple of divisors and its computation via a log resolution, which is a generalization of \cite[Theorem 2.2.1]{DL98}. Let $\underline D = (D_1,\ldots,D_r)$ be a tuple of divisors on $X$, not necessarily effective. For a divisor $D' =\sum_Z a_Z Z$, we write $D'_{\mathrm{red}}:=\sum_{a_Z\neq0}Z$
	and identify it with its support. Let $\rho: Y \to X$ be a log resolution of $\bigcup_{i=1}^r (D_i)_{\mathrm{red}}$ and such that $\rho^{-1}(\bigcup_{i=1}^r(D_i)_{\mathrm{red}}) = \bigcup_{i\in S} E_i$ is a simple normal crossings divisor. Suppose the relative canonical divisor $K_{\rho} = \sum_{i\in S}(\nu_i -1) E_i$ and $\rho^*D_i = \sum_{j\in S} N_{ji} E_j$ for $i=1,\ldots,r$.

	The multivariable motivic zeta function of $\underline D$ is the generating series of motivic volumes of multi-contact loci of $\underline D$, i.e.
	\begin{displaymath}
			Z_{\underline D}^{\mathrm{mot}}(T_1,\ldots,T_r) := \sum_{ (\alpha_1,\ldots,\alpha_r) \in \mathbb Z^r} \mu(\mathscr X_{(\alpha_1,\ldots,\alpha_r)}(X,\underline D)) \cdot T_1^{\alpha_1}\cdots T_r^{\alpha_r} \in \widehat{\mathcal{M}}_{\mathbb{C}}\{T_1^{\pm},\dots,T_r^{\pm}\},
	\end{displaymath}  
	where $\widehat{\mathcal{M}}_{\mathbb{C}}\{T_1^{\pm},\dots,T_r^{\pm}\}$ is the group of bilateral power series with coefficients in the completed Grothendieck ring of \ref{subsGr}. 
		Via the log resolution $\rho$, it can be expressed as
	\begin{equation}\label{Zmot_DL1}
		Z^{\mathrm{mot}}_{\underline{D}}(T_1,\ldots,T_r)=\sum_{I \subset S} [E_I^{\circ}] \prod_{j\in I}\frac{(\mathbb{L}-1)\mathbb{L}^{-\nu_j} T_1^{N_{j1}} \cdots T_r^{N_{jr}} }{1-\mathbb{L}^{-\nu_j}T_1^{N_{j1}} \cdots T_r^{N_{jr}}},
	\end{equation}
	where $E_I^\circ = (\bigcap_{i\in I}E_i) \setminus (\bigcup_{j\notin I} E_j)$. Equation \eqref{Zmot_DL1} is subtle since $\widehat{\mathcal M}_{\mathbb C}\{T_1^{\pm},\ldots,T_r^{\pm}\}$ is not a ring. This will be explained in \ref{convergent}. This expression leads to the notion of polar hyperplanes for the zeta function, intuitively coming from the denominators after performing the allowed cancellations in (\ref{Zmot_DL1}), see \ref{subsPol}. 	All polar hyperplanes are of type $\sum_{i=1}^r N_{ji} \cdot s_i + \nu_j = 0$ for some $j\in S$.
	
	\subs{\bf Remarkable hyperplanes.} \label{Rmk_subs}We associate some more numerical data to $\rho$.
	\begin{definition}
		For $I\subset S$, define $E_I = \bigcap_{i\in I} E_i$ and $E_I^\circ = E_I \setminus (\bigcup_{j\notin I} E_j)$. For a stratum $E_I \neq \varnothing$, define 
		\begin{align*}
			N_I & := \sum_{j\in I} \mathbb Z \cdot (N_{j1},\ldots,N_{jr}) \subset \mathbb Z^r,\\
			F_I & := \bigcap_{j\in I} \{(s_1,\ldots, s_r) \in \mathbb R^r \mid \sum_{i=1}^r N_{ji} \cdot s_i + \nu_j > 0\}.
		\end{align*}
	\end{definition}
	\begin{convention}
		\begin{enumerate}
			\item When using $N_I$ or $F_I$, we always assume $E_I \neq \varnothing$.
			\item The lattices denoted by $L$ are assumed to have full rank, i.e. $\mathrm{rank}\, L = r$; the groups $N_I$ need not have full rank. We use $\mathbb N=\mathbb Z_{\geq0}$.
			\item We will often use multi-index notation in the rest of the paper. For example, $\bm T := (T_1,\ldots,T_r)$, and for $\bm a = (a_1,\ldots,a_r) \in \mathbb Z^r$, we denote $\vert \bm a\vert := |a_1|+\cdots+|a_r|$ and $\bm T^{\bm a} = T_1^{a_1}\cdots T_r^{a_r}$.
		\end{enumerate}
	\end{convention}
	\begin{definition}
		Let $L \subset \mathbb Z^r$ be a sublattice   and $\bm a\in \mathbb Z^r$. We define the \textit{$(\bm a,L)$-lct polytope} as
		\begin{displaymath}
			\mathrm{lct}^{\bm a,L}(\underline D) = \bigcap_{\bm a \in N_I+L, E_I\neq \varnothing} F_I.
		\end{displaymath}
		We adopt the convention that $\bigcap_{\varnothing} F_I = \mathbb R^r$.
		Denote by $\mathrm{lctf}^{\bm a,L}(\underline D)$ the set of hyperplanes spanned by facets of the closure of $\mathrm{lct}^{\bm a,L}(\underline D)$. We call a hyperplane in $\mathbb R^{r}$ \textit{remarkable for $\underline D$} if it belongs to $\mathrm{lctf}^{\bm a,L}(\underline D)$ for some $(\bm a,L)$. Denote by $\mathrm{Rmk}(\underline D)$ the set of remarkable hyperplanes.
		%		
		%		hyperplanes spanned by the facets of $\mathrm{lct}^{\bm a,L}(\underline D)$ are called \textit{remarkable hyperplanes of $\underline D$}.
	\end{definition}
	\begin{theorem}\label{ThmA}
		Let $X$ be a smooth irreducible complex variety and $\underline D = (D_1,\ldots,D_r)$ be a tuple of divisors on $X$, not necessarily effective. Then:
		\begin{enumerate}
			\item $\mathrm{Rmk}(\underline D)$ is independent of the choice of log resolution. \label{A1}
			\item The remarkable hyperplanes are polar hyperplanes of $Z_{\underline D}^{\mathrm{mot}}(T_1,\ldots,T_r)$.\label{A2}
			\item For $N \in \mathbb Z_{>0}$ divisible enough, we have
			\begin{displaymath}
				\mathrm{Rmk}(\underline D) = \bigcup_{\bm a\in \{0,\ldots,N-1\}^r } \mathrm{lctf}^{\bm a,(N\cdot \mathbb Z)^r}(\underline D).
			\end{displaymath} \label{A3}
		\end{enumerate}
	\end{theorem}
	\begin{remark}\label{rmk:enlarged-lct}$\;$
		\begin{enumerate}
			\item When $r = 1$ and $D_1$ is effective, then $\mathrm{Rmk}(\underline D)$ coincides with the set of negatives of remarkable numbers defined in \cite{BPSZ}.
			%, $\mathrm{Rmk}(\underline D)$ coincides with 
			%consists of the points $\{-\theta\}$, where $\theta$ runs through the remarkable numbers defined below.
			\item Taking $\bm a = \bm 0$, $N$ divisible enough, and $D_i \geq 0$ for all $i$, it is not difficult to see that the negative of the Euclidean interior of the usual lct polytope, cf. \cite{BVVW24}, is the intersection of $\mathbb R_{<0}^r$ and $\mathrm{lct}(\underline D) := \mathrm{lct}^{\bm 0,(N\cdot \mathbb Z)^r}(\underline D) = \bigcap_{i\in S} F_{\{i\}}$. We call $-\mathrm{lct}(\underline D)$ the \textit{enlarged lct polytope} and denote by $\mathrm{lctf}(\underline D)$ the set of hyperplanes spanned by facets of the closure of $\mathrm{lct}(\underline D)$. These hyperplanes are automatically remarkable. Moreover, if a facet of the usual lct polytope meets $\mathbb R_{>0}^r$ and spans the hyperplane $H$, then $-H\in\mathrm{lctf}(\underline D)$.
			%is contained in $\mathrm{lct}(\underline D) := \mathrm{lct}^{\bm 0,(N\cdot \mathbb Z)^r}(\underline D) = \bigcap_{i\in S} F_{\{i\}}$. However, $\mathrm{lct}^{\bm 0,(N\cdot \mathbb Z)^r}(\underline D)$ is much larger since the usual lct polytope is restricted to $\mathbb R_{\geq0}^r$. Our open sets, which may be unbounded, are more precisely open polyhedra. We call $-\mathrm{lct}(\underline D)$ the \textit{enlarged lct polytope} and denote by $\mathrm{lctf}(\underline D)$ the set of hyperplanes spanned by facets of the closure of $\mathrm{lct}(\underline D)$. These hyperplanes are automatically remarkable. 
			\item For $x\in X$, we can define the set of {\it local remarkable hyperplanes} $\mathrm{Rmk}_x(\underline D)$ by retaining only the strata $E_I$ whose image contains $x$. Equivalently, $\mathrm{Rmk}_x(\underline D)=\mathrm{Rmk}(\underline D|_U)$ for every sufficiently small neighborhood $U$ of $x$. Indeed, the images $\rho(E_I)$ are closed and there are only finitely many of them. %An intersection of the remarkable sets over all neighborhoods is not required.
		\end{enumerate}
	\end{remark}

	\subs{\bf Intrinsic interpretation via contact loci.} Similar to remarkable numbers, remarkable hyperplanes can be intrinsically read from codimensions of contact loci, see \ref{VP_subs} for definition of codimensions. 
	\begin{definition}
		Let $\eta \neq \mathbb R_{\geq 0} \cdot (\bm 0,1)$ be a ray in $\mathbb R^{r} \times \mathbb R_{\geq 0}$ and $(\bm a,b) \in \eta \setminus \{\bm 0\}$. The {\it hyperplane in $\mathbb R^r$ corresponding to $\eta$} is $\{\bm s \in \mathbb R^r \mid \bm a \cdot \bm s + b = 0\}$.
	\end{definition}
	\begin{theorem}\label{intrinsic_cl}\label{ThmB}
		With notation as in \ref{Rmk_subs}, let $\bm a\in\mathbb Z^r$
		and let $L\subseteq\mathbb Z^r$ be a full-rank sublattice. Set $\mathcal T = \{\bm p \in \bm a + L \mid \mathscr X_{\bm p}(X,\underline D) \neq \varnothing\}$. Let
		\begin{displaymath}
			K^{\bm a,L}_0(\underline D) := \bigcup_{\bm p\in \mathcal T} \big((\bm p,\mathrm{codim}\, \mathscr X_{\bm p}(X,\underline D)) + \mathbb R_{\geq 0} \cdot (\bm 0,1) \big) \subset \mathbb R^{r} \times \mathbb R_{\geq 0}.		
		\end{displaymath}
		Then the closed convex hull of $K_0^{\bm a,L}(\underline D)\cup \{(\bm 0,0)\}$, denoted by $K^{\bm a,L}(\underline D)$, is a finitely generated cone in $\mathbb R^{r} \times \mathbb R_{\geq 0}$. Let $R^{\bm a,L}_0(\underline D)$ be the set of extremal rays of $K^{\bm a,L}(\underline D)$ and $R^{\bm a,L}(\underline D) = R^{\bm a,L}_0(\underline D) \setminus \{\mathbb R_{\geq 0}\cdot (\bm 0,1)\}$. Then $\mathrm{lctf}^{\bm a,L}(\underline D)$ is exactly the set of hyperplanes corresponding to rays in $R^{\bm a,L}(\underline D)$.
	\end{theorem}
	
	\subs{\bf Specialization.} For $M = (m_{ij}) \in \mathbb Z^{p\times r}$, we denote $$\sigma_M : \mathbb R^p \to \mathbb R^r, (t_1,\ldots,t_p) \mapsto (\sum_{i=1}^p m_{i1}t_i,\ldots,\sum_{i=1}^p m_{ir}t_i)$$ and denote by $\underline D^M := (\sum_{i=1}^r m_{1i} D_i, \ldots, \sum_{i=1}^r m_{pi} D_i)$ the specialization of $\underline D$.
	\begin{theorem}\label{ThmC}
		With notation as in \ref{Rmk_subs} and above, we have $$\mathrm{Rmk}(\underline D^M) \subset \sigma_M^{-1}(\mathrm{Rmk}(\underline D)) \setminus \{\varnothing,\mathbb R^p\}.$$
	\end{theorem}
%	For $p=1$, Theorem \ref{ThmC} shows that the specializations of
%	the remarkable hyperplanes of $\underline D$ include the negatives
%	of all remarkable numbers of $\underline D^M$.
%	When $r>1$, this inclusion can be strict, showing that remarkable
%	hyperplanes capture information that is not reflected in the
%	remarkable numbers of the specialization\textcolor{red}{, see Example \ref{ex_rmkhyp_no_nonnegative_specialization}.}
	
	\vspace{0.5em}

%	The next result suggests that remarkable hyperplanes can be useful even for a regular map $f=f_1\cdot\ldots\cdot f_r$ to detect combinatorially-interesting divisorial valuations associated to $f$ that are not sources of remarkable numbers of  $f$ by some combinatorial mishap with the log resolution, that is, some $E_j$ such that  $\sum_{i=1}^r N_{ji} \cdot s_i + \nu_j = 0$ is a remarkable hyperplane for $F=(f_1,\ldots,f_r)$, but such that $\nu_j/N_j$ is not a remarkable number of $f$. This was one of the motivations behind this paper.

	\subs{\bf Complete splitting of multibranch plane curve germs.} We determine  the remarkable hyperplanes for complete splittings of multibranch plane curve singularities. 	
	The next result shows that, for multibranch plane curve
	singularities, the remarkable hyperplanes of the complete
	splitting detect numerical data of all rupture components, including those
	that do not give rise to remarkable numbers.
	Overcoming this limitation of remarkable numbers was one
	of the motivations for this work.

	\begin{theorem}\label{ThmD}
		Suppose $r > 1$ and $(C_i,\bm 0),i=1,\ldots,r$ are distinct reduced germs of unibranch curves at $\bm 0 \in \mathbb C^2$. Let $\rho : (X,E) \to (\mathbb C^2,\bm 0)$ be the minimal log resolution of $(C_1\cup \ldots \cup C_r,\bm 0)$. Then the set of local remarkable hyperplanes of $(C_1,\ldots,C_r)$ is
		the set of hyperplanes $\{\bm s\in \mathbb R^r \mid \sum_{i=1}^r N_{ji} s_i + \nu_j = 0\}$ for $E_j$ meeting $\geq 3$ other $E_k$ with $k \in S\setminus \{j\}$, and $\{s_i + 1 = 0\} , i=1,\ldots,r$. 	In particular, all polar hyperplanes of the local motivic zeta function at $\bm 0$ are remarkable.
	\end{theorem}
	In contrast, we recall that $1$ is never a remarkable number for a non-smooth unibranch plane curve singularity, and that although all rupture components of an irreducible plane curve singularity give remarkable numbers, this is not necessarily so for multibranch plane curve singularities, see \cite{BPSZ}. Theorem \ref{ThmD} then says that all rupture components are detectable by remarkable hyperplanes of the complete splitting.
	
\smallskip

	After adjoining a generic line, the remarkable hyperplanes determine the embedded topological type of $(C_1\cup \ldots \cup C_r,\bm 0)$.
	\begin{theorem}\label{thm_topological_type}
		Let $C_0=(\ell=0)$ be a smooth curve germ transverse to
		each $C_j$, $1\leq j\leq r$. Let $f_j$ be a reduced local equation for $C_j$. Then the set
		\[
		\mathrm{Rmk}_{\bm 0}\bigl(
		\operatorname{div}(\ell),
		\operatorname{div}(f_1),\ldots,\operatorname{div}(f_r)
		\bigr),
		\]
		with its coordinates labelled and the coordinate
		corresponding to $\ell$ distinguished, determines the
		labelled embedded topological type of
		$(C_1\cup\cdots\cup C_r,\bm 0)$.
	\end{theorem}
	
	\subs{\bf Case $r=1$ and rational functions.}\label{case_r1} Theorem \ref{ThmA} produces  new results even when $r = 1$. It gives for example a subset of the poles of the motivic zeta function of a rational function. In this case, we adopt the notation from \ref{Rmk_subs}, but replace $\underline D$ with $D$, and $N_{j1}$ with $N_{j}$.

	We call a subset $I \subset S$ \textit{positive (resp. negative) admissible} if $E_I \neq \varnothing$ and there exists $j \in I$ with $N_j > 0$ (resp. $N_j < 0$).      For each stratum $E_I \neq \varnothing$, we define
	\begin{itemize}
		\item $N_I := \mathrm{gcd}_{j\in I}(|N_j|)$,
		\item $\alpha_I := \min_{i\in I,N_i > 0} \{\nu_i/N_i\}$ if $I$ is positive admissible,
		\item $\beta_I := \min_{j\in I,N_j < 0}\{ -\nu_j/N_j\}$ if $I$ is negative admissible.
	\end{itemize}
	We use non-negative gcds, with $\gcd(a,0)=|a|$ and $\gcd(\varnothing)=0$. Thus $N_I>0$ for every admissible $I$. 
Unlike in the case of regular functions, two different tracks of remarkable numbers arise.
	\begin{definition}\label{def_signed_remarkable}
		A rational number $\theta$ is called a \textit{positive (resp. negative) remarkable number} (with respect to $D$) if there exists a positive (resp. negative) admissible subset $I \subset S$ such that $\theta = \alpha_I$ (resp. $\theta = -\beta_I$), and for every positive (resp. negative) admissible subset $J$ with $\alpha_J<\alpha_I$ (resp. $\beta_J<\beta_I$), $N_J$ does not divide $N_I$. Both positive and negative remarkable numbers are called remarkable numbers.
	\end{definition}
	\begin{theorem}\label{ThmE}
		Let $X$ be a smooth irreducible complex variety and $D$ be a divisor on it, not necessarily effective. Then:
		
		(1) The set of remarkable numbers is independent of the choice of log resolution for $D$.
		
		(2) If $\theta$ is a remarkable number, then $-\theta$ is a pole of $\Zmot_{D}(T)$.
	\end{theorem}
	These remarkable numbers can be reinterpreted as:
	 \begin{definition}\label{def_signed_thresholds}
		Let $d \in \mathbb Z$. Fix a log resolution of $D$ as above. Then the positive and negative $d$-log canonical thresholds of $D$ are 
		\begin{align*}
			& \mathrm{lct}^{(d)}_+(D) := \min\{ \alpha_I \mid I \subset S \text{ is positive admissible and } N_I \text{ divides }d\} \text{ and } \\
			& \mathrm{lct}^{(d)}_-(D) := -\min\{ \beta_I \mid I \subset S \text{ is negative admissible and } N_I \text{ divides }d\}.
		\end{align*}
		A threshold is left undefined when its indexing family is empty.
	\end{definition}
	\begin{proposition}\label{prop_signed_thresholds}
		The positive (resp. negative) remarkable numbers of $D$ are precisely the defined values of $\mathrm{lct}^{(d)}_+(D)$ (resp. $\mathrm{lct}^{(d)}_-(D)$) as $d$ ranges through $\mathbb Z_{>0}$. Thus the set of remarkable numbers of $D$ is
		\[
		\bigl\{\mathrm{lct}^{(d)}_+(D),\ \mathrm{lct}^{(d)}_-(D)\mid d\in\mathbb Z_{>0}\bigr\},
		\]
		with undefined values omitted. Allowing all $d\in\mathbb Z$ gives the same set.
	\end{proposition}

	\subs{\bf Status of the monodromy conjecture.} The strong monodromy conjecture for remarkable numbers was recently proved using AI in \cite{BSZ_SMC}.  Appendix \ref{app:smc} here, by subset of the authors of this paper, verifies Theorem \ref{MC_rmk} obtained using AI, giving a proof of the strong monodromy conjecture for remarkable hyperplanes of a tuple of regular functions $F=(f_1,\ldots, f_r)$.
	
	Regarding rational functions, several zeta functions for rational or meromorphic functions occur in the literature. Gusein-Zade, Luengo, and Melle-Hern\'andez \cite{GLM98} studied the Milnor monodromies and their zeta functions at zero and infinity. Raibaut \cite{Rai12,Rai13} constructed motivic zeta functions and motivic Milnor fibers of rational functions at indeterminacy points. Gonz\'alez Villa and Lemahieu \cite{GVL14} introduced a topological zeta function describing the behavior of a meromorphic function as its value tends to zero and proved its monodromy conjecture for plane meromorphic germs; the version at infinity follows by replacing the function with its reciprocal.

	Veys and Z\'u\~niga-Galindo \cite{VZG17} developed local zeta functions and oscillatory integrals for meromorphic functions over characteristic-zero local fields. Their motivic zeta function \cite[Definition 3.17 and Theorem 3.18]{VZG17} records the order of $f/g$ along arcs, allowing positive, zero, and negative orders. For $D=\divi(f)-\divi(g)$, it agrees with the motivic zeta function considered here after applying coefficientwise
	the natural homomorphism from their coefficient ring to
	$\widehat{\mathcal M}_{\mathbb C}$.
	 They also describe its topological specialization, which need not coincide with the topological zeta function of \cite{GVL14}; see \cite[Section 3.5]{VZG17}. 
	
	Takeuchi \cite{Ta23} studied Bernstein--Sato polynomials for
	meromorphic functions and related their roots to Milnor monodromy.
	An alternative definition was considered by \`Alvarez Montaner,
	Gonz\'alez Villa, Le\'on-Cardenal, and N\'u\~nez-Betancourt
	\cite{AGLN25}.

	One could ask for several versions of the monodromy conjecture for rational functions, and therefore for remarkable numbers of rational functions, see \cite[Conjecture~1]{GVL14} for one formulation. Another formulation would be the following. Let $h=f/g$ be a rational function, with $f$ and $g$
	locally coprime, $D=\divi(f)-\divi(g)$, and $x\in D_{\mathrm{red}}$ be a fixed point. The monodromy conjecture predicts that a negative (resp. positive) pole $\theta$ of the motivic zeta function $Z_{D,x}^{\mathrm{mot}}(T)$ gives an eigenvalue $\exp(2\pi i\theta)$ of the Milnor monodromy of $h$ (resp. $1/h$) at some point on $f^{-1}(0)$ (resp. $g^{-1}(0)$) near $x$.  The latter is the monodromy of $h$ at infinity, with local target parameter $1/h$. We did not find counterexamples to these versions.
	
	One would be tempted to ask for stronger versions, for example, whether for some integer $m\geq0$ one has
	\[
	b_{h,x,m}(\theta)=0\quad\text{if }\theta<0,
	\qquad
	b_{h^{-1},x,m}(-\theta)=0\quad\text{if }\theta>0.
	\]
	Here $b_{h,x,m}$ denotes Takeuchi's local analytic Bernstein--Sato polynomial of order $m$ \cite[Definition 1.2]{Ta23}, computed from a coprime local presentation of $h$; the reciprocal is treated in the same way. If correct, these root assertions would imply the respective monodromy assertions by \cite[Theorem 1.4]{Ta23} since the set of monodromy eigenvalues is invariant under complex conjugation. However, with the assistance of AI we constructed a counterexample to these root assertions in the plane curve case, see Example~\ref{counter_SMC_rational}.
The corresponding monodromy assertions remain valid for this example. 
One could ask for analogues of the strong monodromy conjecture for rational functions
by replacing $b_{h,x,m}$ with the $b$-functions $b_{f/g, x}^{\alpha},b_{f/g,x}^{A}$ defined in \cite[Theorem A]{AGLN25}. This is weaker than the version refuted by \cite[Remark 3.16]{AGLN25}. There are however counterexamples to these assertions as well, see \ref{counter_SMC_Alvarez}. %The data and codes supporting these counterexamples are available in \cite{BCKSZ-data}.  
%		Thus it is natural to ask the following question.
%		\begin{question}
%			Do we have
%				\begin{displaymath}
%					b_{h,x}^1(-\theta)=0\quad\text{if }\theta>0,
%					\qquad
%					b_{h^{-1},x}^1(\theta)=0\quad\text{if }\theta<0?
%				\end{displaymath}
%		\end{question}

\subsection*{AI use.}  The work on this article, not including the Appendix and the examples, was largely finished before Summer 2026 without AI. Subsequently AI was used to double-check and correct some of the arguments in the original draft, as well as for editing (e.g. the original pictures were drawn by hand). AI then assisted in the construction of Example~\ref{counter_SMC_rational}, the examples in \ref{counter_SMC_Alvarez}, and the examples in Section~\ref{sec_regular_examples}.
The Appendix was completed in November 2026 and AI was essential for it. M. Čmrlec declined to be a co-author for the Appendix. The authors have verified all AI-generated content.	
	
\subsection*{\bf Acknowledgement.} We thank Nero Budur for proposing the problem to generalize remarkable numbers to hyperplanes, and for useful discussions. Quan Shi gratefully acknowledges the hospitality of KU Leuven during the course of this work. This work was supported by KU Leuven grants Methusalem METH/21/03 and DB/25/013, FWO grant G0B3123N, an``El-yurt umidi" Foundation fellowship, BJNSF Grant 1252009, NSFC Grant 12671056, NSFC Grant 125B2004. 
	
	\section{Preliminary material}

	\subs{\bf Notation.}\label{not_pre} Let $X$ be a smooth irreducible complex variety and $\underline D = (D_1,\ldots,D_r)$ be an ordered collection of divisors on $X$. Let $\rho: Y \to X$ be a log resolution of $\bigcup_{i=1}^r (D_i)_{\mathrm{red}}$ such that $\rho^{-1}(\bigcup_{i=1}^r(D_i)_{\mathrm{red}}) = \bigcup_{i\in S} E_i$ is a simple normal crossings divisor. Write the relative canonical divisor as $K_{\rho} = \sum_{i\in S}(\nu_i -1) E_i$, and write $\rho^*D_i = \sum_{j\in S} N_{ji} E_j$ for $i=1,\ldots,r$. Set $\underline E=(E_j)_{j\in S}$,
	$\bm N_j=(N_{j1},\ldots,N_{jr})$, and
	$\bm\nu=(\nu_j)_{j\in S}$.

	We may also keep extra components in a simple normal crossings divisor, including exceptional components with $\bm N_i=\bm 0$. The same change-of-variables formulas below apply to this larger divisor. This convention is useful when some components cancel after specialization.
	
	\subs{\bf Grothendieck ring.}\label{subsGr} Let $K_0(\operatorname{Var}_{\mathbb{C}})$ be the Grothendieck ring of varieties. It is the quotient of the free abelian group generated by isomorphism classes of complex varieties by relations
	$[X]=[V]+[X\backslash V]$ for all varieties $X$ and closed subvarieties $V\subset X$. Moreover, its multiplication is defined as $[X]\cdot [Y]:=[X\times Y]$ for two varieties. Write $\mathbb{L}=[\mathbb{A}^1]$ for the class of the affine line in $K_0(\operatorname{Var}_{\mathbb{C}})$.

	Let $\mathcal M_{\mathbb C}:=K_0(\operatorname{Var}_{\mathbb C})[\mathbb L^{-1}]$. Its dimension filtration is the decreasing filtration of abelian groups
	\[
	\mathcal F^m\mathcal M_{\mathbb C}:=
	\left\langle [V]\mathbb L^{-j}\mid j\in\mathbb Z,\ \dim V-j\leq-m\right\rangle_{\mathbb Z}.
	\]
	It satisfies $\mathcal F^a\mathcal F^b\subset\mathcal F^{a+b}$. Put
	$\operatorname{vdim}(u):=\inf\{d\in\mathbb Z\mid u\in\mathcal F^{-d}\mathcal M_{\mathbb C}\}$, with $\operatorname{vdim}(0)=-\infty$.
	Then $e^{\operatorname{vdim}}$ is a submultiplicative non-archimedean seminorm. We use the separated completion
	\[
	\widehat{\mathcal M}_{\mathbb C}:=\varprojlim_{m\to\infty}
	\mathcal M_{\mathbb C}/\mathcal F^m\mathcal M_{\mathbb C}.
	\]
	The quotients here are quotients of abelian groups; multiplication on the completion is induced by the multiplicative filtration. We do not assume that $\mathcal M_{\mathbb C}\to\widehat{\mathcal M}_{\mathbb C}$ is injective. The filtration and virtual dimension extend to the completion.

	\subs{\bf Jet space, arc space, and motivic measure.} For $m\in\mathbb N$, we denote the $m$-jet space and arc space of $X$ by $\mathcal{L}_m(X)=\mathrm{Hom}(\operatorname{Spec} \mathbb{C}[t]/(t^{m+1}),X)$ and $\mathcal{L}(X)=\mathrm{Hom}(\operatorname{Spec} \mathbb{C}[[t]],X)$, respectively. They can be endowed with natural scheme structures. Natural projections $\mathbb{C}[[t]]\rightarrow \mathbb{C}[t]/(t^{m+1})$ and $\mathbb{C}[t]/(t^{l+1})\rightarrow \mathbb C[t]/(t^{m+1})$ for $l\geq m$ induce morphisms of schemes 
	\[\pi_m:\mathcal{L}(X)\rightarrow \mathcal{L}_m(X) \qquad \pi^l_m:\mathcal{L}_l(X)\rightarrow \mathcal{L}_m(X)\]
	We write $\mathcal L_\infty(X):=\mathcal L(X)=\varprojlim_m\mathcal L_m(X)$. We say a subset $C\subset \mathcal{L}(X)$ is a cylinder if there exists some $m\in\mathbb{N}$ such that $\pi_m(C)\subset \mathcal{L}_m(X)$ is constructible and $C=\pi_m^{-1}(\pi_m(C))$. A morphism $f : Z \to X$ of varieties induces morphisms of arc/jet spaces, denoted by $f_{k} : \mathcal L_{k}(Z) \to \mathcal L_k(X)$, where $k\in \mathbb Z_{\geq 0}\cup \{\infty\}$.
	
For $n=\dim X$, our normalization of motivic measure is
	\[
	\mu(\pi_m^{-1}(B))=[B]\mathbb L^{-mn}
	\]
	for a constructible set $B\subset\mathcal L_m(X)$. Thus $\mu(\mathcal L(X))=[X]$. A thin set, contained in $\mathcal L(V)$ for a proper closed subvariety $V\subsetneq X$, has measure zero. We also use disjoint countable unions of cylinders and thin sets for which the virtual dimensions of their measures tend to $-\infty$; their measures are the convergent sums of the measures of the pieces. These are the measurable unions needed below, cf. \cite{DL98,CL}.

	\subs{\bf Contact loci.} For an arc $\gamma : \operatorname{Spec}\mathbb C[[t]] \to X$ and an effective divisor $E$, the contact order of $E$ along $\gamma$ is defined to be the order of $t$ in $f\circ \gamma$, where $f$ is a local defining function of $E$ around $\gamma(0) := \pi_0(\gamma)$. Define the $m$-contact locus of $E$ by
	\[\mathscr{X}_m(X,E)=\left\{\gamma\in\mathcal{L}(X) \mid \ord_{\gamma}E=m \right\}\]
	For $m\geq1$, this is the cylinder
	$\mathscr X_m(X,E)=\pi_{m-1}^{-1}(\mathcal L_{m-1}(E))\setminus\pi_m^{-1}(\mathcal L_m(E))$;
	for $m=0$ it is $\pi_0^{-1}(X\setminus E)$. Here the jet spaces of $E$ are viewed inside the corresponding jet spaces of $X$.

	When $E$ is not necessarily effective, let $E_+$ and $E_-$ be effective divisors with no common prime divisors and such that $E = E_+-E_-$. Then for $m\in\mathbb{Z}$, we define the $m$-contact locus of $E$ as
	\[\mathscr{X}_m(X,E)=\bigsqcup_{k\geq \max\{0,-m\}} \mathscr{X}_{m+k}(X,E_+)\cap \mathscr{X}_k(X,E_-).\] 
	This set need not be a cylinder when $E_-\neq 0$, but it is always a countable disjoint union of cylinders.
	
	More generally, given an ordered collection of divisors $\underline{F}=(F_i)_{i\in S}$  and a multi-index $\bm{m} = (m_i)_{i\in S}$, we define the multi-contact locus as
	\[\mathscr X_{\bm m}(X,\underline{F})= \bigcap_{i\in S} \mathscr X_{m_i}(X,F_i).\]

	Recall from \ref{not_pre} that $\underline D$ is an ordered collection of divisors on $X$. Using the proof of \cite[Theorem A]{ELM04}, we have the following structure theorem for $\mathscr X_{\bm m}(X,\underline D)$.
	\begin{theorem}\label{ELM_upgrade}
		With the notation of \ref{not_pre}, for $\bm m \in \mathbb Z^r$, we have
		\begin{displaymath}
			\mathscr X_{\bm m}(X,\underline D) = \bigsqcup_{I\subset S} \bigsqcup_{\bm{u} \in \mathcal Q_{I}} \rho_{\infty}(\mathscr X_{\bm{u},I}(Y,\underline E)) \sqcup Z,
		\end{displaymath}
		where $\mathcal Q_I = \{\bm{u} \in \mathbb N^S \mid \operatorname{supp}(\bm u)=I,\ \sum_{j\in S} u_j\cdot N_{ji} = m_i \text{ for all }i\}$, $\mathscr X_{\bm{u},I}(Y,\underline E) = \{\gamma \in \mathscr X_{\bm{u}}(Y,\underline E) \mid \gamma(0) \in E_I^\circ\}$, and $Z$ is a thin set of $\mathcal L(X)$. Moreover, each $\rho_{\infty}(\mathscr X_{\bm{u},I}(Y,\underline E))$ is a cylinder in $\mathcal L_{\infty}(X)$ with motivic volume $[E_I^\circ] \cdot \mathbb L^{-\bm{u}\cdot \bm \nu} \cdot (\mathbb L-1)^{\vert I\vert}$.
	\end{theorem}
	\begin{proof}
	An arc outside the thin exceptional set has a unique lift to $Y$. Its finite contact orders $u_i$ determine its stratum $E_I^\circ$, with $u_i>0$ exactly for $i\in I$. On each fixed contact-order piece $\mathscr X_{\bm{u},I}(Y,\underline E)$, the Jacobian order is
	constant and equals $
	e=\sum_{i\in I}u_i(\nu_i-1)$.
	Hence its image is a cylinder by
	\cite[Corollary 1.7]{ELM04}.
	The change-of-variables formula multiplies its measure
	by $\mathbb{L}^{-e}$. Now since the measure of $\mathscr X_{\bm{u},I}(Y,\underline E)$ on $Y$ is $[E_I^\circ](\mathbb L-1)^{|I|}\mathbb L^{-\sum_i u_i}$, we conclude the proof.
	\end{proof}
	Since all $\nu_i$ are positive, we obtain:
	\begin{proposition}\label{decay_vol}
		For all $\bm m \in \mathbb Z^r$, the motivic volume $\mu(\mathscr X_{\bm m}(X,\underline D))$ is well-defined. Moreover, we have $\mathrm{vdim}\, \mu(\mathscr X_{\bm m}(X,\underline D) ) \to -\infty$ as $\vert \bm m\vert \to +\infty$.
	\end{proposition}
	\begin{proof}
	The case $S=\varnothing$ is immediate. Otherwise put $C:=\max_i|\bm N_i|$ and $c:=\min_i\nu_i>0$. A piece indexed by $\bm u$ has virtual dimension at most $n-\sum_i\nu_i u_i\leq n-c\sum_i u_i$. There are only finitely many $\bm u\in\mathbb N^S$ with bounded sum. This proves summability for each $\bm m$. If $C>0$, then $|\bm m|\leq C\sum_i u_i$, giving the bound $n-c|\bm m|/C$. If $C=0$, all nonzero contact vectors have empty loci.
	\end{proof}
	Let
	\begin{displaymath}
		\widehat{\mathcal{M}}_{\mathbb{C}}\{\bm{T}^{\pm}\} = \widehat{\mathcal{M}}_{\mathbb{C}}\{T_1^{\pm},\dots,T_r^{\pm}\} := \{ \sum_{\bm{\alpha}\in\mathbb Z^r} a_{\bm{\alpha}}\bm{T}^{\bm{\alpha}} \mid a_{\bm{\alpha}}  \in \widehat{\mathcal{M}}_{\mathbb{C}} \}
	\end{displaymath}
	be the abelian group of bilateral power series. The \textit{motivic zeta function of $\underline D$} is defined as the generating series of motivic volumes of $\mathscr X_{\bm m}(X,\underline D)$ for $\bm m \in \mathbb Z^r$, i.e.
	\begin{displaymath}
		Z_{\underline D}^{\mathrm{mot}}(\bm{T}) := \sum_{\bm{\alpha} \in \mathbb Z^r} \mu(\mathscr X_{\bm{\alpha}}(X,\underline D)) \cdot \bm{T}^{\bm{\alpha}} \in \widehat{\mathcal{M}}_{\mathbb{C}}\{\bm{T}^{\pm}\}.
	\end{displaymath}
	However, some subtleties may occur since $\widehat{\mathcal{M}}_{\mathbb{C}}\{\bm{T}^{\pm}\}$ is not a ring in the natural way. To this end, we will introduce a ring contained in $\widehat{\mathcal{M}}_{\mathbb{C}}\{\bm{T}^{\pm}\}$ and prove that $Z_{\underline D}^{\mathrm{mot}}(\bm{T})$ belongs to it in the next subsection.

	\subs{\bf Ring of convergent bilateral power series.}\label{convergent}
	
	\begin{definition-proposition}
		The ring of convergent bilateral power series over $\widehat{\mathcal{M}}_{\mathbb{C}}$ is
		\[\widehat{\mathcal{M}}_{\mathbb{C}}\{\!\{T,T^{-1}\}\!\}=\left\{\sum_{i\in\mathbb{Z}} c_iT^i \mid c_i \in \widehat{\mathcal M}_{\mathbb C},\  \lim_{|i|\rightarrow \infty} \mathrm{vdim}(c_i)=-\infty  \right\}\]
		equipped with termwise addition and ordinary multiplication.
		\[\left(\sum c_iT^i\right)\cdot \left(\sum d_iT^i\right)=\sum_{k} \left(\sum_i c_id_{k-i} \right) T^{k} \]
	\end{definition-proposition}
	\begin{proof}
	The coefficient dimensions of each convergent series have a common upper bound, say $C$ and $D$. For fixed $k$, the dimensions of $c_i d_{k-i}$ tend to $-\infty$ as $|i|\to\infty$, so the convolution coefficient exists. Given $m$, choose $B$ so that $\operatorname{vdim}(c_i)<-m-D$ and $\operatorname{vdim}(d_i)<-m-C$ for $|i|>B$. If $|k|>2B$, at least one factor of every $c_i d_{k-i}$ is in its indicated tail, and hence every product has dimension below $-m$. This proves closure under multiplication. The same bounds show summability of the double families used in associativity, and justify their rearrangement. The argument applies without change to multi-indices with the $\ell^1$ norm.
	\end{proof}
	Similarly, we can define the ring of multivariable convergent bilateral power series by
	\[\widehat{\mathcal{M}}_{\mathbb{C}}\{\!\{\bm{T}^{\pm}\}\!\}=\left\{\sum_{\bm{\alpha}\in\mathbb{Z}^r} c_{\bm{\alpha}}\bm{T}^{\bm{\alpha}} \mid c_{\bm{\alpha}} \in \widehat{\mathcal M}_{\mathbb C},\  \lim_{|\bm{\alpha}|\rightarrow \infty} \mathrm{vdim}(c_{\bm{\alpha}})=-\infty  \right\},\]
	equipped with the natural multiplication. By Proposition \ref{decay_vol}, we have $Z_{\underline D}^{\mathrm{mot}}(\bm{T}) \in \widehat{\mathcal{M}}_{\mathbb{C}}\{\!\{\bm{T}^{\pm}\}\!\}$.
	
	As in the case of effective divisors, $Z_{\underline D}^{\mathrm{mot}}(\bm{T})$ admits a Denef-Loeser type formula in terms of a log resolution. This result is somewhat well known to experts, but is difficult to find in the literature, so we include a proof for completeness.
	\begin{theorem}
		With notation as in \ref{not_pre}, we have
		\begin{equation}\label{Zmot_DL}
			Z^{\mathrm{mot}}_{\underline{D}}(\bm{T})=\sum_{I \subset S} [E_I^{\circ}] \prod_{j\in I}\frac{(\mathbb{L}-1)\mathbb{L}^{-\nu_j}\bm{T}^{\bm{N}_{j}}}{1-\mathbb{L}^{-\nu_j}\bm{T}^{\bm{N}_{j}}},
		\end{equation}
		where $E_I^\circ = (\bigcap_{i\in I}E_i) \setminus (\bigcup_{j\notin I} E_j)$.
	\end{theorem}
	\begin{proof}
		By Theorem \ref{ELM_upgrade}, we have
		\begin{equation}\label{Zmot_mid}
			Z^{\mathrm{mot}}_{\underline{D}}(\bm{T}) =\sum_{I \subset S} [E_I^{\circ}] \left(\sum_{\bm \alpha\in\mathbb{N}_{>0}^{I}} \prod_{j\in I} (\mathbb{L}-1)\mathbb{L}^{-\alpha_j\nu_j}\bm{T}^{\alpha_j \cdot \bm N_{j}}    \right).
		\end{equation}
		We conclude the proof by expanding the right hand side of \eqref{Zmot_DL} into a convergent bilateral power series and comparing it with the right hand side of \eqref{Zmot_mid}. 
	\end{proof}
	\begin{remark}
		The main difficulty in proving \eqref{Zmot_DL} is that $\widehat{\mathcal M}_{\mathbb C}\{\bm{T}^{\pm}\}$ is not a ring. If $r = 1$, this issue can also be avoided by using constructible motivic integration, see \cite{CL}.  
	\end{remark}

	\subs{\bf Virtual Poincar\'e specialization.} \label{VP_subs}Consider the ring $\mathbb{Z}[w,w^{-1}]$ equipped with the degree function $\mathrm{ord}(\sum_{i\in \mathbb Z} c_iw^i) = \max\{b \mid c_b \neq 0\}$ and $\mathrm{ord}(0) = -\infty$. Its completion for convergence of the degree to $-\infty$ is $R_0 := \mathbb{Z}((w^{-1}))$, the ring of Laurent series in $w^{-1}$.
	
	Recall the virtual Poincaré specialization $\operatorname{VP}:\widehat{\mathcal{M}}_{\mathbb{C}}\rightarrow R_0$ defined by
	\[[X]\mapsto \sum_{m,l} (-1)^{l+m} \dim \operatorname{Gr}^W_m(H^l_c(X,\mathbb{C}))\cdot w^m\]
	where $W$ is the weight filtration on compactly supported cohomology. For any variety $X$ of dimension $d$, $\operatorname{VP}([X])$ is a polynomial in $w$ of degree $2d$ with leading coefficient equal to the number of irreducible components of $X$ of dimension $d$. For $u\in\mathcal M_{\mathbb C}$, one has $
	\mathrm{ord}(\operatorname{VP}(u))
	\leq 2\operatorname{vdim}(u)$.
	Hence $\operatorname{VP}$ extends uniquely and continuously
	to a ring homomorphism $\widehat{\mathcal M}_{\mathbb C} \to R_0$.
	The same bound holds on the completion.

	The $\mathrm{VP}$ specialization can naturally be extended to
	\[\operatorname{VP}:\widehat{\mathcal{M}}_{\mathbb{C}}\{\!\{\bm T^{\pm}\}\!\}\rightarrow R_0\{\!\{\bm{T}^\pm\}\!\}\]
	where 
	\[ R_0\{\!\{\bm{T}^\pm\}\!\}=\left\{\sum_{\bm{\alpha} \in \mathbb{Z}^r } c_{\bm{\alpha}}\bm T^{\bm\alpha} \; \mid c_{\bm{\alpha}} \in R_0, \ \lim_{\vert \bm{\alpha} \vert \rightarrow\infty} \mathrm{ord}\, c_{\bm{\alpha}} = -\infty  \right\}\]
	is a ring equipped with the natural addition and multiplication.

	For technical use, we define the codimension of an element $A \in \widehat{\mathcal M}_{\mathbb C}$ as
	\begin{displaymath}
		\codim A := -\mathrm{ord}(\mathrm{VP}(A))/2 + \dim X.
	\end{displaymath}
	For $\bm m \in \mathbb Z^r$, we define the codimension of the $\bm m$-contact locus as $C_{\bm m} := \codim \mu(\mathscr X_{\bm m}(X,\underline D))$. This generalizes the usual definition for contact loci of regular functions.

	\subs{\bf Polar loci.}\label{subsPol} Even for a single regular function, the notion of a pole of the
	motivic zeta function requires a careful definition;
	see \cite[Section~2.7]{BPSZ}.
	As in \cite{BPSZ}, we need these definitions only to state our
	strongest results.
	Most of our arguments concern the virtual Poincar\'e specialization,
	which is a rational function in $w,T_1,\ldots,T_r$.
	Every polar hyperplane of this specialization is also a polar
	hyperplane of the motivic zeta function.

	%Recall that even for a single regular function, defining correctly a useful and natural notion of poles of the motivic zeta function requires a pedantic approach, see \cite[2.7]{BPSZ}. This is used in \cite{BPSZ} and here only for stating the strongest claims. In practice, our results are, like in \cite{BPSZ}, about the poles of classical rational functions over $\bC$. %So one might skip the pedantic definitions below at a first reading.

	Let $\mathcal D$ be the multiplicative set generated by
	$1-\mathbb L^{-a}\bm T^{\bm b}$, where $a\in\mathbb Z_{>0}$ and $\bm b\in\mathbb Z^r\setminus\{\bm0\}$, in $\widehat{\mathcal M}_{\mathbb C}[\bm T^{\pm}]$. Its image under $\mathrm{VP}$ is generated by $1-w^{-2a}\bm T^{\bm b}$. Factors with $\bm b=\bm 0$ are already units in the completed coefficient rings and produce no hyperplanes.

	Adjoin compatible rational powers, with $\mathrm{VP}(\mathbb L^{1/N})=w^{2/N}$, and set
	\begin{align*}
	\mathcal M[\mathbb L^{\mathbb Q}]&:=\widehat{\mathcal M}_{\mathbb C}[\mathbb L^{1/N}\mid N>0],\\
	\mathcal A&:=\big(\mathbb C\otimes_{\mathbb Z}\mathcal M[\mathbb L^{\mathbb Q}]\big)
	[(1-\xi\mathbb L^q)^{-1}\mid q\in\mathbb Q\setminus\{0\},\ \xi\text{ a root of unity}],\\
	R_{0,\mathbb Q}&:=\bigcup_{N>0}\mathbb Z((w^{-1/N})),\qquad
	A:=\bigcup_{N>0}\mathbb C((w^{-1/N})).
	\end{align*}
	The field $A$ contains the image of the completed coefficients. Let $\mathcal D_{\mathbb Q}$ be generated by $1-\mathbb L^{-q}\bm T^{\bm b}$ for $q>0$ rational and $\bm b\neq\bm 0$, and use the same notation for its images. Put
	\[
	\mathcal R:=\mathcal M[\mathbb L^{\mathbb Q}][\bm T^{\pm}][\mathcal D_{\mathbb Q}^{-1}],\quad
	\mathcal R':=\mathcal A[\bm T^{\pm}][\mathcal D_{\mathbb Q}^{-1}],\quad
	R:=R_{0,\mathbb Q}[\bm T^{\pm}][\mathcal D_{\mathbb Q}^{-1}],\quad
	R':=A[\bm T^{\pm}][\mathcal D_{\mathbb Q}^{-1}].
	\]
	These rings fit into the natural commutative diagram
	\begin{equation}\label{eqComDia}
	\xymatrix@C=12pt{
	\widehat{\mathcal M}_{\mathbb C}\brT\ar[d]^{\mathrm{VP}}&
	\widehat{\mathcal M}_{\mathbb C}[\bm T^{\pm}][\mathcal D^{-1}]\ar@{_{(}->}[l]\ar[r]\ar[d]&\mathcal R\ar[r]\ar[d]&\mathcal R'\ar[d]\\
	R_0\brT&R_0[\bm T^{\pm}][\mathcal D^{-1}]\ar@{_{(}->}[l]\ar@{^{(}->}[r]&R\ar@{^{(}->}[r]&R'.}
	\end{equation}
	The leftward maps use the convergent geometric expansions. All bottom horizontal maps are injective. In particular, no injectivity assertion about the uncompleted motivic ring is involved.

	We start with some technical lemmas to define polar hyperplanes.
%	\begin{notation}
%		For $\bm \beta_1 = (\beta_{11},\ldots,\beta_{1r}),\bm \beta_2 = (\beta_{21},\ldots,\beta_{2r}) \in \mathbb N^2$, we call $\bm \beta_1 \geq \bm \beta_2$ if $\beta_{1i} \geq \beta_{2i}$ for all $i = 1,\ldots,r$.
%	\end{notation}
	\begin{definition}\label{def_P}
		Let $(\bm a,b) = (a_1,\ldots,a_r,b) \in (\mathbb Z^r\setminus \{\bm 0\} )\times \mathbb Q_{>0}$ with ${\mathrm{gcd}}(a_1,\ldots,a_r) = 1$,  $\xi$ be a root of unity, and $P =  \sum_{\bm \alpha \in \mathbb Z^r} p_{\bm \alpha} \cdot \bm T^{\bm \alpha} \in \mathcal A[T_1^{\pm},\ldots,T_r^{\pm}]$ be a polynomial. %We define $\deg P := \max\{\vert \beta_1\vert +\ldots+\vert \beta_r\vert \mid p_{\bm \beta} \neq 0\}$. 
		For $\bm \beta \in \mathbb Z^r$, let $S^P_{\bm \beta} := \sum_{k \in \mathbb Z} \xi^k p_{\bm \beta - k\cdot \bm a} \mathbb L^{-bk}$. Note that $\xi^{-1}\mathbb L^{b} \cdot S_{\bm \beta + \bm a}^P = S_{\bm \beta}^P$.
	\end{definition}
	\begin{lemma}\label{division_lem0}
	With the notation of Definition \ref{def_P}, $P$ is divisible by $1-\xi\mathbb L^{-b}\bm T^{\bm a}$ if and only if $S_{\bm\beta}^P=0$ for every $\bm\beta\in\mathbb Z^r$. The same criterion holds over $(A,w^2)$ and also for a nonprimitive nonzero vector $\bm a$.
	\end{lemma}
	\begin{proof}
	Put $U=\bm T^{\bm a}$ and group the terms of $P$ according to the cosets of $\mathbb Z\bm a$ in $\mathbb Z^r$. For a representative $\bm\beta$, the corresponding group of terms is
	\[
	\sum_{k\in\mathbb Z}p_{\bm\beta-k\bm a}\bm T^{\bm\beta-k\bm a}
	=\bm T^{\bm\beta}q_{\bm\beta}(U),
	\qquad
	q_{\bm\beta}(U):=\sum_{k\in\mathbb Z}p_{\bm\beta-k\bm a}U^{-k}.
	\]
	These sums are finite because $P$ is a Laurent polynomial. Multiplication by $1-\xi\mathbb L^{-b}\bm T^{\bm a}$ only shifts exponents by $\bm a$, so it does not mix different cosets. Thus $P$ is divisible by this binomial if and only if each $q_{\bm\beta}(U)$ is divisible by $1-\xi\mathbb L^{-b}U$. Indeed, there are only finitely many nonzero coset components, so the quotients obtained for the individual components can be added to give a Laurent polynomial.

	Since $\xi^{-1}\mathbb L^b$ is a unit, we have
	\[
	1-\xi\mathbb L^{-b}U\mid q_{\bm\beta}(U)
	\quad\Longleftrightarrow\quad
	q_{\bm\beta}(\xi^{-1}\mathbb L^b)=0.
	\]
	By the definition of $q_{\bm\beta}$,
	\[
	q_{\bm\beta}(\xi^{-1}\mathbb L^b)
	=\sum_{k\in\mathbb Z}\xi^k p_{\bm\beta-k\bm a}\mathbb L^{-bk}
	=S_{\bm\beta}^P.
	\]
	Changing the representative from $\bm\beta$ to $\bm\beta+\ell\bm a$ multiplies this value by the unit $\xi^\ell\mathbb L^{-b\ell}$. Hence vanishing for one representative of each coset is equivalent to vanishing for every $\bm\beta\in\mathbb Z^r$, proving the criterion. The same argument applies over $(A,w^2)$ and does not require $\bm a$ to be primitive.
	\end{proof}

	\begin{lemma}\label{division_lem1}
		Let $(a_1,\ldots,a_r,b),(c_1,\ldots,c_r,d) \in (\mathbb Z^r\setminus \{\bm 0\} )\times \mathbb Q_{>0}$ be such that $\mathrm{gcd}(a_1,\ldots,a_r) = \mathrm{gcd}(c_1,\ldots,c_r) = 1$ and $\xi,\zeta$ be roots of unity. Suppose
		\begin{enumerate}
			\item $P \in \mathcal A[T_1^{\pm},\ldots,T_r^{\pm}]$ is not divisible by $1 - \xi\mathbb L^{-b}T_1^{a_1}\cdots T_r^{a_r}$ in $\mathcal A[T_1^{\pm},\ldots,T_r^{\pm}]$. 
			\item  Either
			\begin{enumerate}
				\item $(c_1,\ldots,c_r)$ is not a multiple of $(a_1,\ldots,a_r)$, or\label{case1}
				\item $(c_1,\ldots,c_r) = \lambda \cdot (a_1,\ldots,a_r)$ but ${(\xi^{\lambda},{\lambda b}) \neq (\zeta,d)}$. %\textcolor{red}{If we assume $\mathrm{gcd}(c_1,\dots,c_r)=1$ the only possibility are $\lambda=\pm 1$. Is this intended?}\textcolor{blue}{Yes, intended. Do not want to split cases.}
				\label{case2}
			\end{enumerate}
		\end{enumerate}
		Then $P \cdot (1 - \zeta\mathbb L^{-d}T_1^{c_1}\cdots T_r^{c_r})$ is not divisible by $1-\xi\mathbb L^{-b}T_1^{a_1}\cdots T_r^{a_r}$.  The same conclusion holds when replacing $(\mathcal A,\mathbb L)$ with $(A,w^2)$.
	\end{lemma}
	
	\begin{proof}
		%		Write $P = \sum_{\bm \alpha \in \mathbb Z^r} p_{\bm \alpha} \cdot \bm T^{\bm \alpha}$. For $\bm \beta \in \mathbb N^r$, let $S_{\bm \beta} := \sum_{k \in \mathbb Z} \xi^k p_{\bm \beta - k\cdot \bm a} \mathbb L^{-bk}$.	
		
		Suppose $P \cdot (1 - \zeta\mathbb L^{-d}T_1^{c_1}\cdots T_r^{c_r})$ is divisible by $1-\xi\mathbb L^{-b}T_1^{a_1}\cdots T_r^{a_r}$. Let $Q = P \cdot (1 - \zeta\mathbb L^{-d}\bm T^{\bm c})$. Then $S^Q_{\bm \beta} = S_{\bm \beta}^P - \zeta\mathbb L^{-d} S_{\bm \beta-\bm c}^P$. By Lemma \ref{division_lem0}, we have $S^Q_{\bm \beta} = 0$ for all $\bm \beta$ and $S^P_{\bm \beta_0} \neq 0$ for some $\bm \beta_0$. Hence, $S_{\bm \beta}^P = \zeta\mathbb L^{-d} S_{\bm \beta-\bm c}^P$ for all $\bm \beta$.

		Since $P$ is a polynomial, there exists a finite subset $Z$ of $\mathbb Z^r / \mathbb Z\cdot \bm a$ such that $S_{\bm \beta}^P \neq 0$ if and only if $[\bm \beta]\in Z$. For case \eqref{case1}, a contradiction appears since $S^P_{\bm \beta_0},S^P_{\bm \beta_0+\bm c},S^P_{\bm \beta_0+2\bm c},\ldots$ are all non-zero but $\bm \beta_0,\bm \beta_0+\bm c,\bm \beta_0+2\bm c,\ldots$ are pairwise distinct in $\mathbb Z^r/\mathbb Z\cdot \bm a$. For case \eqref{case2}, we have $S_{\bm \beta_0}^P = \zeta \mathbb L^{-d} S_{\bm \beta_0-\bm c}^P = \zeta \mathbb L^{-d} S_{\bm \beta_0 - \lambda \bm a}^P = \zeta \mathbb L^{-d} \xi^{-\lambda} \mathbb L^{\lambda b} S_{\bm \beta_0}^P$. However, $1- \zeta \xi^{-\lambda} \mathbb L^{\lambda b-d}$ is invertible in $\mathcal A$ under our assumption. Hence, $S^P_{\bm \beta_0} = 0$, a contradiction. The second assertion follows from the same argument.
	\end{proof}
	\begin{lemma}\label{division_lem2}
	Let $\bm a\in\mathbb Z^r\setminus\{\bm 0\}$, $b\in\mathbb Q_{>0}$, $\lambda\in\mathbb Z_{>0}$, and $P\in\mathcal A[\bm T^{\pm}]$. Put $\zeta_j:=e^{2\pi j\sqrt{-1}/\lambda}$. If $1-\zeta_j\mathbb L^{-b}\bm T^{\bm a}$ divides $P$ for every $j=0,\ldots,\lambda-1$, then $1-\mathbb L^{-\lambda b}\bm T^{\lambda\bm a}$ divides $P$. The same holds over $(A,w^2)$.
	\end{lemma}
	\begin{proof}
	The factors $1-\zeta_jU$, with $U:=\mathbb L^{-b}\bm T^{\bm a}$, generate pairwise comaximal ideals: an appropriate scalar linear combination of any two is a nonzero complex constant. Their product is $1-U^\lambda$, so divisibility by all factors implies divisibility by the product.
	\end{proof}

	\begin{definition}
		We call a hyperplane of the form $H := \{(s_1,\ldots,s_r) \in \mathbb R^r \mid \sum_{i=1}^r a_i s_i + b = 0\} $ with $(a_1,\ldots,a_r,b) \in (\mathbb Z^r\setminus \{\bm 0\} )\times \mathbb Z_{>0}$ \textit{an allowable hyperplane}.  Define the \textit{fundamental tuple of $H$} to be $(a_1/a,\ldots,a_r/a,b/a)$, where $a = \mathrm{gcd}(a_1,\ldots,a_r) > 0$. 
	\end{definition}
	The following are well-defined due to Lemma \ref{division_lem1} and since elements in $\mathcal D_{\mathbb Q}$ are not zero divisors in both $\mathcal A[T_1^{\pm},\ldots,T_r^{\pm}]$ and $A[T_1^{\pm},\ldots,T_r^{\pm}]$.
	
	\begin{defn}\label{defPole} Let $H$ be an allowable hyperplane, $(a_1,\ldots,a_r,b)$ be its fundamental tuple, and $\xi\in \mathbb C$ be a root of unity.
		
		(1) If $0\neq F\in \mathcal A[T_1^{\pm},\ldots,T_r^{\pm}]$, respectively $\in A[T_1^{\pm},\ldots,T_r^{\pm}]$, denote by $n(F(T),H,\xi)$ the unique natural number $n$ such that $(\mathbb L^{-b}T_1^{a_1}\cdots T_r^{a_r}-\xi)^n\mid F(T)$ but $(\mathbb L^{-b}T_1^{a_1}\cdots T_r^{a_r}-\xi)^{n+1}\nmid F(T)$ in $\mathcal A[\bm T^{\pm}]$, respectively $(w^{-2b}T_1^{a_1}\cdots T_r^{a_r}-\xi)^n\mid F(T)$ but $(w^{-2b}T_1^{a_1}\cdots T_r^{a_r}-\xi)^{n+1} \nmid F(T)$.
		
		(2) If $0\neq Z\in \mathcal R'$ (respectively $\in R'$), write $Z=F/G$ with $0\neq F\in \mathcal A[T_1^{\pm},\ldots,T_r^{\pm}]$ (respectively  $\in A[T_1^{\pm},\ldots,T_r^{\pm}]$) and $G\in \cD_\bQ \subset \mathcal A[T_1^{\pm},\ldots,T_r^{\pm}]$ (respectively  $\subset A[T_1^{\pm},\ldots,T_r^{\pm}]$). We call $H$ {\it  a polar hyperplane} of $Z$ if $M_H:=\max_{\xi}(n(G,H,\xi)-n(F,H,\xi))>0$, where the maximum is over the roots of unity $\xi$. If $H$ is a polar hyperplane, we call $M_H$ {\it the order of $H$ as a polar hyperplane.} 
		
		(3) If $0\neq Z \in \widehat{\mathcal M}_{\mathbb C}[T_1^{\pm},\ldots,T_r^{\pm}][\cD^{-1}]$, respectively $\in R_0[T_1^{\pm},\ldots,T_r^{\pm}][\cD^{-1}]$, we say an allowable hyperplane $H$ is {\it a polar hyperplane of $Z$ of order $n$} if it is so
		for  the image of $Z$ in $\cR'$, respectively in $R'$.
	The zero element has no polar hyperplanes, including when a nonzero motivic class maps to zero in $\mathcal R'$. Binomial multiplicities are finite: after extending a primitive exponent vector to a lattice basis, this is ordinary division by a monic linear polynomial in one Laurent variable.
	\end{defn}
	
	We  give next an equivalent formulation of the definition of poles and pole orders. Note that any element $Z$ of $\widehat{\mathcal M}_{\mathbb C}[T_1^{\pm},\ldots,T_r^{\pm}][\cD^{-1}]$ (resp. $R_0[T_1^{\pm},\ldots,T_r^{\pm}][\cD^{-1}]$, $\cR'$, $R'$) can be written as
	$$
	Z = {Q}/{\prod_{i=1}^p (1-{\mathbb L}^{-c_i}\bm T^{\bm d_i})^{e_i}},
	$$
	where $c_i\in\mathbb Q_{>0}$, $e_i \in \mathbb Z_{>0}$, $\bm d_i \in \mathbb Z^r \setminus \{\bm 0\}$, $[\bm d_i,c_i] \in \mathbb P^r_{\mathbb Q}$ are pairwise distinct, and $Q \in \widehat{\mathcal M}_{\mathbb C}[T_1^{\pm},\ldots,T_r^{\pm}]$ (resp. $R_0[T_1^{\pm},\ldots,T_r^{\pm}]$, $\mathcal A[T_1^{\pm},\ldots,T_r^{\pm}]$, $A[T_1^{\pm},\ldots,T_r^{\pm}]$). We will call this a {\it standard expression} for $Z$. In the numerical rings, each $\mathbb L^{-c_i}$ in such an expression means $w^{-2c_i}$. Denominators defining the same hyperplane can be put over a common multiple of their exponent vectors; all remaining factors are absorbed into $Q$.

	\begin{lemma}\label{estimation of poles}
		Let $Z = {Q}/{\prod_{i=1}^p (1-{\mathbb L}^{-c_i}\bm T^{\bm d_i})^{e_i}}$ in $\mathcal R'$ (resp. $R'$) be a standard expression. Set $H_i := \{ (s_1,\ldots,s_r) \in \mathbb R^r \mid \sum_{j=1}^r d_{ij} s_j + c_i = 0 \}$.
		
		(1)	Then the polar hyperplanes of $Z$ are in $\{H_i\}_{i=1}^p$, and each $H_i$ that is polar has order at most $e_i$.

		(2) If $Q$ is not divisible by $1 - \mathbb{L}^{-c_i}\bm T^{\bm d_i}$ in $\mathcal A[T_1^{\pm},\ldots,T_r^{\pm}]$ (resp. $A[T_1^{\pm},\ldots,T_r^{\pm}]$) then $H_i$ is a polar hyperplane of order exactly $e_i$.%Then  $A^{\mathrm{td}}$ is also an element of $R[T][{\prod_{i=1}^r (1-\mathcal L^{-c_i}T^N)^{-e_i}}]$.
		%		
		%		(3) If $Q$ is not divisible by $1 - \mathcal{L}^{-c_i}\bm T^{\bm d_i}$ in $\cA[T_1^{\pm},\ldots,T_r^{\pm}]$ (resp. $A[T_1^{\pm},\ldots,T_r^{\pm}]$) for all $i$, then for any other standard expression $Z = Q'/\prod_{i=1}^{p'} (1- {\mathbb L}^{-c'_i}\bm T^{d_i'})^{e'_i}$  and for any $i=1,\dots, p$, there exists $i$ such that $H_i = H_i' := \{(s_1,\ldots ,s_r)\in \mathbb Q^r \mid \sum_{j=1}^r d_{ij}' s_j + c_i' = 0 \}$ and $e_i \leq e'_j$.
	\end{lemma}
	\begin{proof}
	Let $\delta_i:=\gcd(d_{i1},\ldots,d_{ir})>0$. The $i$-th denominator factors into $1-\zeta\mathbb L^{-c_i/\delta_i}\bm T^{\bm d_i/\delta_i}$, for $\zeta^{\delta_i}=1$. Factors belonging to different $H_i$ cannot affect each other's multiplicities, by Lemma \ref{division_lem1}. This proves (1). If the whole binomial does not divide $Q$, Lemma \ref{division_lem2} supplies one factor which does not divide $Q$, and its denominator multiplicity is $e_i$. This proves (2).
	\end{proof}
		
		%		 Follows from Lemma \ref{division_lem1} and the facorization $1-\mathbb L^{-c_i}\bm T^{\bm d_i} = \prod_{j=1}^{\delta_i}(1- e^{2\pi j\sqrt{-1} /\delta_i} \mathbb L^{-c_i/\delta_i} \bm T^{\bm d_i/\delta_i}$, where $\delta_i = \mathrm{gcd}(d_{i1},\ldots,d_{ir})$.
		%		
		%		(3) Let $\tilde N$ be a common multiple of $N$ and $N'$. Then $Z(T) = \tilde Q(T) / \prod_{i=1}^r(1- {\mathbb L} ^{-\tilde Nc_i/N}T^{\tilde N})^{e_i}$, where $\tilde Q = Q\cdot \prod_{i=1}^r(\sum_{j=0}^{\tilde N/N-1}  {\mathbb L}^{-c_i j}T^{Nj})^{e_i}$. Then $\tilde Q$ is not divisible by each $1-{\mathbb L}^{-\tilde Nc_i/N}T^{\tilde N}$ since $1- {\mathbb L}^{-\tilde Nc_i/N}T^{\tilde N}$ for all $i$ are pairwise coprime in $\cA[T]$ (resp. $A[T]$). So we may assume $N = N'$. Applying coprimality again, one deduces that for each $i=1,...,r$, there exists $j\in \{1,\dots,r'\}$ such that $c_i = c'_{j}$ and $e_i \leq e'_{j}$.

		\begin{definition}
	\label{def_standard_polar_order}For $n\in\mathbb Z_{\geq0}$, say that $Z$ has polar order at most $n$ along $H$ if it has a standard expression in which the factor associated to $H$ is absent, or has exponent at most $n$. It has polar order $n>0$ if it has order at most $n$ but not at most $n-1$. For motivic or numerical zeta functions we apply this definition to their images in $\mathcal R'$ or $R'$.
	\end{definition}

	\begin{corollary}
		Let $Z$ and $H$ be as in Definition \ref{defPole}. Then $H$ is a polar hyperplane of order $n$ of $Z$ as in Definition \ref{defPole} if and only if $Z$ has polar order $n$ along $H$ in the sense of Definition \ref{def_standard_polar_order}.
	\end{corollary}
	\begin{proof}
	We may assume $Z \neq 0$. If the factor for $H$ is absent, both orders are zero. Otherwise let $e$ be its exponent and set
	\[k:=\min_\xi n(Q,H,\xi),\]
	where $\xi$ runs through the phases of that binomial. Pairwise comaximality shows that its $k$-th power divides $Q$. Cancelling it lowers the exponent to $\max\{0,e-k\}$, which is precisely the maximum of the residual phase multiplicities in Definition \ref{defPole}. Lemma \ref{division_lem1} shows that other hyperplanes do not change this number.
	\end{proof}

	\subs{\bf Specialization of polar hyperplanes.} We suppose $r\geq 2$ in this subsection. For $\bm m = (m_1,\ldots,m_r) \in \mathbb Z^r \setminus \{\bm 0\}$, define $\sigma_{\bm m} : \mathbb R \to \mathbb R^r, s \mapsto (m_1s,\ldots,m_rs)$ to be the specialization map and $\tau_{\bm m} : \mathbb Z^r \to \mathbb Z$ to be the right multiplication by $\bm m$, i.e. $\tau_{\bm m}(\bm a):=\bm a\cdot\bm m$. For $Z\in \mathcal R'$ or $R'$, we denote by $Z_{\bm m}$ its specialization along $\bm m$, i.e. $Z_{\bm m}(T):=Z(T^{m_1},\ldots,T^{m_r})$. 
	\begin{lemma}\label{generic_separation}
	Let $\bm e\in\mathbb Z^r$ be primitive and nonzero, and let $\bm a_1,\ldots,\bm a_p$ have distinct images in $\mathbb Z^r/\mathbb Z\bm e$. There is a subset $U\subset\mathbb P^{r-1}(\mathbb Q)$ dense in the Euclidean topology such that, for every integral representative $\bm m$ of a point of $U$, $\bm e\cdot\bm m\neq0$ and the integers $\bm a_i\cdot\bm m$ are distinct modulo $(\bm e\cdot\bm m)\mathbb Z$.
	\end{lemma}
	\begin{proof}
	A unimodular change of coordinates reduces to $\bm e=(1,0,\ldots,0)$. Write $\bm a_i=(\theta_i,\bm a_i')$, so the $\bm a_i'$ are distinct. Choose an integer $B>2+2\sum_i|\bm a_i'|$ and use the vectors
	\[
	\bm m = (B^r u_1,\ B^r u_2+1,\ B^r u_3+B,\ldots,B^r u_r+B^{r-2}),\qquad \bm u=(u_1,\ldots,u_r)\in\mathbb Z^r,\quad u_1\neq0.
	\]
	For distinct $i,j$, the scalar product of $\bm a_i'-\bm a_j'$ with $(1,B,\ldots,B^{r-2})$ is nonzero by uniqueness of balanced base-$B$ expansion, and has absolute value below $B^r$. Thus the corresponding scalar products are distinct modulo $B^r$, hence also modulo $B^r u_1$. The projective directions of these vectors are dense in the projective space. Both conditions in the statement are unchanged if $\bm m$
	is replaced by a nonzero rational multiple that is still integral.
	\end{proof}
	\begin{lemma}\label{non-division}
		Let $H$ and $(\bm a,b)$ be as in Definition \ref{defPole}, and let $\xi$ be a root of unity. If a Laurent polynomial $P\in\mathcal A[\bm T^{\pm}]$ is not divisible by $1-\xi\mathbb L^{-b}\bm T^{\bm a}$, then there exists a real-projectively dense subset $U$ of $\mathbb P_{\mathbb Q}^{r-1}$ such that for all $\bm m \in \mathbb Z^r\setminus \{\bm 0\}$ with $[\bm m] \in U$, $\bm a \cdot \bm m \neq 0$ and $P_{\bm m}$ is not divisible by $1-\xi\mathbb L^{-b}T^{\bm a\cdot \bm m}$.
		
		The same conclusion holds when replacing $(\mathcal A,\mathbb L)$ with $(A,w^2)$.
	\end{lemma}
	\begin{proof}
		Write $P = \sum_{\bm \alpha \in \mathbb Z^r} p_{\bm \alpha}\cdot \bm T^{\bm \alpha}$, set $Z_0 := \{\bm \alpha \mid p_{\bm \alpha} \neq 0\}$, and let $Z_0'$ be the image of $Z_0$ in $\mathbb Z^r/\mathbb Z\cdot \bm a$. By Lemma \ref{division_lem0}, there exists $\bm \beta_0 \in Z_0$ such that $S^P_{\bm \beta_0} \neq 0$. Choose $\bm\beta_1,\ldots,\bm\beta_k\in Z_0$ representing the elements of
		$Z_0'\setminus\{[\bm\beta_0]\}$.
		%Let $\{\bar{\bm \beta}_1,\ldots,\bar{\bm \beta}_k\}$ be the set of elements in $Z_0'$ with $\bm \beta_i - \bm \beta_0 \notin \mathbb Z\cdot \bm a$. 
		Now by Lemma \ref{generic_separation}, there exists a dense subset $U_0$ of $\mathbb P_{\mathbb Q}^{r-1}$ such that for all $\bm m \in \mathbb Z^r \setminus \{\bm 0\}$ with $[\bm m] \in U_0$, the integers $\bm \beta_0\cdot \bm m,\bm \beta_1\cdot \bm m,\ldots,\bm \beta_k \cdot \bm m$ are pairwise distinct {in $\mathbb{Z}/(\bm{a}\cdot \bm{m}) \mathbb{Z}$}. So $S_{\bm \beta_0\cdot \bm m}^{P_{\bm m}} = S_{\bm \beta_0}^{P} \neq 0$. Using Lemma \ref{division_lem0} and noting that $\{[\bm b] \mid \bm b\cdot \bm a = 0\}$ is nowhere dense, we conclude that $P_{\bm m}$ is not divisible by $1-\xi\mathbb L^{-b}T^{\bm a\cdot \bm m}$. The second assertion follows from the same argument.
	\end{proof}
	\begin{corollary}\label{poles_specialization}
		Let $Z$ and $H$ be as in Definition \ref{defPole}. Then $H$ is a polar hyperplane of $Z$ if and only if there exists a dense subset $U$ of $\mathbb P_{\mathbb Q}^{r-1}$ such that for all $\bm m \in \mathbb Z^r\setminus \{\bm 0\}$ with $[\bm m] \in U$, $\bm a \cdot \bm m \neq 0$ and $-b/(\bm a \cdot \bm m)$ is a pole of  $Z_{\bm m}$.
	\end{corollary}
	\begin{proof}
	Suppose $H$ is polar. In a standard expression whose exponent for $H$ is its pole order $e>0$, the corresponding binomial has the form $1-\mathbb L^{-\delta b}\bm T^{\delta\bm a}$. At least one of its factors $1-\xi\mathbb L^{-b}\bm T^{\bm a}$ does not divide the numerator. By Lemma \ref{non-division}, this nondivisibility persists on a real-projectively dense set of integral directions, with the same phase $\xi$. Remove the finitely many proper projective hyperplanes on which $\bm a\cdot\bm m=0$ or a different candidate hyperplane has the same intersection with the specialization line. At least one factor over $-b/(\bm a\cdot\bm m)$ then remains uncancelled. Conversely, when $H$ is not polar it is absent from some standard denominator. Away from the same finite collection of coincidence hyperplanes, its intersection point cannot be a pole of a specialization. This proves the assertion for both coefficient rings.
	\end{proof}

	\section{Remarkable numbers for rational functions}\label{sec_signed_numbers}
	
	In this section, we study the remarkable numbers for a single divisor, not necessarily effective.
	\subs{\bf Top-degree specialization.}
	\begin{definition}
		Recall that $R_0 := \mathbb Z(\!(w^{-1} )\!)$ is the ring of Laurent series in $w^{-1}$. 
		For $0 \neq g\in R_0$, suppose $\mathrm{ord}(g) = u$ and let $c \cdot w^u$ be the leading term of $g$, i.e. $g = c w^{u} + c_1 w^{u-1}+ c_2w^{u-2}+\cdots$. Then we define the \textit{top degree part of $g$} to be $g^{\mathrm{td}} := cw^u$. Further set $0^\mathrm{td} = 0$.
		
		Define the \textit{non-negative semiring} $R_{0,+} \subset R_0$ by
		\begin{displaymath}
			R_{0,+} = \{\rho \in R_0 \mid \text{the coefficient of } \rho^{\mathrm{td}} \text{ is non-negative}\}.
		\end{displaymath}
	\end{definition}

	Recall 
	\begin{displaymath}
		R_0 \{\!\{ T^{\pm} \}\!\} := \{\sum_{i\in \mathbb Z} a_i T^i \mid a_i \in R_0, \lim_{\vert i\vert \to +\infty} \mathrm{ord}(a_i) = -\infty \}.
	\end{displaymath}
	Then $R_0[T^{\pm}]$ can be embedded into $R_0\{\!\{ T^\pm \}\!\}$ in a natural way.
	
	For a series $Z = \sum_{i\in \mathbb Z} a_iT^i$, denote by
	\begin{displaymath}
		Z^{\mathrm{td}} := \sum_{i \in \mathbb Z} (a_i)^{\mathrm{td}} \cdot T^i
	\end{displaymath}
	the \textit{top degree subseries of $Z$.} This definition is consistent with \cite{BPSZ}.
	%	The \textit{td} arises from the top degree zeta function. Although here \textit{lt} may be more precise,  the convention is kept in this paper.
	
	Recall that the factors $1-w^{-a}T^b$, with $a \in \mathbb Z_{>0}$ and $b\in\mathbb Z$, are invertible in $R_0\{\!\{T^{\pm}\}\!\}$ by their convergent geometric expansions. Consequently the Laurent-polynomial embedding extends to their localization.

	\begin{notation}
		For two series $Z_1,Z_2 \in R_0\{\!\{ T^{\pm} \}\!\}$, we write $Z_1 =_* Z_2$ if $Z_1-Z_2 \in R_0[T^{\pm}]$. And we call elements of $R_0[T^{\pm}]$ Laurent polynomials. 
	\end{notation}
	\begin{notation}
		For a nonzero Laurent polynomial $P \in R_0[T^{\pm}]$, denote by $\deg_{T}(P)$ and $\deg_{T^{-1}}(P)$ its degrees in $T$ and $T^{-1}$, respectively, and set $\mathrm{span}(P) = \deg_{T}(P) + \deg_{T^{-1}}(P)$. For example, $\mathrm{span}(w^{-1}T^{-2} + w^3T^5) = 2+5 = 7$.
	\end{notation}
	
\begin{definition}
	Let $Q = \sum_{i} q_iT^i \in R_0[T,T^{-1}]$ be a Laurent polynomial and fix $0\neq c,d\in \bZ$. We define the Laurent polynomial $Q_{c,d}^1 = \sum_i \bar q_i T^i \in R_0[T,T^{-1}]$ by setting 
	
	$$\bar q_i := \begin{cases}
		0 & \text{ if } q_i = 0;\\
		0 & \text{ if } q_i\neq 0\text{ and }\exists\, j\neq i\text{ such that }q_j\neq 0,\, d\text{ divides }i-j,\\
		& \text{ and } \mathrm{ord}(q_i)-\mathrm{ord}(q_j)<\left(\lfloor \frac{i}{d}\rfloor-\lfloor \frac{j}{d}\rfloor\right)c;\\
		q_i^\td & \text{ otherwise}.
	\end{cases}$$
	
	% \textcolor{brown}{Is it $-\mathrm{ord}(q_{i})+\mathrm{ord}(q_{j})$ or $\mathrm{ord}(q_{i})-\mathrm{ord}(q_{j})$}
	
	We define inductively the polynomial $Q_{c,d}^r := (Q-\sum_{j=1}^{r-1} Q_{c,d}^j)^1_{c,d}$ for $r>1$. In particular, we have $[-\deg_{T^{-1}}(Q_{c,d}^r),\deg_T(Q_{c,d}^r)] \subset [-\deg_{T^{-1}}(Q),\deg_T(Q)]$ for all $r\geq 1$ such that $Q_{c,d}^r \neq 0$. Moreover, $Q$ can be written as the sum $\sum_{r=1}^{\infty} Q_{c,d}^r$, which at the coefficient of each power of $T$ converges in the $w^{-1}$-adic sense.
\end{definition}
The following lemma can be easily checked.
\begin{lemma}\label{divide_infinite}
	Let $P\in R_0[T^{\pm}]$, $0\neq c,d\in\mathbb Z$, and $h\geq1$. The binomial power $(1-w^cT^d)^h$ divides $P$ if and only if its quotient, expanded in convergent bilateral series, has finite support. If $B(P):=\max\{|i| \mid [T^i]P\neq0\}$, with $B(0):=0$, a Laurent-polynomial quotient is supported in $|i|\leq B(P)+h|d|$. If every layer $P_{c,d}^j$ is divisible by this binomial power, then $P$ is divisible by it as well.
\end{lemma}
\begin{proof}
	The inverse is the usual geometric expansion when $c<0$; when $c>0$, factor out $-w^cT^d$ and expand in $w^{-c}T^{-d}$. Finite support is equivalent to polynomial division. Comparison of the smallest and largest exponents in a product with $(1-w^cT^d)^h$ gives the displayed support bound. All the layer quotients therefore have a common finite support. Their coefficients converge: polynomial division on this fixed support uses only finitely many additions and multiplications by fixed powers of $w$. Passing to their sum proves the final assertion.
\end{proof}

\begin{lemma}\label{key_lemma}
	Let $N,a_i,b_j,u_i,v_j$ be positive integers with $a_1<\cdots<a_m$ and $b_1<\cdots<b_n$, and let $P\in R_0[T^{\pm}]$. Empty products are allowed. For
	\[
	Z:=\frac{P}{\prod_{i=1}^m(1-w^{-a_i}T^N)^{u_i}\prod_{j=1}^n(1-w^{-b_j}T^{-N})^{v_j}},
	\]
	one has
	\begin{equation}\label{eq1}
	Z^{\mathrm{td}}=_*
	\left(\frac{P}{\prod_i(1-w^{-a_i}T^N)^{u_i}}\right)^{\mathrm{td}}+
	\left(\frac{P}{\prod_j(1-w^{-b_j}T^{-N})^{v_j}}\right)^{\mathrm{td}}.
	\end{equation}
	Moreover, suppose $P=P_0(T^N)$ and neither $1-w^{-a_1}T^N$ nor $1-w^{-b_1}T^{-N}$ divides $P$. Then
	\begin{equation}\label{eq5}
	Z^{\mathrm{td}}=_*
	\frac{P_{-a_1,N}^{j_{0+}}}{(1-w^{-a_1}T^N)^{u_1}}+
	\frac{P_{-b_1,-N}^{j_{0-}}}{(1-w^{-b_1}T^{-N})^{v_1}},
	\end{equation}
	where $j_{0+}$ and $j_{0-}$ are the first layers not divisible by the corresponding full denominator powers. If one product is empty, the corresponding summand in \eqref{eq5} and its associated assumption are omitted.
\end{lemma}
\begin{proof}
	The assertion is immediate for $P=0$. Otherwise decompose $P=\sum_{e=0}^{N-1}T^eP_e(T^N)$; the residue classes have disjoint support and may be treated separately. Thus it suffices to treat the case $P=P_0(T^N)$.
	Put $U=T^N$ and regard $P$ as a Laurent polynomial in $U$. Since distinct denominator factors of $Z$ are comaximal, we can write $Z$ as
	\[
	Z=L(U)
	+\sum_{i=1}^{m}\sum_{\ell=1}^{u_i}
	\frac{A_{i,\ell}}{(1-w^{-a_i}U)^\ell}
	+\sum_{j=1}^{n}\sum_{\ell=1}^{v_j}
	\frac{B_{j,\ell}}{(1-w^{-b_j}U^{-1})^\ell},
	\]
	where $L(U)\in R_0[U,U^{-1}]$ and
	$A_{i,\ell},B_{j,\ell}\in R_0$.

	For sufficiently large $k>0$, the decomposition shows the coefficient of $U^k$ in $Z$ equals
	\[
	\sum_i w^{-a_i k}
	\sum_{\ell=1}^{u_i}A_{i,\ell}
	\binom{k+\ell-1}{\ell-1}.
	\]
	For each $i$ for which some $A_{i,\ell}$ is nonzero, put
	$h_i=\max_\ell\operatorname{ord}(A_{i,\ell})$.
	The contribution of this $i$ has leading term
	$w^{h_i-a_i k}p_i(k)$, where $p_i$ is a nonzero polynomial,
	since the binomial polynomials have distinct degrees.
	Thus $p_i(k)\neq0$ for all sufficiently large $k$,
	and the lower terms have degree at most
	$h_i-1-a_i k$, uniformly in $k$.
	
	If all $A_{i,\ell}$ vanish, then
	$P/\prod_i(1-w^{-a_i}U)^{u_i}$ is a Laurent polynomial.
	Both this polynomial and $Z$ have zero coefficients of $U^k$
	for all sufficiently large $k>0$.
	Otherwise, let $i_0$ be the smallest index for which some
	$A_{i_0,\ell}$ is nonzero.
	Since $a_1<\cdots<a_m$, the contribution of $i_0$ determines
	the top-degree part of the coefficient of $U^k$ in $Z$
	for all sufficiently large $k>0$.
	
	In this latter case, cancel the full denominator powers
	at the exponents $a_i<a_{i_0}$ and denote the resulting
	numerator by $Q$.
	Put $a=a_{i_0}$, $u=u_{i_0}$ and $V=w^{-a}U$.
	Every other remaining denominator factor has the form
	$1-w^{-d}V^{\pm1}$, with $d>0$, and
	\[
	(1-w^{-d}V^{\pm1})^{-1}
	=1+\sum_{q\geq1}w^{-dq}V^{\pm q}.
	\]
	The remainder of each $V^{\pm q}$ modulo $(1-V)^u$
	has integer coefficients. Thus all terms after $1$ still have
	negative $w$-degree after taking remainders.
	Multiplying by these inverses therefore leaves the highest
	nonzero $w$-degree part of the remainder of $Q(w^aV)$ modulo $(1-V)^u$ unchanged.
	This part comes from the first layer of $Q$ not divisible by
	$(1-w^{-a}U)^u$.
	
	Removing the negative-direction denominator factors therefore
	leaves the top-degree part of the coefficient of $U^k$ in $Z$
	unchanged for all sufficiently large $k>0$.
	Applying the same argument to the coefficient of $U^{-k}$
	in $Z$ proves \eqref{eq1}.
	Under the additional assumptions of \eqref{eq5},
	the selected exponents are $a_1,b_1$ for the nonempty products,
	and no smaller powers are cancelled.
	Thus $Q$ is the original numerator, giving precisely
	the layers stated in \eqref{eq5}.
\end{proof}

% In particular \eqref{eq5} implies that taking $\mathrm{td}$ will not add new poles:
\begin{corollary}\label{td_not_add_poles}
	Let $N,a_1,\ldots,a_m,b_1,\ldots,b_n,u_1,\cdots,u_m,v_1,\ldots,v_n$ be positive integers with $a_1 < \ldots < a_m$ and $b_1 < \ldots < b_n$. Consider the series
	\begin{displaymath}
		Z = \frac{P}{\prod_{i=1}^m (1-w^{-a_i} T^{N})^{u_i} \prod_{j=1}^n (1-w^{-b_j} T^{-N})^{v_j}},
	\end{displaymath} 
	where $P \in R_0[T,T^{-1}]$. Then
	\begin{displaymath}
		Z^{\mathrm{td}} \in R_0[T,T^{-1}] \cdot \left(\prod_{i=1}^m (1-w^{-a_i} T^{N})^{u_i} \prod_{j=1}^n (1-w^{-b_j} T^{-N})^{v_j}\right)^{-1}.
	\end{displaymath}
\end{corollary}
\begin{proof}
	Split $P$ into its residue classes modulo $N$, and in each summand first cancel every common denominator factor. Apply Lemma \ref{key_lemma} to the remaining positive and negative tails, omitting empty tails. Their denominators divide the original denominator, proving the statement. Applying this to a standard expression with minimal exponents also shows that taking top degree creates no new polar values and cannot increase their orders. To use a common $N$, replace each denominator by a multiple exponent and absorb the finite geometric numerator; this does not change its hyperplane or its multiplicity.
\end{proof}

\subs{\bf Top-dimensional zeta function.} Let $X$ be an $n$-dimensional smooth irreducible complex variety and $D$ be a divisor on $X$, not necessarily effective. Let $\rho : Y \to X$ be a log resolution  of $D$, with $K_\rho = \sum_{j\in S} (\nu_j-1) E_j$ and $\rho^*D = \sum_{j\in S} N_j E_j$. The motivic zeta function of $D$ is computed by
\begin{displaymath}
	Z_D^{\mathrm{mot}}(T) = \sum_{I \subset S} [E_I^\circ] \cdot \prod_{j\in I} \frac{(\mathbb L-1)\mathbb L^{-\nu_j}T^{N_j}}{1-\mathbb L^{-\nu_j}T^{N_j}}.
\end{displaymath} 
Applying the virtual-Poincar\'{e} specialization, we have
\begin{align*}
	(\mathrm{VP}(\Zmot_D(T)))^\td & = \bigg(\sum_{I\subset S} (c_Iw^{2n}+\text{lower order terms in } w)\prod_{j\in I}\frac{(w^2)^{-\nu_j}T^{N_j}}{1-(w^2)^{-\nu_j}T^{N_j}}\bigg )^\mathrm{td} \\
	& = \bigg(\sum_{I\subset S} c_I\cdot w^{2n} \cdot \prod_{j\in I}\frac{(w^2)^{-\nu_j}T^{N_j}}{1-(w^2)^{-\nu_j}T^{N_j}}\bigg )^\mathrm{td}
\end{align*}
where $n:=\dim X$ and $c_I$ is the number of irreducible components of the smooth, pure-dimensional stratum $E_I^\circ$, with $c_I=0$ for an empty stratum. Positivity of these leading coefficients justifies the second equality. Set $Z_D^{\mathrm{td}}:=(\mathrm{VP}(Z_D^{\mathrm{mot}}))^{\mathrm{td}}$. Denote by $Z_0 := \sum_{I\subset S} c_I\cdot w^{2n} \cdot \prod_{j\in I}\frac{w^{-2\nu_j}T^{N_j}}{1-w^{-2\nu_j}T^{N_j}}$  the series inside $(-)^{\mathrm{td}}$ in the second line.

By Corollary \ref{td_not_add_poles}, poles of $(\mathrm{VP}(\Zmot_D(T)))^\td$ are always poles of $\mathrm{VP}(\Zmot_D(T))$. In this subsection, $Z_0^{\mathrm{td}}$ and its poles will be computed, thus a subset of poles of $\mathrm{VP}(\Zmot_D(T))$ will be given. We combine the framework of \cite{BPSZ} with Lemma \ref{key_lemma}.

%	
%	Algorithm 1
%	
%	Algorithm 2

%	\begin{lemma}
	%		Let $N,a_1,\ldots,a_m,b_1,\ldots,b_n,u_1,\cdots,u_m,v_1,\ldots,v_n$ be positive integers with $a_1 < \ldots < a_m$ and $b_1 < \ldots < b_n$.
	%	\end{lemma}

\begin{lemma}\label{block1}
	Let $\varnothing\neq I\subset S$ with $E_I\neq\varnothing$, and put
	$Z_I:=c_Iw^{2n}\prod_{i\in I}w^{-2\nu_i}T^{N_i}/(1-w^{-2\nu_i}T^{N_i})$.
	Choose $N>0$ divisible by all nonzero $|N_i|$. Then
	\begin{equation}\label{td_block_1}
	Z_I^{\mathrm{td}}=_*
	\frac{P_+}{(1-w^{-2N\alpha_I}T^N)^{m_I}}+
	\frac{P_-}{(1-w^{-2N\beta_I}T^{-N})^{n_I}},
	\end{equation}
	where $m_I:=\#\{i\in I \mid N_i>0,\ \nu_i/N_i=\alpha_I\}$ and $n_I:=\#\{i\in I \mid N_i<0,\ -\nu_i/N_i=\beta_I\}$. A term is omitted when its admissible sign is absent. Each occurring numerator belongs to $R_{0,+}[T^{\pm}]$ and its exponent residues modulo $N$ are exactly $\{0,N_I,\ldots,N-N_I\}$. If all $N_i$ are zero, $Z_I^{\mathrm{td}}$ is constant and both terms are absent.
\end{lemma}
\begin{proof}
	Replace a nonzero factor with exponent $N_i$ by a denominator with exponent $\operatorname{sign}(N_i)N$, using the finite geometric numerator. Factors with $N_i=0$ are scalar units whose top degree is $w^{-2\nu_i}$. The resulting top-degree numerator is
	\[
	Q:=c_Iw^{2n}\prod_{i\in I,\,N_i=0}w^{-2\nu_i}
	\prod_{i\in I,\,N_i\neq0}\left(\sum_{j=1}^{N/|N_i|}(w^{-2\nu_i}T^{N_i})^j\right).
	\]
	The residue classes of exponents of $T$ in $Q$ form the subgroup generated by the $N_i$ in $\mathbb Z/N\mathbb Z$, namely the multiples of $N_I$. Every nonzero coefficient of $Q$ is a Laurent polynomial in $w$ with positive coefficients. For each nonzero residue summand of $Q$, the first layer cannot vanish on substituting $T^N=w^{2N\alpha_I}$, or $T^{-N}=w^{2N\beta_I}$, respectively; hence the relevant binomial does not divide that layer. Lemma \ref{key_lemma} gives \eqref{td_block_1}, with $P_+=Q^1_{-2N\alpha_I,N}$ and $P_-=Q^1_{-2N\beta_I,-N}$. The scalar units omitted above have leading part $1$ after their displayed powers of $w$ are extracted, so do not affect top degree.
\end{proof}

The following are bilateral versions of the coefficient comparisons used in \cite[Lemmas 3.9 and 4.5]{BPSZ}.
\begin{lemma}\label{block2}
	For $1\leq i\leq p$ and $1\leq j\leq q$, let
	\[
	A_i:=\frac{b_iT^{u_i+l_iN}}{(1-w^{-\theta_i}T^N)^{m_i}},\qquad
	B_j:=\frac{c_jT^{v_j+z_jN}}{(1-w^{-\sigma_j}T^{-N})^{n_j}},
	\]
	where $N,m_i,n_j,\theta_i,\sigma_j>0$ are integers, $l_i,z_j\in\mathbb Z$, $0\leq u_i,v_j<N$, and $b_i,c_j$ are positive integer multiples of powers of $w$. Put
	\[
	\alpha_i:=\theta_i/N,\quad\beta_j:=\sigma_j/N,\quad
	\lambda_i:=\operatorname{ord}(b_i)+l_i\theta_i,\quad
	\rho_j:=\operatorname{ord}(c_j)-z_j\sigma_j.
	\]
	For a residue $u$, let $\mathcal S_u:=\{i \mid u_i=u\}$ and $\mathcal T_u:=\{j \mid v_j=u\}$. Within $\mathcal S_u$, retain the smallest $\alpha_i$, and among those retain the largest $\lambda_i$; denote the resulting set by $\mathcal S_u''$. Within $\mathcal T_u$, retain the smallest $\beta_j$, and among those retain the largest $\rho_j$; denote the result by $\mathcal T_u''$. An empty set contributes nothing. Then
	\[
	\left(\sum_i A_i+\sum_jB_j\right)^{\mathrm{td}}=_*
	\sum_{u=0}^{N-1}\left(\sum_{i\in\mathcal S_u''}A_i+\sum_{j\in\mathcal T_u''}B_j\right).
	\]
\end{lemma}
\begin{proof}
	At exponent $u+kN$ with $k\gg0$, the degree of the coefficient of $A_i$ is $\lambda_i-k\theta_i$. For $k\ll0$, the corresponding degree for $B_j$ is $\rho_j+k\sigma_j$. The binomial coefficients multiplying these monomials are positive. Therefore the smallest slope, and then the largest intercept, determines the top degree in either tail. Terms tied in both remain and cannot cancel.
\end{proof}
\begin{lemma}\label{block3}
	Let
	\begin{equation}\label{eqZforPoles}
	Z=W+\sum_{i=1}^p\sum_{j=1}^{m_i}\frac{P_{i,j}}{(1-w^{-a_i}T^N)^j}
	+\sum_{i=1}^q\sum_{j=1}^{n_i}\frac{Q_{i,j}}{(1-w^{-b_i}T^{-N})^j},
	\end{equation}
	where $N,m_i,n_i,a_i,b_i\in\mathbb Z_{>0}$, the $a_i$ are pairwise distinct and the $b_i$ are pairwise distinct, and all numerators and $W$ are Laurent polynomials over $R_0$. Suppose $P_{i,m_i}$ and $Q_{i,n_i}$ are nonzero and belong to $R_{0,+}[T^{\pm}]$. Then the poles are precisely $-a_i/(2N)$ of order $m_i$ and $b_i/(2N)$ of order $n_i$. Empty sums are allowed.
\end{lemma}
\begin{proof}
	Work in the coefficient field $A$. At $T=w^{a_i/N}$, the numerator $P_{i,m_i}$ does not vanish: each nonzero summand has positive leading coefficient, so the finite sum has nonzero leading coefficient. The binomial $1-w^{-a_i}T^N$ has a simple zero there; every denominator of a different slope is nonzero. Terms of smaller denominator multiplicity cannot cancel its highest-order pole. Similarly evaluate the negative tail at $T=w^{-b_i/N}$. These witnesses give the stated orders, and the displayed denominator gives the reverse bound.
\end{proof}

\begin{proof}[Proof of Proposition \ref{prop_signed_thresholds}]
Suppose first that $\alpha_I$ is a positive remarkable number. By Definition \ref{def_signed_remarkable}, every positive admissible $J$ with $N_J\mid N_I$ satisfies $\alpha_J\geq\alpha_I$. Since $I$ itself occurs in this minimum and $N_I>0$, we obtain
\[
\mathrm{lct}^{(N_I)}_+(D)=\alpha_I.
\]
Conversely, let $d\in\mathbb Z$ and suppose $\mathrm{lct}^{(d)}_+(D)$ is defined. Choose a positive admissible $I$ with $N_I\mid d$ and $\alpha_I=\mathrm{lct}^{(d)}_+(D)$. If a positive admissible $J$ satisfied $\alpha_J<\alpha_I$ and $N_J\mid N_I$, then $N_J\mid d$, contradicting the defining minimum. Thus $\alpha_I$ is remarkable.

For the negative case, apply the same argument to the positive quantities $\beta_I$: the defining condition for $-\beta_I$ is equivalent to
\[
\beta_I=\min\{\beta_J \mid J\text{ negative admissible},\ N_J\mid N_I\}.
\]
Hence $\mathrm{lct}^{(N_I)}_-(D)=-\beta_I$. Conversely, any defined value $\mathrm{lct}^{(d)}_-(D)=-\beta_I$ is remarkable by the same divisibility contradiction. Every remarkable number has therefore been obtained with $d=N_I>0$, while the argument for the reverse inclusion allowed every integer $d$. This also proves the last assertion, including $d=0$.
\end{proof}

\subs{\bf Proof of Theorem {\ref{ThmE}}.}
Choose a positive integer $N$ divisible by every nonzero $|N_i|$. By Lemmas \ref{block1} and \ref{block2}, each nonempty positive residue tail of $Z_D^{\mathrm{td}}=Z_0^{\mathrm{td}}$ has slope
\[
2\min\{\alpha_I \mid N_I\mid d,\ I\text{ positive admissible}\}=2\mathrm{lct}^{(d)}_+(D),
\]
and each nonempty negative residue tail has slope $-2\mathrm{lct}^{(d)}_-(D)$, for $0\leq d<N$. The retained numerators have positive monomial coefficients. Grouping residues with the same slope gives an expression of the form \eqref{eqZforPoles}, with positive nonzero highest-order numerators. Lemma \ref{block3} therefore identifies its poles with the negatives of the attained signed thresholds.

By Proposition \ref{prop_signed_thresholds}, these thresholds are exactly the remarkable numbers. The set of poles of the intrinsic series $Z_D^{\mathrm{td}}$ is resolution independent. Corollary \ref{td_not_add_poles} puts these poles among those of $\mathrm{VP}(Z_D^{\mathrm{mot}})$, and hence among the motivic poles by specialization. This proves both assertions. \hfill$\Box$

\begin{remark}\label{rmk_local}
	For a closed subset $\Sigma\subset X$, local remarkable numbers and hyperplanes are defined by retaining strata with $\rho(E_I)\cap\Sigma\neq\varnothing$. They agree with the values on sufficiently small neighborhoods of $\Sigma$. The same pole assertions hold for the local zeta function, formed with arcs centered in $\Sigma$: in the resolution formula its coefficient is $[E_I^\circ\cap\rho^{-1}(\Sigma)]$ instead of $[E_I^\circ]$. Empty exact strata are omitted. At a point of a stratum meeting the fiber, the exact stratum has a larger index set; the latter has a larger lattice and smaller defining polyhedron, so this omission leaves the intersections defining the local lct polyhedra unchanged. The leading coefficients of all the remaining virtual Poincar\'e polynomials are positive. Thus the proofs below of the facet-pole statement apply verbatim; the coefficient degrees may vary, but do not affect the smallest-slope argument.
\end{remark}

The following lemma is reserved for later use.
\begin{lemma}\label{key_lem_pos}
	Let
	\[
	Z=\sum_{i=1}^q\frac{P_i(T)}{\prod_{j=1}^{n_i}(1-w^{-2c_{ij}}T^{d_{ij}})},
	\]
	where $0\neq P_i\in R_{0,+}[T^{\pm}]$, $c_{ij}\in\mathbb Z_{>0}$, and $d_{ij}\in\mathbb Z$. Constant factors with $d_{ij}=0$ are permitted. If positive $d_{ij}$ occur, then
	$-\min_{d_{ij}>0}c_{ij}/d_{ij}$ is a pole. If negative $d_{ij}$ occur, then $\min_{d_{ij}<0}(-c_{ij}/d_{ij})$ is a pole.
\end{lemma}
\begin{proof}
	Let $a:=\min_{d_{ij}>0}c_{ij}/d_{ij}$ and work over $A$ at $T=w^{2a}$. A denominator of this slope vanishes simply. Every other denominator evaluates to $1-w^{-h}$ for some $h>0$, so has positive leading part $1$. After extracting the common local parameter $1-w^{-2a}T$, a vanishing factor contributes the positive integer $d_{ij}$ to its remaining unit. Choose the largest number of vanishing factors in a summand. Its highest-order coefficient is a sum of quotients with positive leading coefficient, since $P_i(w^{2a})$ has positive leading coefficient. It cannot vanish, proving the first assertion. Replacing $T$ by $T^{-1}$ proves the second.
\end{proof}

%\textcolor{blue}{ 
	
%	In the regular case, i.e., when $g=1$, we have the following useful bound on the remarkable numbers.}
%
%\begin{proposition}
%All remarkable numbers of $f$ lie in (0,1].
%\end{proposition}
%    
%\begin{proof}
%Since $C_{m}=\mathrm{codim}\, \mathscr {X}_m(f)\le m$, it follows from \cite[Theorem 1.7]{BPSZ} that every remarkable number of $f$ is at most $1$. Since remarkable numbers are positive, the claim follows.    
%\end{proof}
%
%\textcolor{blue}{We expect that the preceding bound admits an analogue for rational functions. However, the argument used in the proof above does not directly extend to this  setting. We therefore leave the following as an open problem.}
%
%\begin{question}
% Do all remarkable numbers of a rational function $f/g$ lie in the interval [-1,1]?  
%\end{question}
%

%\textcolor{blue}{
	Let $f/g$ be a rational function on $X$, $D=\mathrm{div}(f)-\mathrm{div}(g)$, and $x \in D_{\mathrm{red}}$. In the regular case, it is known that the negative of the smallest remarkable number is a root of the Bernstein-Sato polynomial. Indeed, the log canonical threshold $\mathrm{lct}(f)$ is the smallest remarkable number, and $-\mathrm{lct}(f)$ is the largest root of $b_{f}(s)$, see \cite{Yan83,Lic89,Kol97}. Takeuchi's result \cite[Corollary 5.4]{Ta23} gives the following
	consequence for remarkable numbers of rational functions. Before proceeding, we note that multiplier ideals can also be defined for rational functions, together with their corresponding jumping numbers, analogously to the regular case, cf. \cite[Section~5]{Ta23} and \cite[Sections~5 and~6]{AGLN25}. The authors of the second paper use a different definition of jumping numbers; however, the jumping numbers considered here coincide with their 0-jumping numbers.
	%}

\begin{proposition}\label{prop_root of b-function}
Let $\theta$ be the smallest (resp. greatest) positive (resp. negative) local remarkable number at $x$. If $\theta$ exists, then $-\theta$ (resp. $\theta$) is the largest root of $b_{f/g,x,0}$ and $b^{A}_{f/g,x}$ (resp. $b_{g/f,x,0}$ and $b_{g/f,x}^A$).  
Here $b_{f/g,x,0}$ and $b_{f/g,x}^{A}$ are the local $b$-functions of $f/g$ defined in \cite[Definition~1.2, with $m=0$]{Ta23} and \cite[Definition~3.12]{AGLN25}, respectively.
\end{proposition}
\begin{proof}
Suppose the smallest positive local remarkable number at $x$ exists. By definition, we have $\theta=\min_{\substack{i\in S,\ N_i>0\\x\in\rho(E_i)}} \{\nu_i/N_i\}$. On the other hand, $\min_{\substack{i\in S,\ N_i>0\\x\in\rho(E_i)}} \{\nu_i/N_i\}$ is a minimal jumping number of $f/g$ at $x$, then we have $-\theta$ is the largest root of $b_{f/g,x,0}$ by \cite[Corollary 5.4]{Ta23}. Since $b_{f/g,x,0}$ divides $b_{f/g,x}^{A}$, cf. \cite[\S 1]{Ta23}, we have $b_{f/g,x}^{A}(-\theta) = 0$. Then $-\theta$ is the largest root by \cite[Corollary 3.20]{AGLN25}.

If $\theta$ is the greatest negative local remarkable number for $f/g$ at $x$, then, by symmetry, $-\theta$ is the smallest positive local remarkable number for $g/f$ at $x$, and the second claim follows from the first.   
\end{proof}

	\section{Remarkable hyperplanes}
	
	Recall that $X$ is a smooth irreducible complex variety and $\underline D = (D_1,\ldots,D_r)$ is an ordered collection of divisors on $X$. $\rho : Y \to X$ is a log resolution of $\bigcup_{i=1}^r (D_i)_{\mathrm{red}}$, with
	\begin{displaymath}
		\rho^*D_i = \sum_{j\in S} N_{ji} E_j,\ K_\rho = \sum_{j\in S} (\nu_j-1)E_j,
	\end{displaymath}
	where $\bigcup_{j\in S} E_j$ is SNC. Note that $N_{ji}$ can be negative.
	
	Recall that in the introduction we defined subgroups $N_I$ and polytopes $F_I$ for a non-empty stratum $E_I$. And for $\bm a \in \mathbb Z^r$ and a sublattice $L$ of $\mathbb Z^r$, we defined the $(\bm a,L)$-lct polytope $\mathrm{lct}^{\bm a,L}(\underline D)$ and $\mathrm{lctf}^{\bm a,L}(\underline D)$, the set of hyperplanes spanned by facets of the closure of $\mathrm{lct}^{\bm a,L}(\underline D)$. The set of remarkable hyperplanes, denoted by $\mathrm{Rmk}(\underline D)$, is the union of $\mathrm{lctf}^{\bm a,L}(\underline D)$ as $(\bm a,L)$ run through all points and sublattices in $\mathbb Z^r$.
	
%	We start with several remarks on remarkable hyperplanes.
	\begin{remark}
		\begin{enumerate}
			\item Since $\mathrm{lct}^{\bm a,L}(\underline D) = \mathrm{lct}^{\bm a+\bm b,L}(\underline D)$ for all $\bm b \in L$, one can replace $\bm a\in \mathbb Z^r$ with $\bm a \in \mathbb Z^r/L$ in the definition of $\mathrm{lct}^{\bm a,L}$ and $\mathrm{lctf}^{\bm a,L}$.
			\item If $r = 1$, let $a = a_1, D = D_1$, then $\mathrm{lct}^{a,L}(D)$ equals the interval 
			\begin{displaymath}
				(b_{-},b_+) = \bigcap_{a' \in a+L} (-\mathrm{lct}^{(a')}_{+}(D), -\mathrm{lct}_{-}^{(a')}(D)),
			\end{displaymath}
			where we set $\mathrm{lct}_{\pm}^{(a')} = \pm \infty$ if they do not exist. Due to the periodicity of $\mathrm{lct}_{\pm}^{(a')}$, this is a finite intersection and so elements in $\mathrm{Rmk}(D)$ are the negatives of remarkable numbers of $D$. When $L = N \cdot \mathbb Z$ for some $N \in \mathbb Z_{>0}$ divisible enough, we have $\mathrm{lct}^{a,L}(D) = (-\mathrm{lct}^{(a)}_{+}(D), -\mathrm{lct}_{-}^{(a)}(D))$ and thus $\mathrm{Rmk}(D)$ is exactly the set of negatives of remarkable numbers of $D$. 
%			\item Equivalently, we can define  $\mathrm{lctf}^{\bm a,L}(\underline D)$ the hyperplanes spanned by facets of $\overline{\mathrm{lct}^{\bm a,L}(\underline D)} = \bigcap_{\bm a\in N_I+L,E_I\neq \varnothing} \overline{F_I}$. This equality holds since $\in \mathrm{lct}^{\bm a,L}(\underline D)$ contains an open neighborhood of the origin and thus non-empty.
		\end{enumerate}
	\end{remark}
	
	We start with a simple fact that will be used multiple times.
	\begin{lemma}\label{finite_intersection}
	Let $\mathfrak C_1,\ldots,\mathfrak C_m\subset\mathbb R^r$ be open half spaces containing $\bm 0$, or equal to $\mathbb R^r$. Put $\mathrm C_I:=\bigcap_{i\in I}\mathfrak C_i$, and denote by $\mathrm{Cf}_I$ the supporting hyperplanes of the facets of $\overline{\mathrm C_I}$. If $\bigcap_{j=1}^k\mathrm C_{I_j}=\mathrm C_J$, then $\mathrm{Cf}_J\subset\bigcup_{j=1}^k\mathrm{Cf}_{I_j}$.
\end{lemma}
\begin{proof}
	Every $\mathrm C_I$ contains a neighborhood of $\bm 0$ and is the intersection of the open half spaces determined by its facets. Intersect these facet descriptions for $I_1,\ldots,I_k$. Removing redundant inequalities can only remove, and cannot create, facet hyperplanes. This proves the assertion.
\end{proof}

	The first goal of this section is to prove the following theorem, which implies Theorem~\ref{ThmA}\eqref{A2}\eqref{A3}. Note that Theorem \ref{ThmA}\eqref{A1} follows directly from Theorem \ref{ThmB}.
	\begin{theorem}\label{thm_rmkhyp}
		\begin{enumerate}
			\item Remarkable hyperplanes are polar hyperplanes of the virtual Poincar\'e specialization of the motivic zeta function of $\underline D$.
			\item For $N \in \mathbb Z_{>0}$ divisible enough, we have
			\begin{displaymath}
				\mathrm{Rmk}(\underline D) = \bigcup_{\bm a\in (\mathbb Z/N\cdot \mathbb Z)^r } \mathrm{lctf}^{\bm a,(N\cdot \mathbb Z)^r}(\underline D) = \bigcup_{\bm a\in \mathbb Z^r } \mathrm{lctf}^{\bm a,(N\cdot \mathbb Z)^r}(\underline D).
			\end{displaymath} 
		\end{enumerate}
	\end{theorem}
%	\begin{remark}
%		At the moment I do not know if $\mathrm{lct}^{\bm a,L}(F)$ enjoys a periodicity similar to $\mathrm{lct}^{(d)}(f)$.
%	\end{remark}
	
	%	\subs{\bf Sketch of the proof.}
	
	\begin{definition}
		Let $P = \sum_{\bm b \in \mathbb Z^r} \lambda_{\bm b} \cdot {\bm T}^{\bm b} \in R_0\brT$. For a sublattice $L \subset \mathbb Z^r$   and $\bm a \in \mathbb Z^r$, we define
		\begin{displaymath}
			P^{\bm a,L} := \sum_{\bm b \in \bm a + L} \lambda_{\bm b} \cdot {\bm T}^{\bm b}.
		\end{displaymath}
	\end{definition}
	The following lemma is straightforward.
	\begin{lemma}\label{separation_lemma}
	Let $L\subset\mathbb Z^r$ have full rank and suppose
	$P=Q\prod_k(1-w^{-v_k}\bm T^{\bm b_k})^{-1}$, with $Q\in R_0[\bm T^{\pm}]$, $\bm b_k\in L$, and $v_k\in\mathbb Z_{>0}$. Then
	\[
	P^{\bm a,L}=Q^{\bm a,L}\prod_k(1-w^{-v_k}\bm T^{\bm b_k})^{-1}.
	\]
	In particular, residue projection creates no new polar hyperplanes.
\end{lemma}
\begin{proof}
	Every denominator expansion is supported in $L$, so coefficient selection gives the formula. To prove the final assertion for actual poles, start with a standard expression whose denominator omits a given nonpolar hyperplane. Replace each exponent vector by a positive multiple in $L$, absorbing the finite geometric sum into the numerator. This introduces no new hyperplane. Applying the formula shows that the projected denominator still omits the given hyperplane.
\end{proof}
	%	\begin{lemma}\label{lct_hyperplane}
		%		Let $r,q,n_1,\ldots,n_q \in \mathbb Z_{>0}$, $\bm N_{ij} = (N_{ij1},\ldots,N_{ijr}) \in \mathbb N^r \setminus \{\bm 0\}$, and $c_{ij} \in \mathbb Z_{>0}$, where $i = 1,\ldots,q,j = 1,\ldots,n_i$. For $i=1,\ldots,q$, suppose $P_i = Q_i/ \prod_{j=1}^{n_i}(1- w^{-2c_{ij}}\bm T^{\bm N_{ij}})$ is such that $Q_i \in R_{0,+} \setminus \{0\}$, then hyperplanes spanned by facets of $\mathrm{lct}(P) := \bigcap_{i,j} \{(s_1,\ldots,s_r) \in \mathbb R^r \mid \sum_{k=1}^r N_{ijk} s_k + c_{ij} \geq 0\}$ are polar hyperplanes of $P = \sum_{i=1}^q P_i$.
		%	\end{lemma}
	%	\begin{proof}
		%	
		%		
		%		
		%	\end{proof}
	
	Now we deal with the virtual-Poincar\'e specialization of the motivic zeta function, i.e.
	\begin{displaymath}
		P := \mathrm{VP}(Z_{\underline D}^{\mathrm{mot}}) = \sum_{I \subset S} \mathrm{VP}(E_I^\circ)\cdot (w^2-1)^{\vert I\vert} \cdot \prod_{j\in I} \frac{w^{-2\nu_j} \bm T^{\bm N_j}}{1-w^{-2\nu_j} \bm T^{\bm N_j}},
	\end{displaymath}
	where we abbreviate $T_1^{N_{j1}} \cdots T_r^{N_{jr}}$ as $\bm T^{\bm N_j}$. For $I\subset S$, denote the corresponding summand by $P_I$. 
	\begin{remark}
		Remember, to consider polar hyperplanes, one should view $T_i$ as $w^{-2s_i}$.
	\end{remark}

	\subs{\bf Proof of Theorem \ref{thm_rmkhyp}(1).} Fix $\bm a \in \mathbb Z^r$ and a sublattice $L\subset \mathbb Z^r$. Then there exists a positive integer $K$ such that $K\cdot \bm N_j \in L$ for all $j$. For each $j\in S$ and $I \subset S$, let $Q_j := \sum_{l=1}^{K} w^{-2l\nu_j} \bm T^{l\cdot \bm N_j}$ and $Q_I := \prod_{j\in I} Q_j$. Then we have
	\begin{displaymath}
		P_I = \mathrm{VP}(E_I^\circ)\cdot (w^2-1)^{\vert I\vert} \cdot  \frac{Q_I}{\prod_{j\in I}(1-w^{-2K\nu_j} \bm T^{K\cdot \bm N_j})}
	\end{displaymath}
	and thus
	\begin{displaymath}
		P_{I}^{\bm a,L} = \mathrm{VP}(E_I^\circ)\cdot (w^2-1)^{\vert I\vert} \cdot  \frac{Q_I^{\bm a,L}}{\prod_{j\in I}(1-w^{-2K\nu_j} \bm T^{K\cdot \bm N_j})}
	\end{displaymath} 
	Note that $Q_I^{\bm a,L} \neq 0$ if and only if $(\bm a + L) \cap N_I \neq \varnothing \iff \bm a\in N_I + L$.
	
	%	 Since $P^{\bm a,L} = \sum_{I\subset S} P_I^{\bm a,L}$, we conclude Theorem \ref{thm_rmkhyp} by Lemma \ref{separation_lemma} and Lemma \ref{lct_hyperplane}. \hfill$\Box$
	%	
	
	Suppose $H=\{\bm s \mid \bm\lambda\cdot\bm s+\delta=0\}$ supports a facet of $\overline{\mathrm{lct}^{\bm a,L}(\underline D)}$, with $\delta>0$. Choose a rational point $\bm s_0$ in its relative interior, outside all other distinct candidate hyperplanes, and also outside $\{\bm s \mid \bm N_j\cdot\bm s=0\}$ for every nonzero $\bm N_j$. These excluded sets cut proper affine hyperplanes in $H$, since $\bm 0\notin H$. For a sufficiently divisible positive integer $q$, put $\bm m:=-q\bm s_0\in\mathbb Z^r$.

	Then $\bm\lambda\cdot\bm m=q\delta>0$. For every selected stratum $I$ and $j\in I$ with $\bm N_j\cdot\bm m>0$, one has
	\[
	\frac{\nu_j}{\bm N_j\cdot\bm m}\geq\frac1q
	=\frac{\delta}{\bm\lambda\cdot\bm m},
	\]
	and equality holds precisely for components whose candidate hyperplane is $H$. Under $\bm T^{\bm b}\mapsto T^{\bm b\cdot\bm m}$ the coefficients of each nonzero $Q_I^{\bm a,L}$ still have positive leading coefficient. The same is true after multiplication by $\mathrm{VP}(E_I^\circ)(w^2-1)^{|I|}$. Factors with $\bm N_j=\bm 0$ are harmless scalar units. Lemma \ref{key_lem_pos} therefore gives a pole at $-1/q$ of $(P^{\bm a,L})_{\bm m}$. By the choice of $\bm s_0$, no candidate hyperplane other than $H$ specializes to this point. Thus $H$ must be polar for $P^{\bm a,L}$, and Lemma \ref{separation_lemma} makes it polar for $P$. The same argument works for local coefficients as in Remark \ref{rmk_local}. \hfill$\Box$

	\subs{\bf Periodicity.} The finite polar set of $P$ gives finitely many witnesses $(\bm a_i,L_i)$ such that $\mathrm{Rmk}(\underline D) = \bigcup_{i=1}^p \mathrm{lctf}^{\bm a_i,L_i}(\underline D)$ by Theorem \ref{thm_rmkhyp}(1).
	\begin{lemma}\label{refinement}
		Suppose $L \subset L' \subset \mathbb Z^r$ are sublattices. Write $L'/L = \{[\bm b_1], \ldots,[\bm b_q]\}$, where $\bm b_1,\ldots,\bm b_q \in L'$. Then
		\begin{equation}\label{translate_lct}
			\mathrm{lct}^{\bm a,L'}(\underline D) = \bigcap_{i=1}^q \mathrm{lct}^{\bm a+\bm b_i,L}(\underline D).
		\end{equation} 
		In particular, $\mathrm{lctf}^{\bm a,L'}(\underline D) \subset \bigcup_{i=1}^q \mathrm{lctf}^{\bm a+\bm b_i,L}(\underline D)$.
	\end{lemma}
	\begin{proof}
	A stratum $I$ satisfies $\bm a\in N_I+L'$ if and only if $\bm a+\bm b_i\in N_I+L$ for at least one coset representative $\bm b_i$ of $L'/L$. Hence the collection of inequalities defining the left side of \eqref{translate_lct} is exactly the union of the collections defining its right side. Their intersections are equal. The facet inclusion follows from Lemma \ref{finite_intersection}.
	\end{proof}

	\subs{\bf Proof of Theorem \ref{thm_rmkhyp}(2).}
	Choose $N>0$ with $L_N:=(N\mathbb Z)^r\subset L_i$ for all the finitely many witnesses above. Lemma \ref{refinement} gives
	\[
	\begin{split}
	\mathrm{Rmk}(\underline D)
	&=\bigcup_{i=1}^p\mathrm{lctf}^{\bm a_i,L_i}(\underline D)\\
	&\subset\bigcup_{i=1}^p\ \bigcup_{\bm b\in L_i/L_N}\mathrm{lctf}^{\bm a_i+\bm b,L_N}(\underline D)
	\subset\bigcup_{\bm a\in\mathbb Z^r/L_N}\mathrm{lctf}^{\bm a,L_N}(\underline D)
	\subset\mathrm{Rmk}(\underline D).
	\end{split}
	\]
	All inclusions are therefore equalities. The same holds for every multiple of $N$. \hfill$\Box$

	\subs{\bf Specializations.} Recall that for $M = (m_{ij}) \in \mathbb Z^{p\times r}$, we denote $$\sigma_M : \mathbb R^p \to \mathbb R^r, (t_1,\ldots,t_p) \mapsto (\sum_{i=1}^p m_{i1}t_i,\ldots,\sum_{i=1}^p m_{ir}t_i)$$ and denote by $\underline D^M := (\sum_{i=1}^r m_{1i} D_i, \ldots, \sum_{i=1}^r m_{pi} D_i)$ the specialization of $\underline D$. Let $\tau_M : \mathbb Z^r \to \mathbb Z^p$ denote left multiplication by $M$. For the sake of distinction, we use $F_{I, \underline D},F_{I,\underline D^M}$ and $N_{I, \underline D},N_{I,\underline D^M}$ to denote the corresponding polytopes and sublattices.
%	\begin{lemma}\label{linear_algebra}
%		Let $\varphi : \mathbb Z^r \to \mathbb Z^p$ be a group homomorphism, $\bm a' \in \mathrm{im}\, \varphi$, $B$ be a subgroup $\mathbb Z^r$, $L' \subset \mathrm{im}\, \varphi$ be a sublattice of $\mathbb Z^p$. If $\bm a' \in \varphi(B) + L'$, then there exist $\bm a\in \mathbb Z^r$  such that $\varphi(\bm a) = \bm a'$ and $\bm a \in B + L$, where $L = \varphi^{-1}(L')$.  
%	\end{lemma}
%	\begin{proof}
%		Since $\ker \varphi \subset L$, we have $$\varnothing \neq \varphi^{-1}(\bm a') \subset \varphi^{-1}(\varphi(B)+L') = \varphi^{-1}(\varphi(B+L)) = B+L + \ker \varphi = B+L.$$
%	\end{proof}

	\begin{lemma}\label{specialization_lct}
		For $\bm a'\in\mathbb Z^p$ and a full-rank sublattice
		$L'\subset\mathbb Z^p$, put $L=\tau_M^{-1}(L')$.
		Then $L$ has full rank. If
		$\bm a'\in\operatorname{im}\tau_M+L'$, choose
		$\bm a\in\mathbb Z^r$ such that
		$\bm a'-\tau_M(\bm a)\in L'$. We have
		\begin{equation}\label{specialization_eq}
			\mathrm{lct}^{\bm a',L'}(\underline D^M)
			=
			\begin{cases}
				\sigma_M^{-1}\bigl(
				\mathrm{lct}^{\bm a,L}(\underline D)\bigr),
				& \bm a'\in\operatorname{im}\tau_M+L',\\
				\mathbb R^p, & \text{otherwise}.
			\end{cases}
		\end{equation}
		In the first case,
		\[
		\mathrm{lctf}^{\bm a',L'}(\underline D^M)
		\subset
		\sigma_M^{-1}\bigl(
		\mathrm{lctf}^{\bm a,L}(\underline D)\bigr)
		\setminus\{\varnothing,\mathbb R^p\}.
		\]
		In the second case,
		$\mathrm{lctf}^{\bm a',L'}(\underline D^M)=\varnothing$.
		Here $\sigma_M^{-1}(\mathcal T)
		=\{\sigma_M^{-1}(H)\mid H\in\mathcal T\}$.
	\end{lemma}
	\begin{proof}
		Choose $N>0$ such that $N\mathbb Z^p\subset L'$.
		Then $N\mathbb Z^r\subset L$, so $L$ has full rank.
		On a common log resolution, we have
		\[
		N_{I,\underline D^M}=\tau_M(N_{I,\underline D}),
		\qquad
		F_{I,\underline D^M}=\sigma_M^{-1}(F_{I,\underline D}).
		\]
		In the first case, since $\bm a'-\tau_M(\bm a)\in L'$, we have
		\[
		\bm a'\in N_{I,\underline D^M}+L'
		\quad\Longleftrightarrow\quad
		\bm a\in N_{I,\underline D}+L.
		\]
		Taking intersections over the corresponding strata proves the first case of \eqref{specialization_eq}.
		If $\bm a'\notin\operatorname{im}\tau_M+L'$, the
		indexing family is empty, giving $\mathbb R^p$.
		
		For the assertion about facets, describe
		$\mathrm{lct}^{\bm a,L}(\underline D)$ by its facet
		half spaces and pull them back along $\sigma_M$.
		Every nontrivial pulled-back inequality still contains
		a neighborhood of $\bm0$. Deleting redundant inequalities
		proves the claimed inclusion.
		%Extra components with $\bm N_i=\bm0$ may be retained
		%on the common resolution, by the convention in \ref{not_pre}.
	\end{proof}
	
	\subs{\bf Proof of Theorem \ref{ThmC}.}
	By definition and Lemma \ref{specialization_lct},
	\begin{align*}
		\mathrm{Rmk}(\underline D^M)
		=
		\bigcup_{\bm a',L'}
		\mathrm{lctf}^{\bm a',L'}(\underline D^M) \subset
		\{\sigma_M^{-1}(H)\mid
		H\in\mathrm{Rmk}(\underline D),\
		\sigma_M^{-1}(H)\neq\varnothing\},
	\end{align*}
	where $\bm a'$ runs through $\mathbb Z^p$ and $L'$
	runs through the full-rank sublattices of $\mathbb Z^p$.
	Every nonempty preimage here is a hyperplane, because
	$\bm0\notin H$. This proves the assertion.
	\hfill$\Box$

	\subs{\bf Intrinsic interpretation via contact loci.}

		\begin{lemma}\label{generator_of_cone}
	Let $\bm c,\bm c_1,\ldots,\bm c_m\in\mathbb Z^r$ and $d,d_1,\ldots,d_m>0$. Put $F:=\{\bm s \mid \bm c\cdot\bm s+d>0\}$ and $F_i:=\{\bm s \mid \bm c_i\cdot\bm s+d_i>0\}$. Then $\bigcap_iF_i\subset F$ if and only if $(\bm c,d)$ belongs to the cone generated by $(\bm c_i,d_i)$ and $(\bm0,1)$. Zero vectors $\bm c_i$ and an empty family are allowed.
\end{lemma}
\begin{proof}
	If $(\bm c,d)$ is not in this closed polyhedral cone, separation supplies $(\bm v,t)$ nonnegative on all its generators and negative on $(\bm c,d)$. Replacing $(\bm v,t)$ by $(\bm v,t+\varepsilon)$ for sufficiently small $\varepsilon>0$ makes $t>0$ and all generator pairings strictly positive while preserving the last strict inequality. Dividing by $t$ gives a point of $\bigcap_iF_i$ outside $F$. Conversely, a nonnegative linear combination of the generator inequalities is strictly positive on $\bigcap_iF_i$, since $d>0$ excludes the zero combination.
\end{proof}
	
	The above lemma  says that the dual cone of a polyhedral cone $C$ is polyhedral cone bounded by hyperplanes associated to extremal rays.	

	\begin{lemma}\label{include_by_cones}
	Let $\bm c_1,\ldots,\bm c_m\in\mathbb Z^r$ and $d_1,\ldots,d_m>0$. Fix $\bm a\in\mathbb Z^r$ and a full-rank lattice $L$. For $\bm p\in\bm a+L$, set
	$S_{\bm p}:=\{\bm b\in\mathbb N^m \mid \sum_i b_i\bm c_i=\bm p\}$,
	$\mathcal T_{\bm a}:=\{\bm p \mid S_{\bm p}\neq\varnothing\}$, and suppose $\mathcal T_{\bm a}\neq\varnothing$. Put $\bm d=(d_1,\ldots,d_m)$. Set $\lambda_{\bm p}:=\min_{\bm b\in S_{\bm p}}\bm b\cdot\bm d$ and
	\[
	K_{\bm a}:=\bigcup_{\bm p\in\mathcal T_{\bm a}}
	\big((\bm p,\lambda_{\bm p})+\mathbb R_{\geq0}(\bm0,1)\big).
	\]
	Then the closed convex hull of $K_{\bm a}\cup\{\bm 0\}$ is the cone generated by $(\bm c_i,d_i)$ and $(\bm0,1)$. The same holds with all coordinates of $\bm b$ required to be positive. It also holds for $m=0$ when $\mathcal T_{\bm a}\neq\varnothing$.
\end{lemma}
\begin{proof}
	Denote the closed convex hull by $\mathcal C_1$ and the stated cone by $\mathcal C_0$. The inclusion $\mathcal C_1\subset\mathcal C_0$ is immediate. Choose one admissible $\bm b_0$, and put $(\bm p_0,\lambda_0):=\sum_i b_{0i}(\bm c_i,d_i)$. Choose $N>0$ with $N\mathbb Z^r\subset L$. For each $i$ and $l\geq1$, the point $(\bm p_0,\lambda_0)+lN(\bm c_i,d_i)$ lies in $K_{\bm a}$, since its height is at least the minimum at its first coordinate. For any fixed $t\geq0$, its multiple by $t/(lN)$ is in the convex hull for large $l$, and converges to $t(\bm c_i,d_i)$. Similarly, use $(\bm p_0,\lambda_0+l)$ to obtain $t(\bm0,1)$. Thus $\mathcal C_1$ contains every generating ray, not just one point on each ray. Its convexity then contains their full conical hull. The argument preserves positivity of the $b_i$, and the vertical-ray argument covers $m=0$.
\end{proof}

	\subs{\bf Proof of Theorem \ref{intrinsic_cl}.}
	By Theorem \ref{ELM_upgrade} and positivity of the leading virtual Poincar\'e coefficients, the codimension at contact vector $\bm p$ is
	\[
	C_{\bm p}=\min\left\{\sum_{i\in I}b_i\nu_i \mid
	E_I^\circ\neq\varnothing,\ b_i\in\mathbb Z_{>0},\ \sum_{i\in I}b_i\bm N_i=\bm p\right\}.
	\]
	For $I=\varnothing$ the sum is zero and only $\bm p=\bm 0$ is allowed. A minimum exists whenever the locus is nonempty, since $\nu_i>0$. For each fixed $I$, admissible positive tuples in the coset $\bm a+L$ exist exactly when $\bm a\in N_I+L$: an integral solution modulo $L$ can be made positive by adding sufficiently large multiples of the exponent of $\mathbb Z^r/L$ to its coordinates.

	Decompose $K_0^{\bm a,L}$ into the corresponding stratumwise epigraphs. By Lemma \ref{include_by_cones}, the closed convex hull of each nonempty epigraph together with $\bm 0$ is the cone generated by $(\bm0,1)$ and $(\bm N_i,\nu_i)$, $i\in I$. The closed convex hull of their finite union is therefore the cone generated by $(\bm0,1)$ and
	\[
	\{(\bm N_i,\nu_i) \mid i\in S_{\bm a,L}\},\qquad
	S_{\bm a,L}:=\{i \mid \exists I\ni i,\ E_I\neq\varnothing,\ \bm a\in N_I+L\}.
	\]
	If there are no admissible strata, the cone is instead $\{\bm 0\}$ and the lct polyhedron is $\mathbb R^r$. In the other cases, Lemma \ref{generator_of_cone} identifies its nonvertical extremal rays with the irredundant inequalities defining $\mathrm{lct}^{\bm a,L}$. These are exactly its facet hyperplanes.  \hfill$\Box$

\vspace{0.5em}

	We conclude this section by specializing the preceding result to the case of a rational function, that is, $r=1$ and $D=\mathrm{div}(f)-\mathrm{div}(g)$. The proof is the same as in \cite[Theorem 1.7]{BPSZ}.

\begin{theorem}
 Let $C_{m}$ denote the codimension of $\mathscr{X}_m(X,D)$ in $\mathcal L(X)$ and set $N=\mathrm{lcm}(|N_i|\mid i\in S,\ N_i\neq0)$, with $\mathrm{lcm}(\varnothing)=1$. Fix $1\leq d\leq N$. Then:
 \begin{enumerate}
               \item $\mathrm{lct}^{(d)}_{+}(D)$ (resp. $\mathrm{lct}^{(d)}_{-}(D)$) exists if and only if $\mathscr{X}_{lN+d}(X,D){\neq \emptyset}$ for some $l \geq 0$ (resp. $l \leq -2$).
               \item If $\mathrm{lct}^{(d)}_{+}(D)$ (resp. $\mathrm{lct}^{(d)}_{-}(D)$) exists, then $\{C_{lN+d}/(lN+d)\}_{l\gg 1}$ (resp. $\{C_{-lN+d}/(lN-d)\}_{l \gg 1}$) is non-increasing as $l$ increases and 
               \begin{displaymath}
                   \lim_{l\to +\infty} C_{lN+d}/(lN+d) = \mathrm{lct}^{(d)}_+(D) \text{ (resp.} \lim_{l \to +\infty} C_{-lN+d}/(lN-d) = -\mathrm{lct}^{(d)}_-(D)\text{)}.
               \end{displaymath}
           \end{enumerate}
\end{theorem}

	\section{Plane curves}\label{Sec_curves}

	 Let $(C_1,\bm 0),\ldots,(C_r,\bm 0)$ be distinct reduced unibranch curve germs in $(\mathbb C^2,\bm 0)$, with reduced defining equations $f_1,\ldots,f_r$. Put $\underline D:=(\operatorname{div}(f_1),\ldots,\operatorname{div}(f_r))$. Throughout this section remarkable hyperplanes and zeta functions are local at $\bm 0$, as in Remark \ref{rmk_local}. These local constructions also apply to analytic plane
	 curve germs. Contact-order conditions are algebraic at
	 each finite jet level, and a finite sequence of point
	 blow-ups gives an embedded resolution with algebraic
	 strata over the origin. The same change-of-variables
	 argument therefore applies to analytic curve germs.

	\vspace{0.5em}

	The minimal log resolution of $f=f_1\cdots f_r$ is obtained by blowing up precisely the points where the reduced total transform is not SNC. It can be organized as follows.
	\begin{algorithm}\label{alg_curve}
	\begin{enumerate}
	\item If the original divisor is already SNC, stop. Otherwise blow up the origin. \label{alg1}
	\item Resolve the individual strict branches: blow up a point where a strict branch is singular, tangent to the reduced exceptional locus, or passes through an exceptional crossing. Repeat as necessary. \label{alg2}
	\item Separate the remaining smooth strict branches by blowing up their common points on the exceptional locus. Stop when the reduced total transform is SNC. \label{alg3}
	\end{enumerate}
	\end{algorithm}

	Let $\rho:Y\to(\mathbb C^2,\bm 0)$ be this resolution and $E:=\operatorname{Exc}(\rho)$ its exceptional divisor, empty when $\rho$ is the identity. The graph $\Gamma$ has one vertex for each exceptional component and an arrowhead for each strict branch; edges record intersections. For each branch, $\rho$ factors through its individual minimal log resolution $\rho_i:Y_i\to(\mathbb C^2,\bm 0)$. The intervening maps are point blowups. Thus the exceptional graph is obtained from the individual graph $\Gamma_i$ by subdividing edges and attaching trees, with the other strict branches added as arrowheads.

	Denote the strict transform of $C_j$ on $Y$ by $\widetilde C_j$, and on an intermediate model by $\overline C_j$. Let $\mathcal E:=\{E_i \mid i\in S_{\mathrm{exc}}\}$ and $\mathcal S:=\{\widetilde C_j \mid 1\leq j\leq r\}$. We also use a component itself as its index in expressions such as $N_{\{R,F\}}$. Take the exceptional index set disjoint from $\{1,\ldots,r\}$ and put $S:=S_{\mathrm{exc}}\sqcup\{1,\ldots,r\}$, with $E_j:=\widetilde C_j$, $N_{ji}:=\delta_{ji}$ and $\nu_j:=1$ for strict-transform indices. For an exceptional component $E_j$, put $N_{ji}:=\operatorname{ord}_{E_j}(f_i)$ and
	\[
	\mathbb H_j:=\{\bm s \mid \bm N_j\cdot\bm s+\nu_j>0\},\qquad
	H_j:=\{\bm s \mid \bm N_j\cdot\bm s+\nu_j=0\}.
	\]
	For a strict branch use $\widetilde{\mathbb H}_j:=\{s_j+1>0\}$ and $\widetilde H_j:=\{s_j+1=0\}$. We also write $H_{E_j} := H_j$.

	\begin{definition}We say that $R \in \mathcal E$ is a \textit{rupture divisor} (resp. \textit{leaf}, \textit{bridge divisor}) if its node has valency $\geq 3$ (resp. $=1$, $=2$) in $\Gamma$. Denote by $\mathcal{RP}$ (resp. $\mathcal{LF},\mathcal{BD}$) the set of rupture divisors (resp. leaves, bridge divisors).
	\end{definition}
	The main theorem of this section is as follows, which implies Theorem \ref{ThmD}.
	\begin{theorem}\label{thm_curves} Suppose $r \geq 2$, then:
		\begin{enumerate}
			\item For all $j=1,\ldots,r$, $\widetilde H_j $ is in $\mathrm{lctf}(\underline D)$, so they are remarkable.\label{curve1}
			\item For all $E_i \in \mathcal{LF} \cup \mathcal{BD}$, if $H_i\neq 
			H_j, \widetilde{H}_k$ for all $E_j \in \mathcal{RP}$ and $\widetilde C_k \in \mathcal S$, then $H_i$ is not a polar hyperplane of $Z_{\underline D}^{\mathrm{mot}}$. In particular, these hyperplanes are not remarkable. \label{curve2}
			\item For all $E_i \in \mathcal{RP}$, $H_i := \{\bm s\in \mathbb R^r \mid \bm N_i \cdot \bm s + \nu_i = 0\}$ is a remarkable hyperplane of $\underline D$. \label{curve3}
		\end{enumerate}
		Consequently, polar hyperplanes of $Z_{\underline D}^{\mathrm{mot}}$ are exactly remarkable hyperplanes of $\underline D$.
	\end{theorem}

	\subs{\bf Proof of parts \eqref{curve1} and \eqref{curve2} of Theorem \ref{thm_curves}.}
	For \eqref{curve1}, take $j=r$. If the resolution is the identity, the assertion is immediate. Otherwise every exceptional vector has all its coordinates positive. Since $r\geq2$, taking $s_r=-1$ and $s_1,\ldots,s_{r-1}$ sufficiently large and positive makes all exceptional inequalities and all other strict-transform inequalities strict. Hence $\widetilde H_r$ supports a facet.

	For \eqref{curve2}, fix $E_i\in\mathcal{LF}\cup\mathcal{BD}$
	satisfying the stated assumptions. For a hyperplane
	$H=\{\bm s\mid\bm a\cdot\bm s+b=0\}$ with
	$\bm a\in\mathbb Z_{\geq0}^r\setminus\{\bm0\}$ and $b>0$,
	and a vector $\bm\lambda\in\mathbb Z_{>0}^r$, put
	\[
	\mathrm{sp}(H,\bm\lambda)
	:=-\frac{b}{\bm a\cdot\bm\lambda}.
	\]
	Let $U_+\subset\mathbb P^{r-1}(\mathbb Q)$ consist of the
	classes of positive integral vectors $\bm\lambda$ satisfying
	\[
	\begin{aligned}
		\mathrm{sp}(H_i,\bm\lambda)
		&\neq\mathrm{sp}(H_j,\bm\lambda)
		&&\text{for all }E_j\in\mathcal{RP},\\
		\mathrm{sp}(H_i,\bm\lambda)
		&\neq\mathrm{sp}(\widetilde H_k,\bm\lambda)
		&&\text{for all }1\leq k\leq r.
	\end{aligned}
	\]
	By hypothesis, each excluded equality defines a proper
	projective hyperplane. Hence $U_+$ is nonempty and open
	in the Euclidean topology.
	
	For such a vector $\bm\lambda$, put
	$f_{\bm\lambda}:=\prod_{k=1}^r f_k^{\lambda_k}$.
	Its reduced zero locus is unchanged, so $\rho$ is also
	its minimal log resolution. Moreover,
	\[
	Z_{\divi(f_{\bm\lambda})}^{\mathrm{mot}}(T)
	=
	Z_{\underline D}^{\mathrm{mot}}
	(T^{\lambda_1},\ldots,T^{\lambda_r}).
	\]
	By \cite[Theorem 3.2]{Vey90}, whose argument also applies
	to motivic zeta functions, cf. \cite[Section 8.2]{Bla24},
	$\mathrm{sp}(H_i,\bm\lambda)$ is not a pole of this
	specialization.
	
	If $H_i$ were a polar hyperplane of
	$Z_{\underline D}^{\mathrm{mot}}$, Corollary
	\ref{poles_specialization} would give a dense subset
	$U\subset\mathbb P^{r-1}(\mathbb Q)$ along which
	$\mathrm{sp}(H_i,\bm\lambda)$ is a pole.
	Since $U_+$ is nonempty and open, $U\cap U_+\neq\varnothing$,
	a contradiction. Thus $H_i$ is not a polar hyperplane,
	and Theorem \ref{ThmA} implies that it is not remarkable.
	\hfill$\Box$

	\subs{\bf Rupture divisors.} We first study how rupture divisors arise in Algorithm
	\ref{alg_curve} and establish their numerical properties.
	These results will be used both to prove Theorem
	\ref{thm_curves}\eqref{curve3} and to determine the embedded
	topological type from the remarkable hyperplanes.
	%In the remaining parts, we will prove Theorem \ref{thm_curves}\eqref{curve3}. To this end, we need to understand how a rupture divisor can be created when running Algorithm \ref{alg_curve} and the numerical properties of these divisors.
	
	The minimal log resolution for each individual branch gives rupture components. We recall how the characteristic exponents change under blowups and how the strict transform meets the exceptional divisors.

	Choose a primitive parametrization $x=t^m$, $y=\sum_{j>m}\alpha_jt^j$ of a singular branch. Put $\beta_0:=m$, $e_0:=m$, and inductively
	\[
	\beta_{i+1}:=\min\{j \mid \alpha_j\neq0,\ e_i\nmid j\},\qquad
	e_{i+1}:=\gcd(e_i,\beta_{i+1}).
	\]
	The process stops at $e_g=1$; $(\beta_0;\beta_1,\ldots,\beta_g)$ is its characteristic exponent sequence. Subtracting integral powers of $x$ below $\beta_1$ makes $\operatorname{ord}_t y=\beta_1$. The strict branch after a blowup has the following exponents after normalizing its coordinates again.

\begin{lemma}\label{exp_blowup}
Let $Y\rightarrow \mathbb{A}^2$ be the blowup of the singular point of a singular branch with Puiseux exponents $(\beta_0;\beta_1,\dots,\beta_g)$. Then the Puiseux exponents of the strict transform are
\begin{align}
	&(\beta_0;\beta_1-\beta_0,\dots,\beta_g-\beta_0) & \beta_1>2\beta_0 \label{bl_ext} \\
	&(\beta_1-\beta_0;\beta_0,\beta_2-\beta_1+\beta_0,\dots,\beta_g-\beta_1+\beta_0) & \beta_1<2\beta_0 \text{ and } (\beta_{1}-\beta_0)\nmid \beta_0 \label{bl_over} \\		
	&(\beta_1-\beta_0;\beta_2-\beta_1+\beta_0,\dots,\beta_g-\beta_1+\beta_0) & \beta_1<2\beta_0 \text{ and } (\beta_{1}-\beta_0)\mid \beta_0 \label{bl_rap}
\end{align}
\end{lemma}
\begin{proof}
	In the chart $(x,y/x)$ the orders are $\beta_0$ and $\beta_1-\beta_0$. If the first is smaller, subtracting $\beta_0$ from the remaining exponents gives \eqref{bl_ext}. Otherwise interchange the two coordinates and reparametrize the coordinate of order $\beta_1-\beta_0$ as a pure power. The successive characteristic exponents become $\beta_0$ and $\beta_j-\beta_1+\beta_0$ for $j\geq2$; the first is omitted exactly when it is divisible by the new multiplicity. This gives \eqref{bl_over} and \eqref{bl_rap}. Equality $\beta_1=2\beta_0$ is impossible for a first characteristic exponent.
\end{proof}

	For the log resolution, one must also record how the strict transform meets the exceptional divisors, since making the branch smooth is not sufficient. A blowup at a point on exactly one exceptional component adds a new vertex joined to the vertex corresponding to that component; we say it \textit{extends the chain}. If the strict transform passes through the intersection of two exceptional divisors, blow up that point. We call such blowups \textit{zig-zags}. In the dual graph, the edge joining the two
	divisors is replaced by two edges through the vertex corresponding
	to the new exceptional divisor. Repeat these zig-zags until the strict transform meets an
	exceptional divisor at a point lying on no other exceptional divisor.
	
\begin{figure}[H]
		\centering
		\includegraphics[width=0.7\linewidth]{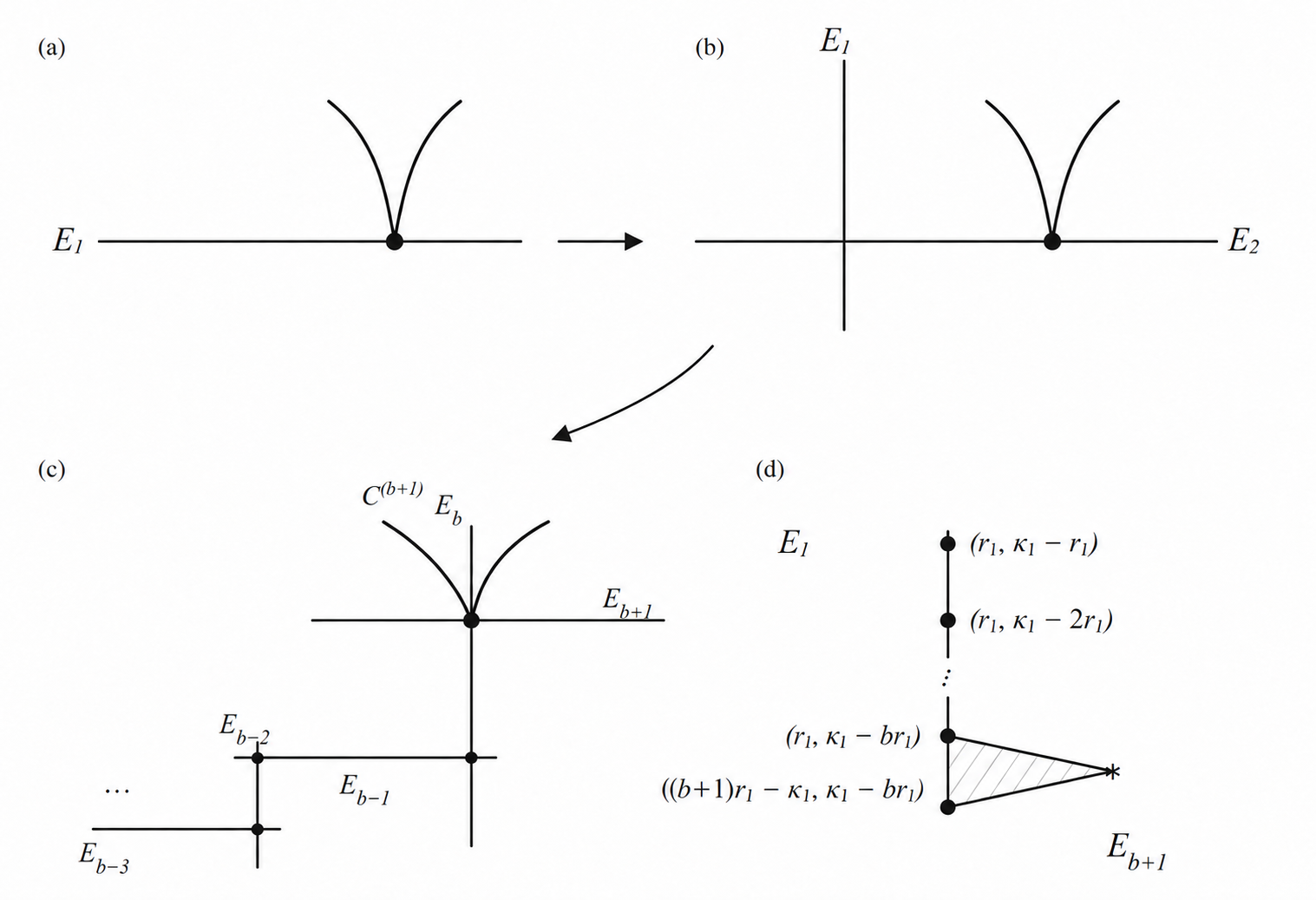}
		\caption[Resolution of unibranch curve: step 1]{extending the chain}
		\label{fig:unistep1}
	\end{figure}

\begin{figure}[H]
		\centering
		\includegraphics[width=0.7\linewidth]{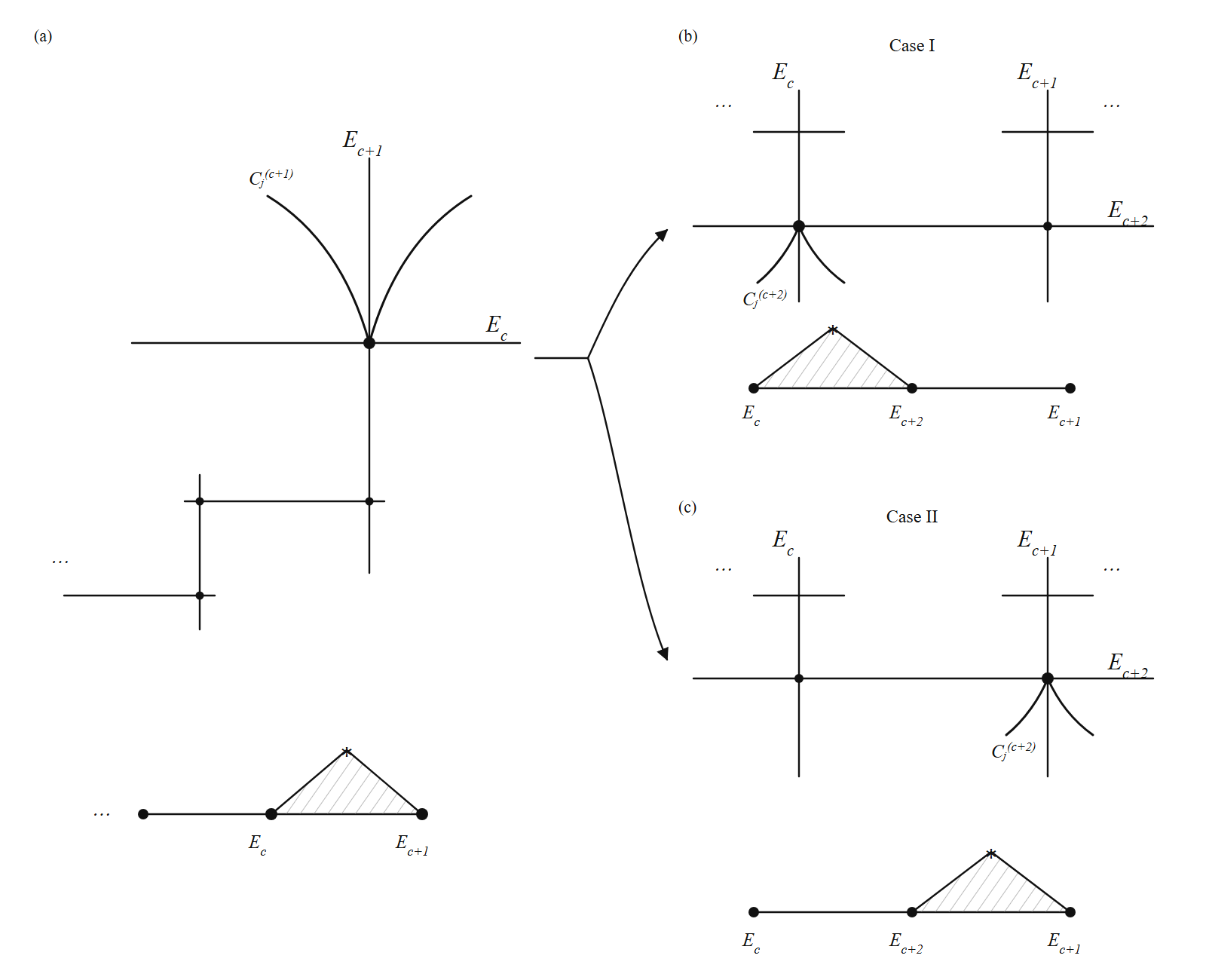}
		\caption[Resolution of unibranch curve: step 1]{zig-zags}
		\label{fig:unistep2}
	\end{figure}

\begin{figure}[H]
		\centering
		\includegraphics[width=0.7\linewidth]{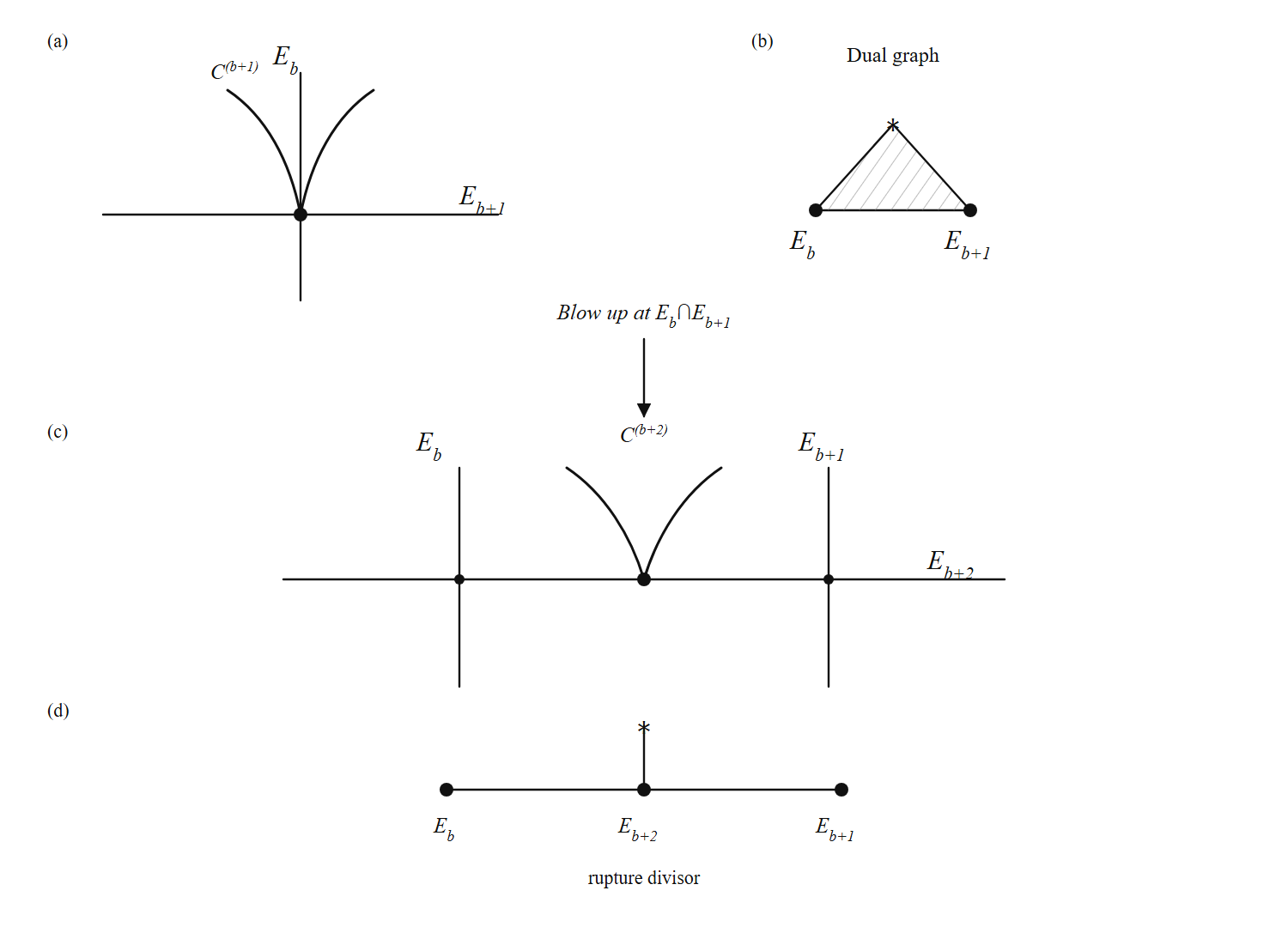}
		\caption[Resolution of unibranch curve: step 1]{rupture divisor}
		\label{fig:unistep3}
	\end{figure}

	The Euclidean algorithm governs this process, including the final blowups needed for transversality. Figures \ref{fig:unistep1}--\ref{fig:unistep3} depict these graph operations schematically. In Figure \ref{fig:unistep1}, $r_1=\beta_0$ and $\kappa_1=\beta_1$; the displayed pairs record local orders during the Euclidean algorithm, not the numerical data $(N_i,\nu_i)$. The superscript in $C^{(k)}$ counts blowups, and $b,c$ index the intermediate stages shown.

\begin{figure}[H]
		\centering
		\includegraphics[width=0.7\linewidth]{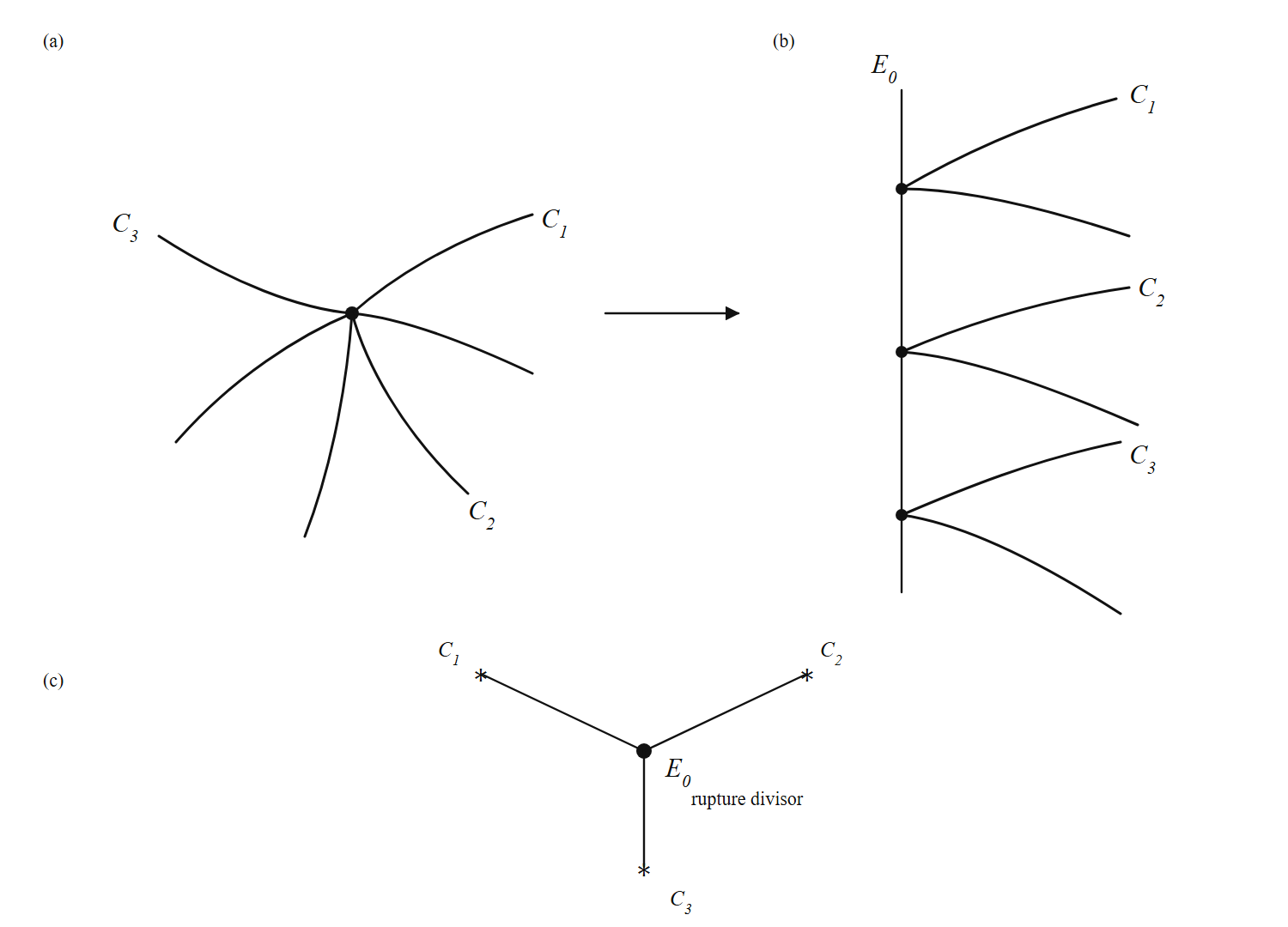}
		\caption[Rupture: case 1]{\ref{mult_1}}
		\label{fig:resol_case1}
	\end{figure}

\begin{figure}[H]
		\centering	
		\includegraphics[width=0.7\linewidth]{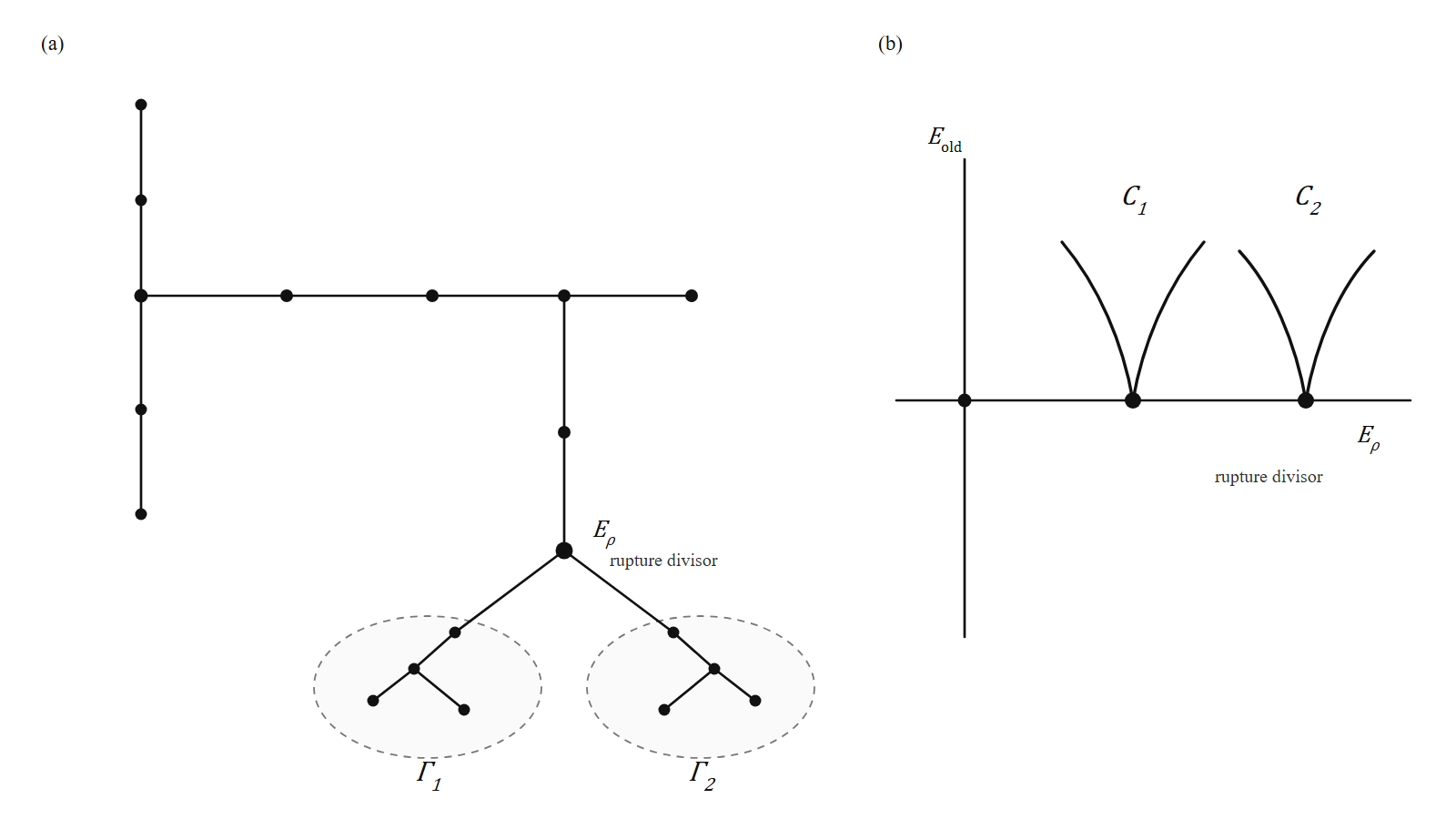}
		\caption[Rupture: case 2]{\ref{mult_2}}
		\label{fig:resol_case2}
	\end{figure}

\begin{figure}[H]
		\centering
		\includegraphics[width=0.7\linewidth]{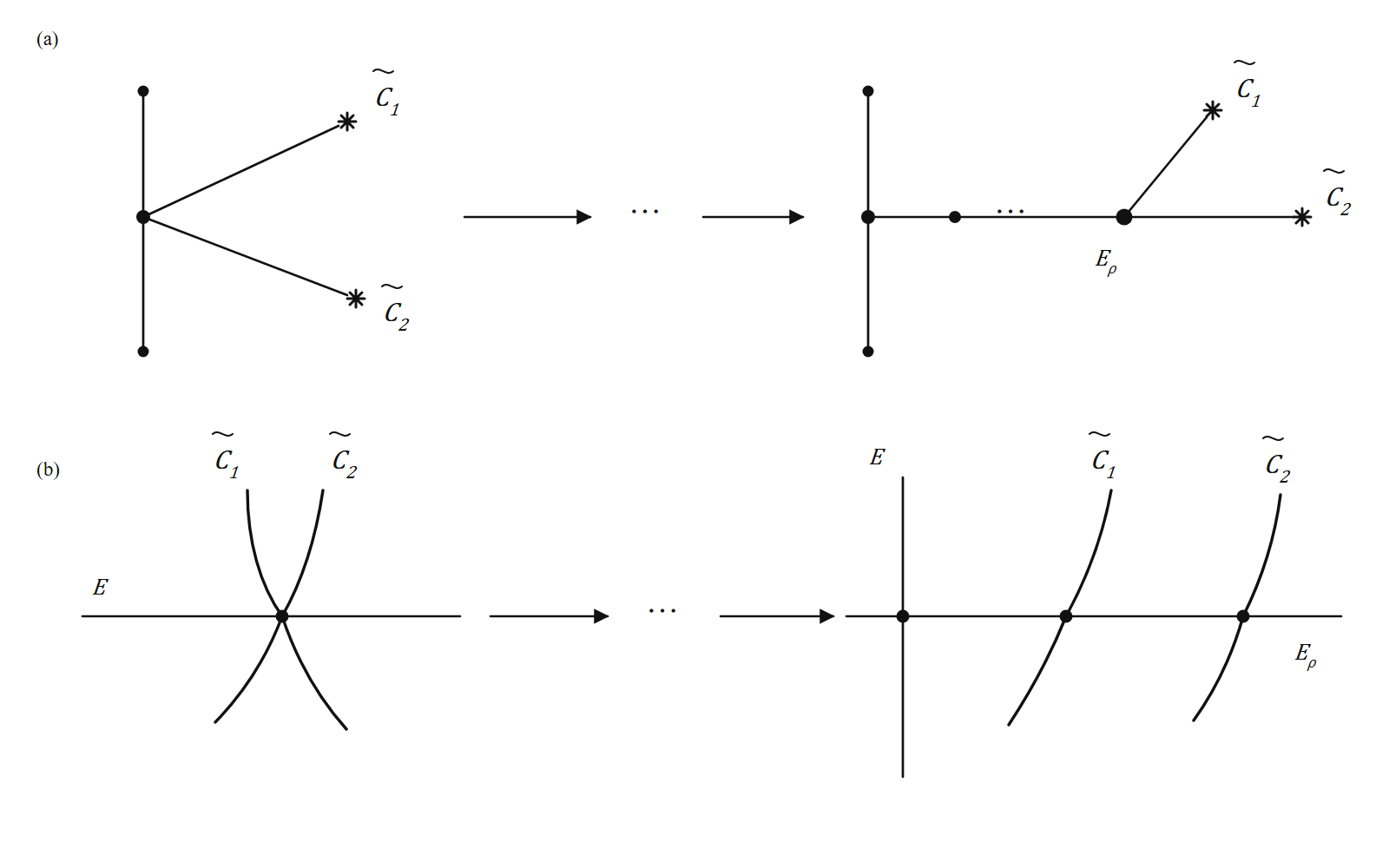}
		\caption[Rupture: case 3]{\ref{mult_3}}
		\label{fig:resol_case_3}
	\end{figure}

	From Algorithm \ref{alg_curve}, rupture divisors are produced in the following four cases:
	\begin{enumerate}[label = (Rup.\arabic*), start=0]
		\item In producing the minimal log resolution of each $C_j$, a rupture divisor of $\Gamma_j$ is again rupture in $\Gamma$. We call them \textit{intrinsic rupture components}. \label{mult_0}
		\item In the first blow up \eqref{alg1}, if there are at least three distinct tangent directions, then this exceptional divisor is rupture, see Figure \ref{fig:resol_case1}. \label{mult_1}
		\item A blowup at a \textit{free center}, i.e. a point lying on exactly one exceptional divisor, during \eqref{alg2} separates strict transforms of branches into at least two groups meeting the new exceptional divisor at distinct points away from the strict transform of the old exceptional divisor, see Figure \ref{fig:resol_case2}. \label{mult_2}
		\item In \eqref{alg3}, one blows up the intersection of strict transforms at the exceptional set until they are separated, see Figure \ref{fig:resol_case_3}.\label{mult_3}
	\end{enumerate}
	

%	\begin{remark}
%	To see that the cases exhaust all rupture components, consider the moment a component is created. A satellite center lies on two exceptional components. If a strict branch leaves that crossing on the new component, the blowup is required already in that branch's individual log resolution and gives an intrinsic rupture. Indeed, a smooth branch transverse to one exceptional component remains so after a blowup at its point, and cannot acquire an exceptional crossing from further such blowups. Subsequent blowups only subdivide its incident edges. At a free center there is one incoming edge, so a nonintrinsic rupture requires at least two distinct outgoing branch directions, giving \ref{mult_2} or \ref{mult_3}. At the origin there is no incoming edge, and at least three outgoing directions are needed, giving \ref{mult_1}. Distinct outgoing directions mean distinct smooth points on the new exceptional component; this is the property used below.
%\end{remark}

\begin{remark}
	To see that \ref{mult_0}--\ref{mult_3} exhaust all rupture
	components, let $R$ be a rupture component and consider the
	model immediately after its creation. Each point where $R$
	meets the other components of the reduced total transform
	contributes exactly one edge incident to $R$ in the final
	graph. Indeed, further blowups over such a point leave only
	one component meeting the strict transform of $R$ above that
	point, and no blowups are needed at the remaining points of $R$.
	
	Suppose that $R$ is created by blowing up the intersection
	of two exceptional components. Since $R$ is rupture, some
	strict branch must meet $R$ at a third point immediately
	after this blowup. The same blowup is required in the
	individual minimal log resolution of that branch, where $R$
	is already rupture. Thus $R$ is of type \ref{mult_0}.
	
	If $R$ is created by blowing up a point on exactly one
	exceptional component, then immediately after the blowup
	at least two strict branches must meet $R$ at distinct
	points away from the strict transform of that component.
	This gives \ref{mult_2} or \ref{mult_3}.
	
	Finally, if $R$ is created by the first blowup at the origin,
	the branches must have at least three distinct tangent
	directions for $R$ to be rupture, giving \ref{mult_1}.
\end{remark}
	
	Now we treat \ref{mult_0}-\ref{mult_3} separately.
		\begin{lemma}\label{finite_comparison}
		Let $\mathfrak C_i = \{\bm s \in \mathbb R^r \mid \bm c_i \cdot \bm s + d_i > 0\},i=1,\ldots,p$ be half spaces, where $\bm c_i \in \mathbb N^r \setminus \{\bm 0\}$ and $d_i > 0$. Suppose $c_{11} > 0$ and $d_1c_{j1} < d_jc_{11}$ for all $j\neq 1$, then {the hyperplane} $\mathrm{Cf}_1 := \{\bm s \in \mathbb R^r \mid \bm c_{{1}} \cdot \bm s + d_{{1}} = 0\}$ is {a hyperplane} spanned by a facet of $\overline{\bigcap_{i=1}^p \mathfrak C_i}$. More generally, if there exists $(\lambda_1,\ldots,\lambda_r) \in \mathbb Q_{\geq0}^r$ such that $\sum_{i=1}^r \lambda_i c_{1i} > 0$ and $d_1\sum_{i=1}^r \lambda_i c_{ji} < d_j\sum_{i=1}^r \lambda_i c_{1i}$ for all $j\neq 1$, then $\mathrm{Cf}_1$ is {a hyperplane} spanned by a facet of $\overline{\bigcap_{i=1}^p \mathfrak C_i}$.
	\end{lemma}
	\begin{proof}
		This is because $(-\lambda_1d_1/(\sum_{i=1}^r \lambda_i c_{1i}),\ldots,-\lambda_rd_1/(\sum_{i=1}^r \lambda_i c_{1i})) \in \mathrm{Cf}_1 \cap \bigcap_{i=2}^p \mathfrak C_i$.
	\end{proof}
	\begin{lemma}\label{avoidance_modulo}
		Suppose $A\subset \mathbb Z^r$ is a subgroup and $\bm a \in \mathbb Z^r \setminus A$. Then there exists a sublattice $L \subset \mathbb Z^r$ such that $\bm a\notin A+L$. Furthermore, if $A_1,\ldots,A_p$ are subgroups of $\mathbb Z^r$ {such that $\bm a\not\in A_i$ for all $i$}, then for $N$ divisible enough, we have $\bm a \notin A_i + (N \cdot \mathbb Z)^r$ for all $i$. 
	\end{lemma}
	\begin{proof}
	Smith normal form, using a unimodular change of ambient coordinates, gives $A=\bigoplus_i b_i\mathbb Z$ with $b_i\geq0$. If $b_i>0$ and $b_i\nmid a_i$, any $N$ divisible by $b_i$ separates that coordinate. If $b_i=0$ and $a_i\neq0$, any $N>|a_i|$ separates it. Taking a common multiple of these choices handles finitely many subgroups simultaneously. Unimodular coordinate changes preserve $N\mathbb Z^r$, and every further multiple of $N$ preserves nonmembership.
	\end{proof}
	\begin{lemma}\label{tail_comparison}
		Let $\pi:Y\to(\mathbb C^2,\bm 0)$ be a sequence of point blowups,
		and denote by $\overline C_j$ the strict transform of $C_j$ on $Y$.
		Let $R$ be an irreducible exceptional component of $\pi$,
		and let $P\in R$ lie on no other exceptional component.
		For $\bm \lambda = (\lambda_1,\ldots,\lambda_r)\in\mathbb Q_{\geq0}^r$, put
		\[
		D_\lambda:=\sum_j\lambda_j C_j,
		\qquad
		B:=\sum_j\lambda_j\overline C_j.
		\]
		Assume that $M_R:=\operatorname{ord}_R(D_\lambda)>0$,
		and set $c:=\nu_R/M_R$.
		If $c(B\cdot R)_P<1$, then $
		\nu_F-c\operatorname{ord}_F(D_\lambda)>0 $ for every irreducible exceptional divisor $F$ over $Y$
		whose center on $Y$ is $P$.	Moreover, $1-c\lambda_j>0$ whenever $P\in\overline C_j$.
	\end{lemma}
\begin{proof}
	Since $R$ is smooth and $c(B\cdot R)_P<1$, the pair $(R,(cB)|_R)$ is klt near $P$. By inversion of adjunction \cite[Theorem 5.50]{KM98}, the pair $(Y,R+cB)$ is plt near $P$.
	Since $cM_R=\nu_R$ and $P$ lies on no other exceptional
	component, the quantity
	$\nu_F-c\operatorname{ord}_F(D_\lambda)$ is the log discrepancy
	of this pair along $F$, and is therefore positive.
	Finally, if $P\in\overline C_j$, then
	$(\overline C_j\cdot R)_P\geq1$, so
	$c\lambda_j\leq c(B\cdot R)_P<1$.
\end{proof}

	\begin{proposition}\label{prop_case1}
	For a rupture divisor $E_i$ of type \ref{mult_0}, $H_i$ is remarkable.
\end{proposition}
\begin{proof}
	Let $\pi_1:Y_1\to(\mathbb C^2,\bm 0)$ be the minimal log
	resolution of $C_1$, and write the common resolution as
	$\rho=\pi_1\circ\tau$, where $\tau:Y\to Y_1$ is a sequence
	of point blowups.
	Denote the components of the reduced total transform of
	$C_1$ on $Y_1$ by $E_j^{(1)}$, $j\in S_1$, using the same
	indices for their strict transforms on $Y$.
	
	Suppose that $R=E_i$ corresponds to the $l$-th rupture
	component of the minimal log resolution of $C_1$. Put $
	d:=e_l, \alpha_R:=\frac{\nu_i}{N_{i1}},$
	and choose $N>d$ divisible by all $N_{j1}$ with $j\in S_1$.
	For each nonempty $J\subset S_1$ with
	$E_J^{(1)}:=\bigcap_{j\in J}E_j^{(1)}\neq\varnothing$, write $N_J^{(1)}:=\sum_{j\in J}\mathbb Z N_{j1}. $
	By \cite[Lemmas 5.2 and 5.3]{BPSZ}, if
	$d\in N_J^{(1)}+N\mathbb Z$, then
	\[
	\nu_j-\alpha_R N_{j1}\geq0\qquad(j\in J),
	\]
	with equality only for $j=i$.
	These lemmas also give an edge incident to $E_i^{(1)}$
	whose two multiplicities have gcd $d$.
	
	We now pass to $Y$ by the point blowups in $\tau$.
	Throughout this argument, consider only the strata of the
	reduced total transform of $C_1$, and put
	$q_j:=\nu_j-\alpha_R N_{j1}$ for each of its components.
	Let $E_*$ be the component created by one blowup.
	If the center lies on exactly one component $E_a$, the new
	strata involving $E_*$ have multiplicity subgroup
	$\mathbb Z N_{a1}$, and
	\[
	q_*=q_a+1.
	\]
	For both new strata $E_*$ and $E_*\cap E_a$, the required
	membership is $d\in\mathbb Z N_{a1}+N\mathbb Z$.
	The same formula applies to the old stratum $E_a$.
	If the center is $E_a\cap E_b$, their multiplicity subgroups
	are contained in
	$\mathbb Z N_{a1}+\mathbb Z N_{b1}$, and
	\[
	q_*=q_a+q_b.
	\]
	For $E_*$ the required membership is
	$d\in\mathbb Z(N_{a1}+N_{b1})+N\mathbb Z$;
	for either $E_*\cap E_a$ or $E_*\cap E_b$ it is
	$d\in\mathbb Z N_{a1}+\mathbb Z N_{b1}+N\mathbb Z$.
	Each of these implies
	$d\in\mathbb Z N_{a1}+\mathbb Z N_{b1}+N\mathbb Z$,
	so the inequalities for the old stratum $E_a\cap E_b$ apply.
	In either case, the relevant old values are nonnegative
	and only the value for $R$ can vanish, so $q_*>0$.
	This proves inductively that the same inequalities hold
	on $Y$.
	
	Moreover, suppose an edge joining $R$ to a component $E_j$
	has multiplicity gcd $d$. Blowing up $R\cap E_j$ preserves
	this gcd on the new edge adjacent to $R$, since
	\[
	\gcd(N_{i1},N_{i1}+N_{j1})
	=\gcd(N_{i1},N_{j1})=d.
	\]
	Starting with the edge on $Y_1$ described above, we therefore
	obtain on $Y$ a component $E_j\neq R$ such that
	\[
	R\cap E_j\neq\varnothing,
	\qquad d\in\mathbb Z N_{i1}+\mathbb Z N_{j1}+N\mathbb Z.
	\]
	
	Finally, let $E_I$ and $N_I\subset\mathbb Z^r$ have their
	original meanings on the common resolution, and set
	\[
	L:=N\mathbb Z\times\mathbb Z^{r-1},
	\qquad \bm a:=(d,0,\ldots,0).
	\]
	Then
	\[
	\bm a\in N_I+L
	\quad\Longleftrightarrow\quad
	d\in\sum_{j\in I}\mathbb Z N_{j1}+N\mathbb Z.
	\]
	The other strict transforms $\widetilde C_j$, $j>1$, have
	first multiplicity zero, so they do not change the subgroup
	on the right; their inequalities take value $1$ at
	$\bm s=(-\alpha_R,0,\ldots,0)$.
	Thus, among the inequalities defining
	$\mathrm{lct}^{\bm a,L}(\underline D)$, the inequality for
	$R$ occurs and is an equality at $\bm s$, while every other
	inequality is strict.
	Therefore $H_i$ supports a facet of
	$\overline{\mathrm{lct}^{\bm a,L}(\underline D)}$ and is
	remarkable.
\end{proof}
	
	\begin{proposition}\label{prop_case2}
	For a rupture divisor $E_i$ of type \ref{mult_1}, $H_i$ is remarkable.
\end{proposition}
\begin{proof}
	Let $R$ be the first exceptional divisor, and choose $C_1,C_2,C_3$ in three distinct tangent directions. Put $m_j:=\operatorname{mult}_{\bm 0}(C_j)$, $\lambda_j:=1/m_j$ for $1\leq j\leq3$, and $\lambda_j:=0$ for $j>3$.
	Use the notation $D_\lambda,B$ of Lemma \ref{tail_comparison} on the first blowup. Then $\operatorname{ord}_R(D_\lambda)=3$, $\nu_R=2$, and $c=2/3$. Each point of $R$ meets at most one chosen strict branch, so $(B\cdot R)_P$ is either $1$ or $0$.
	Lemma \ref{tail_comparison} makes all other inequalities defining $\mathrm{lct}(\underline D)$ strict at $\bm s=-c\bm\lambda$, while the inequality of $R$ is an equality. Thus $H_R$ supports a facet of $\overline{\mathrm{lct}(\underline D)}$. The three directions are illustrated in Figure \ref{fig:pf2}.
\begin{figure}
			\centering		
			\includegraphics[width=0.7\linewidth]{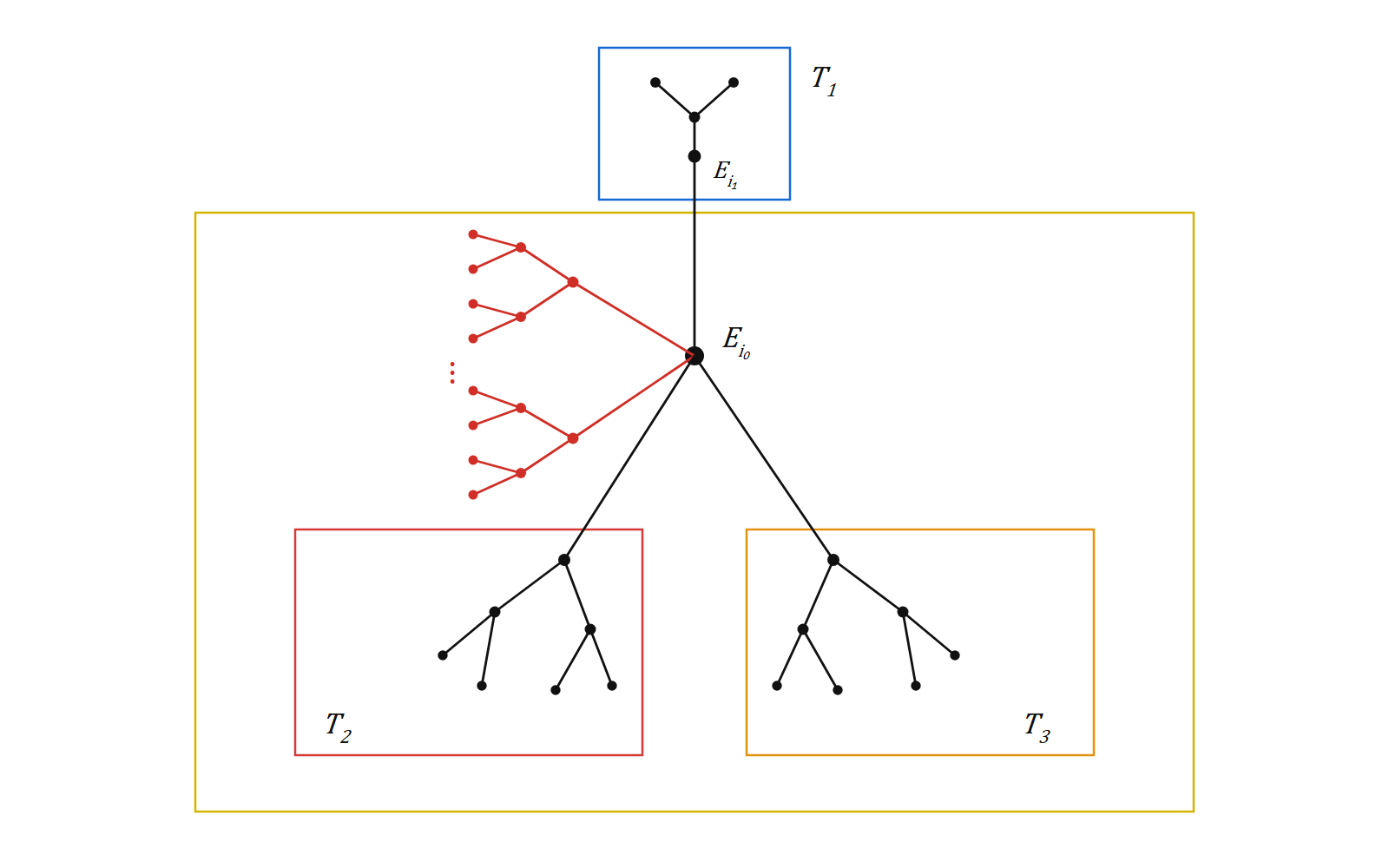}
			\caption{$\Gamma$ based on the first exceptional divisor}
			\label{fig:pf2}
		\end{figure}
\end{proof}

%\textcolor{blue}{I simplified the proof using Lemma \ref{exp_blowup}}.
\begin{lemma}\label{same_eucl_alg}
	Let $Y\to(\mathbb C^2,\bm 0)$ be a sequence of point blowups. In this lemma, index the exceptional components as $E_1,\ldots,E_n$ in the order in which they are created, and put $N_{ij}:=\operatorname{ord}_{E_i}(C_j)$. Suppose the strict transforms of $C_1,C_2$ meet the reduced exceptional locus only on $R:=E_n$, at smooth points of that locus. Put $m_j:=(R\cdot\overline{C}_j)>0$, where $\overline{C}_j$ is the strict transform of $C_j$ on $Y$. There are positive integers $a_i$ such that
	\[
	N_{ij}=a_i m_j\qquad (j=1,2).
	\]
	This proportionality persists for components obtained by blowing up away from the two strict transforms. If both strict transforms are smooth and transverse to $R$, then $m_1=m_2=1$. Finally, writing $A:=a_n$ and $B:=\nu_R$, one has $2A>B$ unless $R$ is the first exceptional divisor, in which case $A=1$ and $B=2$.
\end{lemma}
\begin{proof}
	Take a curvette of $R$, namely a smooth curve transverse to $R$ at a general point and disjoint from the other exceptional components, and let $a_i$ be the vanishing orders of its image. The exceptional intersection matrix $Q=(E_i\cdot E_k)$ is negative definite by \cite[p.~6]{Mum61}, and hence invertible. Write $\bm u_n:=(0,\ldots,0,1)^{\mathsf t}$. Intersecting the total transforms of principal divisors with each $E_i$ gives
	\[
	Q(N_{1j},\ldots,N_{nj})^{\mathsf t}=-m_j \bm u_n,\qquad
	Q(a_1,\ldots,a_n)^{\mathsf t}=-\bm u_n.
	\]
	Thus $N_{ij}=m_j a_i$. Extra blowups away from the strict branches preserve the relation by addition of the multiplicities of the components through the center.

	For the last assertion let $b_i\geq0$ be the multiplicities of the curvette image at the successive blowup centers, and let $E_i^*$ denote the total transform on $Y$ of the exceptional divisor of the $i$-th blowup. Then $\sum_i b_i E_i^*=\sum_i a_iE_i$. Its intersection vector is $-\bm u_n$, whereas $E_i^*\cdot E_k^*=-\delta_{ik}$. Pairing these identities with $E_i^*$ shows that the coefficient of $R$ in $E_i^*$ is $b_i$. Since $K_{Y/\mathbb C^2}=\sum_iE_i^*$, we obtain
	\[
	A=\sum_i b_i^2,\qquad B=1+\sum_i b_i.
	\]
	If $R$ is not the first exceptional divisor, at least two $b_i$ are positive integers. Hence $2\sum_i b_i^2>1+\sum_i b_i$. The first-blowup values are immediate.
\end{proof}

\begin{figure}[htbp]
			\centering			
			\includegraphics[width=0.7\linewidth]{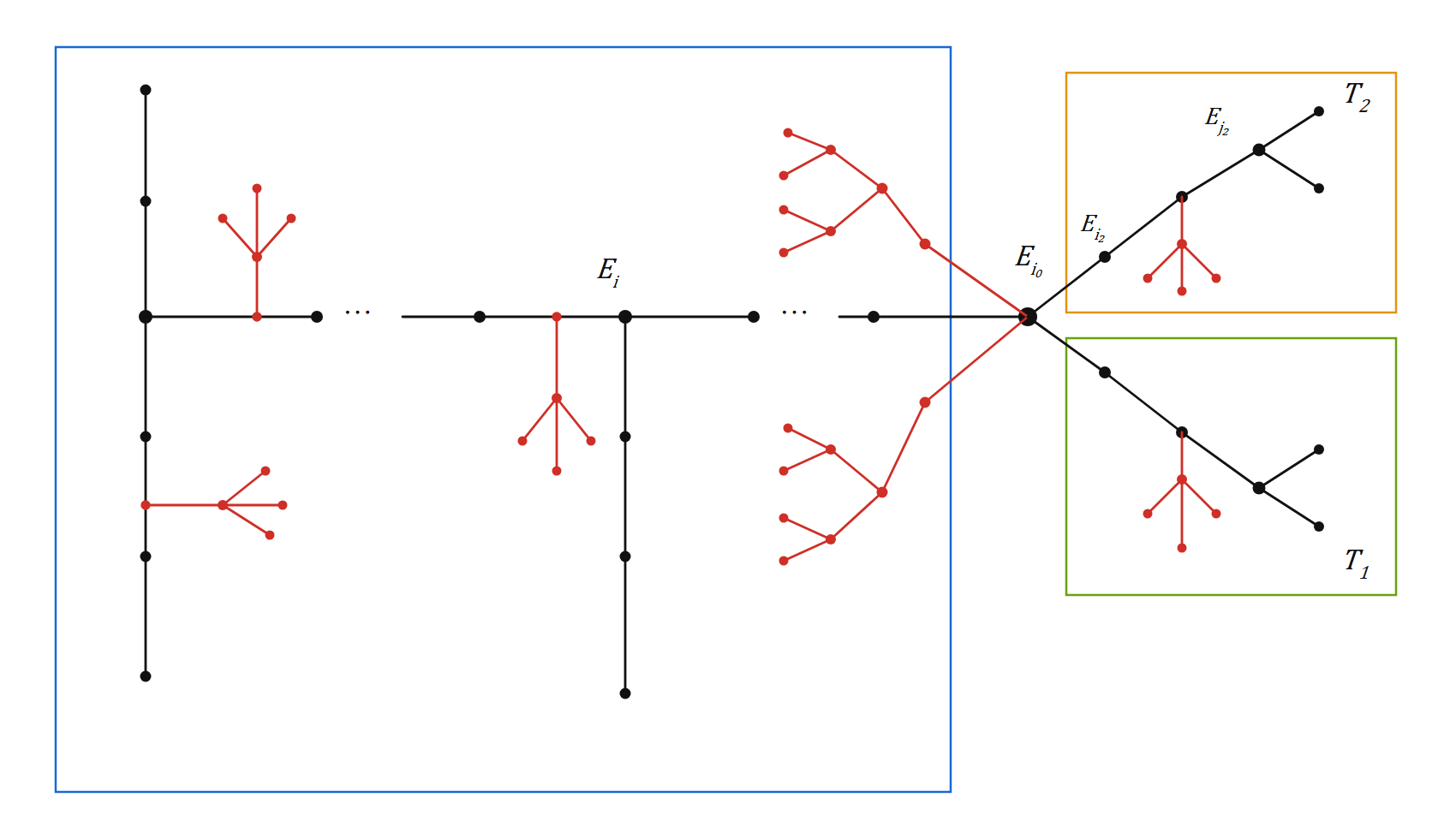}
			\caption{$\Gamma$ based on $R$, denoted by $E_{i_0}$ in the figure}
			\label{fig:pf3}
		\end{figure}

	\begin{proposition}\label{prop_case3}
	For a rupture divisor $E_i$ of type \ref{mult_2}, $H_i$ is remarkable. The same argument applies to every rupture divisor not of type \ref{mult_0} that is created by blowing up a point on exactly one exceptional component.
\end{proposition}
\begin{proof}
	Put $R:=E_i$. On the model just after its creation, choose two branches $C_1,C_2$ whose strict transforms $\overline C_1,\overline C_2$ meet $R$ at distinct points $P_1,P_2$ lying on no other exceptional component. Let $T_1,T_2$ be their components of the final graph after deleting $R$; include the edge from $R$ when speaking of a stratum in either tree. By Lemma \ref{same_eucl_alg}, put
	\[
	N_{R1}=Am_1,\quad N_{R2}=Am_2,\quad \nu_R=B,\quad 2A>B.
	\]
	Every component outside $T_1\cup T_2$ has its first two vanishing orders proportional to $(m_1,m_2)$; this includes the other strict branches, whose two orders are zero.

	Let $F_1$ be the neighbor of $R$ in $T_1$, possibly $\widetilde C_1$. Since it lies over $P_1$, its first two orders have the form
	\[
	(N_{F_1,1},N_{F_1,2})=(qAm_1+h,qAm_2),\qquad q \in \mathbb N,\ h \in \mathbb Z_{> 0}.
	\]
	Here $q=\operatorname{ord}_{F_1}(R)$ and $h=\operatorname{ord}_{F_1}(\overline C_1)$ on the intermediate model; for $F_1=\widetilde C_1$ they are $0$ and $1$. The positive value of $h$ persists under subdivisions of the edge toward $P_1$. Define the full integer vector
	\[
	\bm a:=N_{R2}\bm N_{F_1}-N_{F_1,2}\bm N_R\in N_{\{R,F_1\}}.
	\]
	Its first two coordinates are $(Am_2h,0)$, so $\bm a\notin N_I$ for every nonempty stratum $E_I$ all of whose components lie outside $T_1\cup T_2$. By Lemma \ref{avoidance_modulo}, choose $L=(N\mathbb Z)^r$ such that $\bm a\notin N_I+L$ for all these strata.
	On the other hand, $\bm a\in N_{\{R,F_1\}}+L$. Since $T_1,T_2$ are connected components of $\Gamma\setminus\{R\}$, every nonempty stratum $E_I$ with $\bm a\in N_I+L$ has all its components in $\{R\}\cup T_1\cup T_2$.

	Set $\lambda_1:=1/m_1$, $\lambda_2:=1/m_2$, and all other $\lambda_j:=0$. Then $M_R=2A$ and $c:=B/(2A)<1$. On the model just after creating $R$, put $B_\lambda:=\lambda_1\overline C_1+\lambda_2\overline C_2$. At each $P_j$ one has $(B_\lambda\cdot R)_{P_j}=1$.
	By Lemma \ref{tail_comparison}, at $\bm s=-c\bm\lambda$ the inequalities of all components in $T_1\cup T_2$ are strict, while that of $R$ is an equality. Thus every inequality defining $\mathrm{lct}^{\bm a,L}(\underline D)$ other than that of $R$ is strict. The inequality of $R$ occurs because $\bm a\in N_{\{R,F_1\}}+L$. Therefore $H_R$ supports a facet of $\overline{\mathrm{lct}^{\bm a,L}(\underline D)}$ and is remarkable. Figure \ref{fig:pf3} illustrates the trees $T_1$ and $T_2$.
\end{proof}

\begin{figure}[htbp]
			\centering		
			\includegraphics[width=0.7\linewidth]{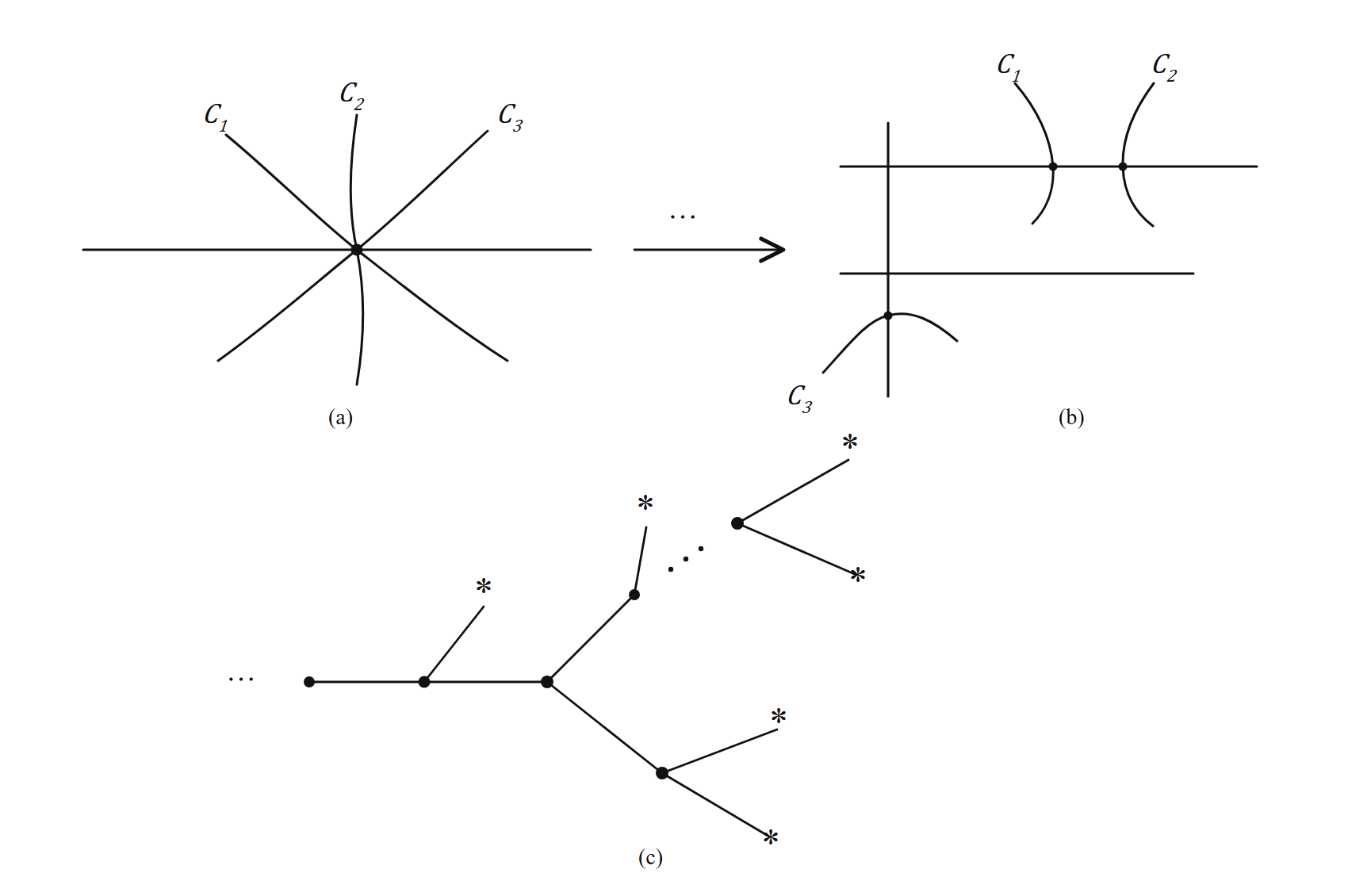}
			\caption{$\Gamma$ after \eqref{alg2}}
			\label{fig:pf41}
		\end{figure}

\begin{figure}[htbp]
			\centering			
			\includegraphics[width=0.7\linewidth]{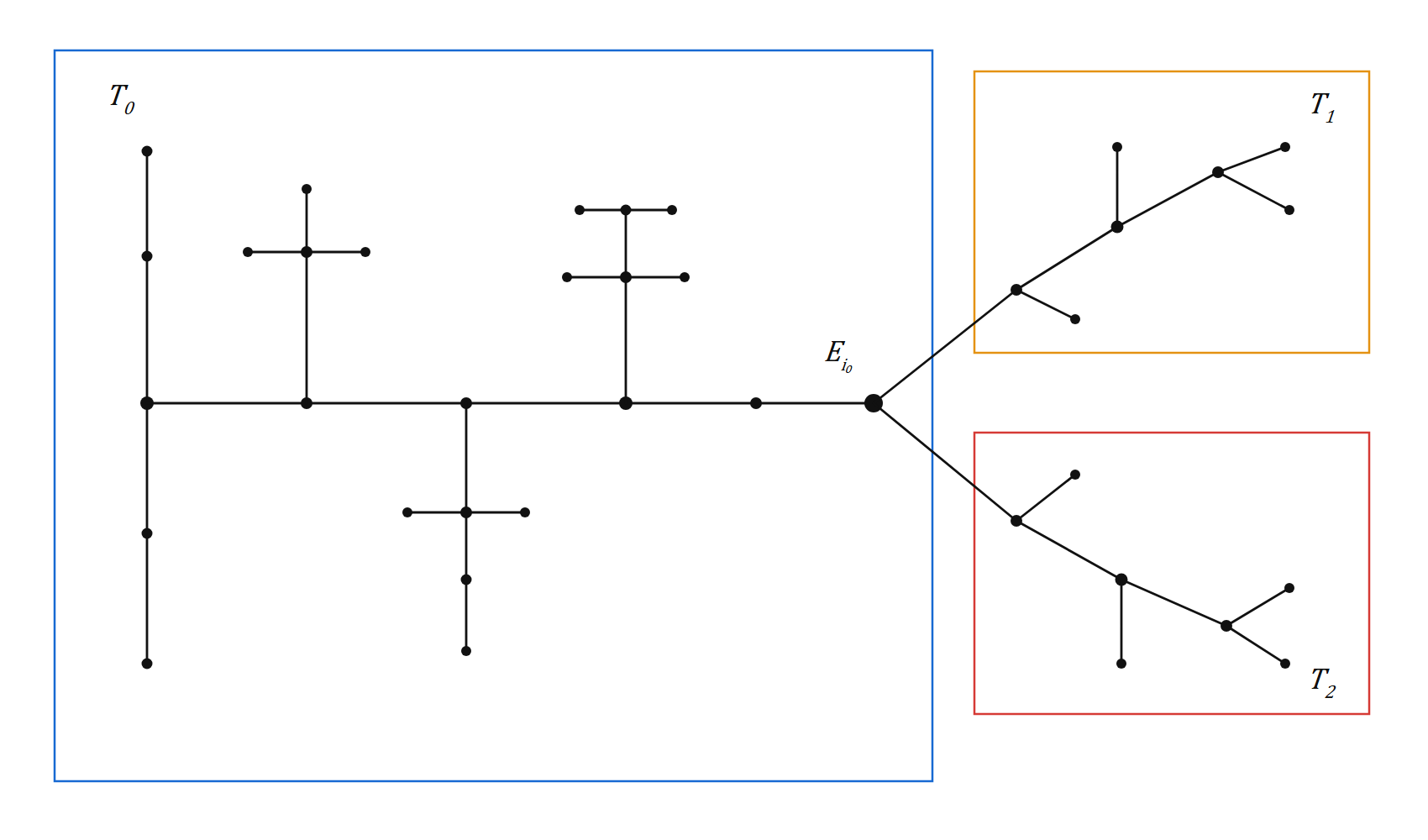}
			\caption{$\Gamma$ based on $R$, denoted by $E_{i_0}$ in the figure}
			\label{fig:pf42}
		\end{figure}

	\begin{proposition}\label{prop_case4}
	For a rupture divisor $E_i$ of type \ref{mult_3}, $H_i$ is remarkable.
\end{proposition}
\begin{proof}
	After the individual branch resolutions, the strict branches still requiring separation are smooth and transverse to the exceptional divisors, as in Figure \ref{fig:pf41}. At a new rupture $R$ choose two branches whose outgoing directions have separated, as in Figure \ref{fig:pf42}. The center is free, and on the model just after creating $R$ these two strict branches meet it transversely at distinct points. Thus the proof of Proposition \ref{prop_case3} applies with $m_1=m_2=1$, proving the claim.
	\end{proof}

	\subs{\bf Proof of Theorem \ref{thm_curves}\eqref{curve3}.} This is a consequence of Proposition \ref{prop_case1}-\ref{prop_case4}. $\hfill\Box$
	

%	\subs{\bf Determination of the embedded topological type.} We now show that adding one transverse smooth branch allows the remarkable hyperplanes to recover the embedded
%	topological type.

	\subs{\bf Proof of Theorem \ref{thm_topological_type}.} Put $C:=C_1\cup\cdots\cup C_r$ and
	$m_j:=\operatorname{mult}_{\bm 0}(C_j)$.
	Remove the hyperplanes $s_j+1=0$ for $0\leq j\leq r$ from the given remarkable
	set, and denote the remaining set by $\mathscr H$.
	Normalize each $H\in\mathscr H$ as
	\[
	s_0+\sum_{j=1}^r w_j(H)s_j+b(H)=0.
	\]
	By Theorem \ref{thm_curves}, these are the hyperplanes of
	the rupture components of $C\cup C_0$, with
	\[
	w_j(H_R)=\frac{N_{Rj}}{N_{R0}},
	\qquad
	b(H_R)=\frac{\nu_R}{N_{R0}},
	\qquad
	N_{R0}:=\operatorname{ord}_R(\ell).
	\]
	Here $b(H_R)>1$ since $(\mathbb C^2,C_0)$ is purely log terminal.
	
	We first recover the multiplicity of each branch at $\bm 0$ by
	\[
	m_j=\max\left(
	\{1\}\cup
	\left\{\frac{w_j(H)}{b(H)-1}:H\in\mathscr H\right\}
	\right).
	\]
	Indeed, inversion of adjunction shows that
	$(\mathbb C^2,C_0+C_j/m_j)$ is log canonical, so
	$\nu_R-N_{R0}-N_{Rj}/m_j\geq0$.
	For a singular branch $C_j$, equality holds at its first
	intrinsic rupture: if $\beta_1$ is its first characteristic
	exponent and $e_1=\gcd(m_j,\beta_1)$, then
	\[
	(N_{R0},N_{Rj},\nu_R)
	=\frac1{e_1}(m_j,m_j\beta_1,m_j+\beta_1).
	\]
	For a smooth branch, the additional value $1$ gives equality.
	
	Next, the distinct vectors
	\[
	\bigl(w_1(H),\ldots,w_r(H)\bigr)^{\mathsf t},
	\qquad H\in\mathscr H,
	\]
	are precisely the columns of the partial polar invariant
	matrix of $C$, cf. \cite[Definition 6.3]{GB00}. These columns are
	\[
	\left(
	\frac{(C_1,\gamma)_{\bm 0}}{\operatorname{mult}_{\bm 0}(\gamma)},
	\ldots,
	\frac{(C_r,\gamma)_{\bm 0}}{\operatorname{mult}_{\bm 0}(\gamma)}
	\right)^{\mathsf t},
	\]
	where $\gamma$ runs through the branches of a general polar
	curve, chosen transverse to $C_0$.
	To see the identification, apply \cite[Theorem 6.2]{GB00}.
	Along each terminal chain occurring there, which does not
	contain the first exceptional component, the intersection
	equations and the invertibility of the chain's intersection
	matrix show that $N_{Ej}/N_{E0}$ equals its value at the
	attaching rupture component.

	This also covers the case where the strict transform
	$\widetilde\gamma$ of $\gamma$ meets a crossing $P=E\cap F$,
	where $E$ and $F$ belong to such a terminal chain together
	with its attaching rupture component $R$.
	Let $\alpha=(E,\widetilde\gamma)_P$ and
	$\beta=(F,\widetilde\gamma)_P$.
	Since $C_0$ is transverse to $\gamma$ and $N_{Ej}/N_{E0} = N_{Fj}/N_{F0}= N_{Rj}/N_{R0}$, we have
	\[
	\frac{(C_j,\gamma)_{\bm 0}}
	{\operatorname{mult}_{\bm 0}(\gamma)}
	=
	\frac{\alpha N_{Ej}+\beta N_{Fj}}
	{\alpha N_{E0}+\beta N_{F0}}
	=
	\frac{N_{Rj}}{N_{R0}},
	\qquad 1\leq j\leq r.
	\]
%	At a smooth point of the exceptional divisor, the same
%	calculation has only one term.

	If $C$ has exactly two tangent directions, the additional
	polar column comes from the first exceptional component,
	which becomes rupture after adjoining $C_0$.
	For an ordinary node, this is the single column
	$(1,1)^{\mathsf t}$.
	
	Thus both the partial polar invariant matrix and the
	multiplicity of each branch at $\bm 0$ are determined.
	The conclusion now follows from \cite[Theorem 6.1]{GB00}.
	\hfill$\Box$
	
	%	\newpage
	%	
	%	
	%	
	%\clearpage

	\section{Examples of rational functions}\label{sec_rational_examples}

	We use the divisor notation of \ref{case_r1}. For $h=f/g$, put $D:=\divi(f)-\divi(g)$ and write $N_{f,i}:=\operatorname{ord}_{E_i}(f\circ\rho)$ and $N_{g,i}:=\operatorname{ord}_{E_i}(g\circ\rho)$, so $N_i=N_{f,i}-N_{g,i}$. Local statements use Remark \ref{rmk_local}. In each of the first two examples, the global remarkable numbers on $\mathbb A^2$ agree with the local remarkable numbers at the origin. All thresholds displayed there are understood for the stated integers $d$.

%	\subs{\bf Plane curves.}
	\begin{example}\label{ex_rational_node}
	Let $f:=y^2-x^2+x^3$ and $g:=y$. Blowing up the origin gives a log resolution of $D$. Let $E_0$ be the exceptional divisor, and let $E_1,E_2$ be the strict transforms of $\{f=0\}$ and $\{g=0\}$, respectively. The numerical data are
	\[
	\begin{array}{c|rrr}
	 &E_0&E_1&E_2\\ \hline
	N_{f,i}&2&1&0\\
	N_{g,i}&1&0&1\\
	N_i&1&1&-1\\
	\nu_i&2&1&1
	\end{array}
	\]
	Indeed, the two tangent directions of $\{f=0\}$ at the origin and the direction of $\{g=0\}$ are distinct. The other intersection of the two curves, at $(1,0)$, is transverse. The positive candidates are $2$ and $1$. Every positive admissible $I$ has $N_I=1$, and the candidate $1$ occurs for $E_1$, so $2$ is not remarkable. The only negative candidate is $-1$, from $E_2$. It follows that the remarkable numbers are $\{-1,1\}$ and, for every $d\in\mathbb Z$,
	\[
	\mathrm{lct}^{(d)}_+(D)=1,\qquad
	\mathrm{lct}^{(d)}_-(D)=-1.
	\]
	Both strict transforms meet the fiber over the origin, so the same computation applies locally.
%
%By Proposition \ref{prop_root of b-function}, the negative of a positive remarkable number, namely $-1$, is a root of Takeuchi's Bernstein-Sato polynomial of $f/g$ of order zero at the origin, i.e., $b_{f/g,\bm0,0}(-1)=0$. Alternatively, this can be deduced directly from the definition. Let $ D_{\mathbb C^2,\bm 0}$ denote the ring of germs of holomorphic differential operators at $\bm 0\in\mathbb C^2$. There exist $P_1(s),\ldots,P_n(s)\in D_{\mathbb C^2,\bm 0}[s]$ such that
%\[
%b_{f/g,\bm 0,0}(s)\left(\frac{f}{g}\right)^s=\displaystyle\sum_{k=1}^{n}P_{k}(s)\left(\frac{f}{g}\right)^{s+k}.
%\]
%Evaluating this functional equation at $s=-1$, the right hand side belongs to $\mathcal O_{\mathbb C^2,\bm 0}[g^{-1}]$, whereas $g/f\notin \mathcal O_{\mathbb C^2,\bm 0}[g^{-1}]$, since $f\nmid g^m$ for every $m\ge1$. Hence, $b_{f/g,\bm 0,0}(-1)=0$.
	\end{example}

	\begin{example}\label{ex_rational_cusp}
	Let $f:=y^2+x^5$, $g:=x^n$, and $D:=\divi(f)-\divi(g)$, where $n\geq1$. Four successive point blowups give a log resolution. Denote their exceptional divisors, in order of creation, by $E_0,E_1,E_2,E_3$, and the strict transforms of $\{f=0\}$ and $\{x=0\}$ by $E_4,E_5$. The dual graph is
	\[
	\xymatrix@R=12pt@C=22pt{
	 &&&E_4\ar@{-}[d]&\\
	 E_5\ar@{-}[r]&E_0\ar@{-}[r]&E_1\ar@{-}[r]&E_3\ar@{-}[r]&E_2
	}
	\]
	and the data are
	\[
	\begin{array}{c|rrrrrr}
	 &E_0&E_1&E_2&E_3&E_4&E_5\\ \hline
	N_{f,i}&2&4&5&10&1&0\\
	N_{g,i}&n&n&n&2n&0&n\\
	N_i&2-n&4-n&5-n&10-2n&1&-n\\
	\nu_i&2&3&4&7&1&1
	\end{array}
	\]
	For example, the four exceptional valuations have respective weights $(1,1),(1,2),(1,3),(2,5)$ on $(x,y)$. The first three blowups follow the successive points of the cusp; the fourth separates its strict transform from the intersection of $E_1$ and $E_2$. The displayed multiplicities follow by taking orders of $y^2+x^5$ and $x^n$, and the log discrepancies are the sums of the weights. Components with $N_i=0$ are retained in the simple normal crossings divisor, as allowed in \ref{not_pre}.

	For $n=1,2,3$ the remarkable numbers and thresholds are
	\[
	\begin{array}{c|c|c|c}
	n&\text{remarkable numbers of }D&\mathrm{lct}^{(d)}_+(D)&\mathrm{lct}^{(d)}_-(D)\\ \hline
	1&-1,\ 7/8&7/8&-1\\[2pt]
	2&-1/2,\ 1&1&-1/2\quad(2\mid d)\\[2pt]
	3&-1/3,\ 1&1&-1/3
	\end{array}
	\]
	The negative threshold in the middle row is undefined for odd $d$; all other entries hold for every $d\in\mathbb Z$.

%	For $n=1$, the smallest positive ratio is $7/8$, at $E_3$. The edge $\{E_1,E_3\}$ has $N_{\{1,3\}}=\gcd(3,8)=1$, so this minimum is attained for every $d$. The negative component $E_5$ has multiplicity $-1$. For $n=2,3$, all positive ratios are at least $1$, and $E_4$ has $(N_4,\nu_4)=(1,1)$, proving the positive entries. When $n=2$, the only negative component is $E_5$, with multiplicity $-2$, and its only neighbor is $E_0$, with multiplicity $0$. Thus every negative admissible stratum has gcd $2$ and $\beta_I=1/2$. When $n=3$, $E_5$ has the smallest negative ratio in absolute value, namely $1/3$, and the edge $\{E_5,E_0\}$ has gcd $\gcd(3,1)=1$. These facts prove the negative entries and, by Proposition \ref{prop_signed_thresholds}, the lists of remarkable numbers. All the strata used meet the fiber over the origin. As a direct check of the parity restriction, for $n=2$ one has
%	\[
%	\frac{y^2+x^5}{x^2}=\left(\frac yx\right)^2+x^3;
%	\]
%	its negative contact orders along arcs are necessarily even.
%
%\textcolor{blue}{For $n=1$, Proposition \ref{prop_root of b-function} implies that $-7/8$ is a root of $b_{f/g,0}$. Unlike in the previous example, this root cannot be deduced directly from the defining functional equation.} 
\end{example}

\begin{example}\label{counter_SMC_rational}
	Consider
	\[
	h(x,y)=\frac{(x+y^3)^2}{xy},\qquad
	D=2\{x+y^3=0\}-\{x=0\}-\{y=0\}.
	\]
	Three successive blowups give
	a log resolution of $D_{\mathrm{red}}$. Let $F,G,H$ denote the strict transforms of
	$\{x+y^3=0\}$, $\{x=0\}$, and $\{y=0\}$, respectively. The dual graph and numerical data are as follows.
	\begin{displaymath}
		\xymatrix@C=1.8em@R=1.2em{
			&&& F \\
			H \ar@{-}[r]
			& E_1 \ar@{-}[r]
			& E_2 \ar@{-}[r]
			& E_3 \ar@{-}[u] \ar@{-}[d] & \\
			&&& G
		}
	\end{displaymath}
	\vspace{0.5em}
	\begin{displaymath}
		\begin{array}{c|rrrrrr}
			&E_1&E_2&E_3&F&G&H\\ \hline
			N_{f,i}&2&4&6&2&0&0\\
			N_{g,i}&2&3&4&0&1&1\\
			N_i&0&1&2&2&-1&-1\\
			\nu_i&2&3&4&1&1&1
		\end{array}
	\end{displaymath}
	\vspace{0.5em}
	Therefore, the remarkable numbers are $-1,1/2,2$. Note that the local topological zeta function of $h$ at the origin is
	\begin{displaymath}
		Z_{h,\bm 0}^{\mathrm{top}}(s) = \frac{s+1}{(1-s)(2s+1)}
	\end{displaymath}
	and $-2$ is not a pole of it. Hence, the pole $-2$ of the motivic zeta function is cancelled along the specialization.
	
	The argument in the proof of \cite[Theorem 12]{GVL14}
	also addresses the weak monodromy conjecture for the motivic
	zeta function considered here in the plane curve case. In particular, three remarkable numbers $-1,1/2,2$ give monodromy eigenvalues. However, we will prove
	\begin{equation}
		b_{h,\bm 0,m}(-2)\neq0
		\quad\text{for every }m\geq0.
		\label{eq:counter-allpoints}
	\end{equation}
%	and 
%	\begin{equation}
%		b_{h,\bm 0}^A(-2)= 0.
%		\label{eq:A-allpoints}
%	\end{equation}
	Put $u=x+y^3$ and
	$g=y(u-y^3)$, so that $h=u^2/g$. In this example, all derivatives below are
	taken in the coordinates $(u,y)$.
	
	\vspace{0.5em}
	\noindent\textit{Case 1 of \eqref{eq:counter-allpoints}: $m = 0$.} Set
	\[
	B_0(s)=(s+1)(2s+1)(2s+5)
	\]
	and define the differential operators
	\[
	\begin{aligned}
	P_1(s)={}&\frac{s}{4}\partial_y^2-2(4s+1)y\partial_u +\frac{y}{2}\bigl((4s+19)u-(2s+5)y^3\bigr)\partial_u^2,
	\\
	P_2(s)={}&\frac1{12}\bigl(9s-5-7u\partial_u\bigr)
	\partial_u\partial_y +\frac{y^2}{12}\bigl(53(s+1)-16u\partial_u\bigr)
	\partial_u^2,
	\\
	P_3(s)={}&\frac1{12}\bigl(2s-5-3(u-y^3)\partial_u\bigr)
	\partial_u^2.
	\end{aligned}
	\]
	Direct differentiation gives
	\[
	B_0(s)h^s=P_1(s)h^{s+1}+P_2(s)h^{s+2}+P_3(s)h^{s+3}.
	\]
	Hence $b_{h,\bm 0,0}(s)\mid B_0(s)$, and $B_0(-2)=3\neq0$.
	
	\vspace{0.5em}
	\noindent\textit{Case 2 of \eqref{eq:counter-allpoints}: $m\geq1$.} Put
	\[
	K_m(s)=(2s-4m+4)(2s-4m+5),\qquad
	B_m(s)=(s+1)(2s+1)K_m(s),
	\]
	and define
	\[
	\begin{aligned}
	Q_m(s)={}&\frac{s+m}{2}\partial_y^2
	\bigl(s-m+1-(u-y^3)\partial_u\bigr)
	\\
	&+3(s+m)(s+1)y
	\bigl(y\partial_y+2s-4m+7\bigr)\partial_u
	\\
	&+\frac{K_m(s)}2y
	\bigl((u-y^3)\partial_u+s+m+2\bigr)\partial_u.
	\end{aligned}
	\]
	Direct differentiation gives
	\[
	B_m(s)g^{-m}h^s=Q_m(s)\bigl(g^{-m}h^{s+1}\bigr).
	\]
	Thus $b_{h,\bm 0,m}(s)\mid B_m(s)$, and
	$B_m(-2)=12m(4m-1)\neq0$.
%	
%	\vspace{0.5em}
%	\noindent\textit{Proof of \eqref{eq:A-allpoints}.} Let $R=\mathbb C\{u,y\}$ be the ring of convergent power
%	series at the origin, and let $D_R$ be its ring of
%	differential operators. Put
%	\[
%	D=\mathbb C[u,y]\langle\partial_u,\partial_y\rangle \subset D_R,
%	\qquad h=f/g,\qquad M=g^{-1}h^s.
%	\]
%	By the definition of $B(s):=b^A_{f/g,\bm 0}(s)$,
%	cf. \cite[Theorem A(2)]{AGLN25}, there exists $P(s)\in D_R[s]$ such that $
%	B(s)h^s=P(s)h^{s+1}$.
%	Since $h$ is weighted homogeneous of degree $2$ for the
%	weights $\operatorname{wt}(u)=3$ and
%	$\operatorname{wt}(y)=1$, retaining the component of $P(s)$
%	of weight $-2$ gives such an equation with $P(s)\in D[s]$. Here we assign weights $-3$ and $-1$ to $u$ and $y$, respectively.	Consequently,
%	\[
%	B(s)g\in J:=\operatorname{Ann}_{D[s]}(M)+D[s]f.
%	\]
%	A computation with \textsc{Singular} \cite{Singular} gives
%	\[
%	\begin{aligned}
%		J={}&D[s]\bigl(3u\partial_u+y\partial_y-2s+4\bigr)+D[s]\bigl((y^4-uy)\partial_y+(s+1)(4y^3-u)\bigr)
%		+D[s]u^2.
%	\end{aligned}
%	\]
%	Let $J_{-2}\subset D$ be the specialization of $J$ at $s=-2$.
%	Using the degree reverse lexicographic order with
%	$u>y>\partial_u>\partial_y$, \textsc{Singular} further computes
%	a Gr\"obner basis $\mathcal G_{-2}$ of $J_{-2}$ and gives
%	\[
%	\operatorname{NF}_{\mathcal G_{-2}}(g)=\frac23uy\neq0.
%	\]
%	Since $B(-2)g\in J_{-2}$, this implies $B(-2)=0$.
\end{example}

\subs{\bf Negative poles of zeta functions need not be roots of the Bernstein--Sato polynomials in \cite{AGLN25}.}\label{counter_SMC_Alvarez}
Let $h=f/g$ be a meromorphic germ at a point $p$ of a smooth complex variety $X$, with $f,g$ nonzero and coprime. Let $R=\mathcal O^{\mathrm{an}}_{X,p}$ be the ring of holomorphic germs at $p$ and let $\mathcal D_R$ be its ring of differential operators. For $\alpha\in\mathbb R_{\geq0}$, consider the free $R[(fg)^{-1}][s]$-module generated by the symbol $
e_\alpha=\frac{f^s}{g^{s+\alpha}}$ such that $
\partial(ae_\alpha)=\left(\partial a+a\left(s\frac{\partial f}{f}-(s+\alpha)\frac{\partial g}{g}\right)\right)e_\alpha$ for all derivatives $\partial$.
Following \cite[Definition 3.12]{AGLN25}, denote by $b^\alpha_{h,p}(s)$, when it exists, and $b_{h,p}^{A}(s)$ the monic polynomials of least degree satisfying, respectively,
\[
P(s)(fe_\alpha)=b^\alpha_{h,p}(s)e_\alpha,
\qquad
Q(s)(he_0)=b_{h,p}^{A}(s)e_0,
\qquad \text{for some } P(s),Q(s)\in\mathcal D_R[s].
\]
The first polynomial exists for all but finitely many $\alpha$, whereas the second always exists, by \cite[Theorem A]{AGLN25}. Any nonzero polynomial occurring on the right of the corresponding functional equation is divisible by the indicated $b$-function.

The following three examples are defined on complex affine spaces, with $p$ chosen to be their origins. They show that, for every $\alpha\geq0$, there is a germ with a nonzero $b^\alpha_{h,\bm0}$ and a negative pole of the motivic zeta function which is not a root of that polynomial. Example~\ref{ex_alvarez_generic} covers $\alpha\notin\{0,1\}$, and Examples~\ref{ex_alvarez_zero} and~\ref{ex_alvarez_one} cover $\alpha=0$ and $\alpha=1$, respectively. Example~\ref{ex_alvarez_zero} also disproves the analogous assertion for $b_{h,\bm0}^{A}$. We use the topological specialization of the motivic zeta function considered here. Each actual topological pole below is therefore also a pole of the motivic zeta function. We used \textsc{Singular}~\cite{Singular} to construct the differential operators $P_0(s)$ and $P_1(s)$ in Examples~\ref{ex_alvarez_zero} and~\ref{ex_alvarez_one}, respectively. The scripts verifying the resulting Bernstein--Sato equations are available in \cite{BCKSZ-data}.

%The functional equations and zeta computations are verified in the supplementary Python file \path{data_and_codes/counter_examples_Alvarez_et_al.py}.

\begin{example}\label{ex_alvarez_generic}
Consider
\[
h(x,y)=\frac{x^2+y^2}{y}.
\]
Blowing up the origin gives one exceptional divisor with $(N,\nu)=(1,2)$, meeting the two strict zero components and the strict pole component at three distinct points. Hence
\[
Z^{\mathrm{top}}_{h,\bm0}(s)
=\frac1{s+2}\left(-1+\frac2{s+1}+\frac1{1-s}\right)
=-\frac{s^2-s+2}{(s-1)(s+1)(s+2)}.
\]
In particular, $-2$ is a pole, with residue $-8/3$. Define
\[
\begin{aligned}
P_\alpha(s)={}&(s+1)\partial_y^2-\frac{x}{2}\partial_x\partial_y^2
+(s-\alpha+2)\partial_x^2-\frac{x}{2}\partial_x^3,\\
B_\alpha(s)={}&(s+1)(s-\alpha+2)(s-\alpha+3).
\end{aligned}
\]
Direct differentiation gives
\[
P_\alpha(s)\bigl((x^2+y^2)e_\alpha\bigr)=B_\alpha(s)e_\alpha.
\]
Thus $b^\alpha_{h,\bm0}(s)\mid B_\alpha(s)$, and
$B_\alpha(-2)=\alpha(1-\alpha)\neq0$ for every $\alpha\geq0$ with $\alpha\notin\{0,1\}$.
\end{example}

\begin{example}\label{ex_alvarez_zero}
Consider
\[
h(x,y)=\frac{f}{g}=\frac{x^2+y^3}{y(x+y)}.
\]
Three successive point blowups resolve $\{fg=0\}$. The exceptional divisors $E_1,E_2,E_3$, numbered in order of creation, have data
\[
\begin{array}{c|rrr}
&E_1&E_2&E_3\\ \hline
N_{f,i}&2&3&6\\
N_{g,i}&2&2&4\\
\nu_i&2&3&5\\
\chi(E_i^\circ)&-1&1&-1
\end{array}
\]
Let $F,G,H$ denote the strict transforms of $\{f=0\}$, $\{x+y=0\}$, and $\{y=0\}$, respectively. The dual graph of the total transform of $\{fg = 0\}$ is
\begin{displaymath}
	\xymatrix@C=1.8em@R=1.2em{
		& G & F & \\
		H \ar@{-}[r]
		& E_1 \ar@{-}[u] \ar@{-}[r]
		& E_3 \ar@{-}[u] \ar@{-}[r]
		& E_2.
	}
\end{displaymath}
 The strata formula thus gives
\[
Z^{\mathrm{top}}_{h,\bm0}(s)
=\frac{s^3+3s^2+5s+5}{(1-s)(s+1)(2s+5)}.
\]
Thus $-5/2$ is a pole, with residue $5/12$.

There is a polynomial-coefficient differential operator $P_0(s)$ of order $7$, given in \cite[\path{Alvarez_et_al_P0.txt}]{BCKSZ-data}, such that
\[
P_0(s)(fh^s)=(s+1)h^s.
\]
Consequently $b^0_{h,\bm0}(s)\mid s+1$. Composing $P_0(s)$ on the right with multiplication by $g$ also gives
\[
(P_0(s) \cdot g)h^{s+1}=(s+1)h^s,
\]
so $b_{h,\bm0}^{A}(s)\mid s+1$. So neither $b_{h,\bm 0}^0$ nor $b_{h,\bm 0}^{A}$ vanishes at $-5/2$.
\end{example}

\begin{example}\label{ex_alvarez_one}
Consider
\[
h(x,y,u,v)=\frac{f}{g}=\frac{x^2+y^3+uv}{y(x+y)}.
\]
First blow up the origin, obtaining $E_1\simeq\mathbb P^3$. Next blow up the vertex $[0:1:0:0]$ of the quadric $\{x^2+uv=0\}\subset E_1$, obtaining $E_2$. Finally blow up the smooth center $E_1'\cap E_2\simeq\mathbb P^2$, where $E_1'$ is the strict transform of $E_1$, obtaining $E_3$. This gives a log resolution $\rho$ of the divisor $\{fg=0\}$. Denote the final strict transforms of $\{f=0\}$, $\{y=0\}$, and $\{x+y=0\}$ by $F,H_1,H_2$, respectively. The numerical data are
\[
\begin{array}{c|rrrrrr}
&E_1&E_2&E_3&F&H_1&H_2\\ \hline
N_{f,i}&2&3&6&1&0&0\\
N_{g,i}&2&2&4&0&1&1\\
N_i&0&1&2&1&-1&-1\\
\nu_i&4&7&11&1&1&1
\end{array}
\]
Writing $\chi_I=\chi(E_I^\circ\cap\rho^{-1}(\bm0))$ and using $F,H_1,H_2$ also as component indices, the nonzero Euler characteristics are
\[
\begin{array}{c|r}
I&\chi_I\\ \hline
\{1\},\ \{3\}&-1\\
\{2\},\ \{1,3\},\ \{3,F\},\ \{1,H_1\},\ \{1,H_2\}&1\\
\{2,3\}&3\\
\{1,3,F\},\ \{1,F,H_1,H_2\}&2
\end{array}
\]
All remaining strata have Euler characteristic zero. Hence
\[
\begin{aligned}
Z^{\mathrm{top}}_{h,\bm0}(s)
={}&-\frac14+\frac1{s+7}-\frac1{2s+11}
+\frac1{2(1-s)}+\frac1{2(s+1)(1-s)^2}\\
&+\frac1{4(2s+11)}+\frac3{(s+7)(2s+11)}
+\frac3{2(s+1)(2s+11)}\\
={}&-\frac{s^4+4s^3+4s^2-22}{2(s-1)^2(s+1)(2s+11)}.
\end{aligned}
\]
Thus $-11/2$ is a pole, with residue $11/24$.

Put $e_1=f^s/g^{s+1}$. There is a polynomial-coefficient differential operator $P_1(s)$ of order $7$, given in \cite[\path{Alvarez_et_al_P1.txt}]{BCKSZ-data}, such that
\[
P_1(s)(fe_1)=(s+1)e_1.
\]
Hence $b^1_{h,\bm0}(s)\mid s+1$.  Thus $ b^1_{h,\bm0}(-11/2)\neq0$.
\end{example}

\subs{\bf Unboundedness of remarkable numbers for rational functions.} \label{ex_unbounded_rational} In the case of regular functions, remarkable numbers always lie in $(0,1]$, cf. \cite[Theorem 1.7(3)]{BPSZ}. However, such boundedness fails badly for rational functions, see the following examples.
	
\begin{example}
		For an integer $m\geq 4$, put
	\[
	h_m:=\frac{(x+y^m)^2}{xy^{m-2}},
	\qquad D_m:=\divi(h_m).
	\]
	We show that $(m+1)/2$ is a remarkable number of $D_m$.
	
	Blow up the origin and then the successive common points of the strict
	transforms of $\{x=0\}$ and $\{x+y^m=0\}$, for a total of $m$ blowups.
	This gives a log resolution. Denote the exceptional divisors, in order
	of creation, by $E_1,\ldots,E_m$, and the strict transforms of
	$\{x+y^m=0\},\{x=0\},\{y = 0\}$ by $F,G,H$, respectively. The dual graph is
	\[
	\xymatrix@R=10pt@C=18pt{
		&&&&F\ar@{-}[d]\\
		H\ar@{-}[r]&E_1\ar@{-}[r]&\cdots\ar@{-}[r]&E_{m-1}\ar@{-}[r]&E_m\\
		&&&&G\ar@{-}[u]
	}
	\]
	
	In the chart $x=u_k y^k$, one has
	\[
	h_m=y^{k-m+2}\frac{(u_k+y^{m-k})^2}{u_k},
	\qquad
	dx\wedge dy=y^k\,du_k\wedge dy.
	\]
	Consequently,
	\[
	(N_{E_k},\nu_{E_k})=(k-m+2,k+1),
	\qquad
	(N_F,\nu_F)=(2,1),
	\qquad
	(N_G,\nu_G)=(-1,1).
	\]
	The strict transform of $\{y=0\}$ has multiplicity
	$-(m-2) < 0$. Thus the only components with positive multiplicity are
	$F,E_{m-1},E_m$, and the positive candidate values are
	\[
	\frac12,\qquad m,\qquad \frac{m+1}{2}.
	\]
	
	The components $E_m$ and $G$ intersect. For the corresponding stratum
	$I=\{E_m,G\}$, we have
	\[
	N_I=\gcd(2,1)=1,
	\qquad
	\alpha_I=\frac{m+1}{2}.
	\]
	If a positive admissible stratum $J$ satisfies
	$\alpha_J<\alpha_I$, then it contains $F$ and has $\alpha_J=1/2$.
	Since $F$ meets only $E_m$, such a stratum is either $\{F\}$ or
	$\{F,E_m\}$, and in both cases $N_J=2$.
	As $2\nmid N_I$, Definition~\ref{def_signed_remarkable} shows that
	$(m+1)/2$ is remarkable.
	
	These values tend to $+\infty$ as $m$ increases. Replacing $h_m$
	with $1/h_m$ gives the remarkable numbers $-(m+1)/2$, so there is
	also no uniform lower bound. All strata used above have images
	containing the origin, and the same conclusions hold for local
	remarkable numbers at the origin.
\end{example}

\subs{\bf Reduced hyperplane arrangements.} Define by $\mathrm{lct}_{+}(D):=\min_{i\in S, N_{i}>0}{(\nu_{i}/N_{i})}$ and $\mathrm{lct}_{-}(D):=-\min_{i\in S, N_{i}<0}{(-\nu_{i}/N_{i})}$ the positive and negative log canonical thresholds of $D$, respectively, whenever they exist.

\begin{proposition}
Let $f,g\in \mathbb{C}[x_{1},...,x_{n}]$ be reduced products of linear polynomials with no common factor. Then the positive remarkable numbers, if any, consist precisely of $\mathrm{lct}_{+}(D)$, and the negative remarkable numbers, if any, consist precisely of $\mathrm{lct}_{-}(D)$. In particular, $f/g$ has at most two remarkable numbers.
\end{proposition}

\begin{proof}
We adapt the argument of \cite[6.7. Hyperplane arrangements]{BPSZ} to the rational setting. Let $S$ be the set of proper edges of the hyperplane arrangement associated with $fg$. Blowing up the (strict transforms of the) edges in increasing dimension yields a log resolution. Moreover, the subsets $I\subset S$ such that $E_{I}\not=\emptyset$ are naturally identified with flags $W_{1}\subsetneq...\subsetneq W_{\mid I\mid}$ of edges. For each edge $W\in S$, $\nu_{W}$ is the codimension of $W$, and $N_{W}$ is the difference between the numbers of hyperplanes of $f$ and $g$ containing $W$. Suppose that $\mathrm{lct}_{+}(D)$ exists, and choose an edge $W_{0}$ such that $\mathrm{lct}_{+}(D)=\nu_{W_{0}}/N_{W_{0}}$. Since $N_{W_{0}}>0$, there is a hyperplane $H$ of $f$ containing $W_{0}$. If $W_{0}\neq H$, then $W_{0}\subsetneq H$ is a flag; if $W_{0}=H$, we take the one-element flag. Let $I$ denote the corresponding set of edges. Since $f$ and $g$ are reduced and have no common factor, $N_{H}=1$, and hence $N_{I}=1$. Moreover, we have $\alpha_{I}=\mathrm{lct}_{+}(D)$. It follows from the definition of a positive remarkable number that no positive remarkable number can be strictly greater than $\mathrm{lct}_{+}(D)$, since $N_{I}=1$ divides $N_{J}$ for every positive admissible $J$. Thus $\mathrm{lct}_{+}(D)$ is the only positive remarkable number. The negative case is analogous. 

\end{proof}

\section{Examples of regular maps}\label{sec_regular_examples}

\begin{example}\label{ex_auxiliary_branch_topology}
The auxiliary branch in Theorem~\ref{thm_topological_type} cannot in
general be omitted. For $c=1,2$, put
\[
f_c=(y^2-x^{15})^2-cx^4y^{35},
\qquad
g_c=y^5-x^{12}-cx^2y^6.
\]
All four germs are irreducible at $\bm0$. The minimal log resolutions
of $\{f_1f_2=0\}$ and $\{g_1g_2=0\}$ each have exactly two rupture
components, with numerical data $(N_{R1},N_{R2},\nu_R)$ as follows:
\[
\begin{array}{c|cc}
&\text{first rupture}&\text{second rupture}\\ \hline
(f_1,f_2)&(60,60,17)&(1066,1066,507)\\
(g_1,g_2)&(60,60,17)&(82,82,39)
\end{array}
\]
The dual graphs are shown below. In both graphs, the indicated chain
lengths count vertices, including the displayed endpoints. The solid arrows
represent the strict transforms $F_c$ of $\{f_c=0\}$ and $G_c$ of
$\{g_c=0\}$. The dashed arrow indicates the strict transform $C_0$
of the auxiliary line $\{x=0\}$, if it is adjoined.
\begin{center}
\begin{minipage}{\linewidth}
\centering
\begin{tikzpicture}[x=1cm,y=1cm,>=stealth,
  every node/.style={font=\small},
  vertex/.style={circle,fill,inner sep=1.3pt}]
  \node[anchor=west] at (-0.7,1.45)
    {$(f_1,f_2)$: $247$ exceptional components};
  \node[vertex] (a1) at (0,0) {};
  \node[vertex] (a6) at (1.4,0) {};
  \node[vertex] (a7) at (2.5,0) {};
  \node[vertex] (r1) at (3.7,0) {};
  \node[above right=3pt] at (r1) {$R_1$};
  \node[vertex] (b) at (3.7,1.05) {};
  \node[vertex] (c1) at (4.8,0) {};
  \node[vertex] (c235) at (6.6,0) {};
  \node[vertex] (c236) at (7.7,0) {};
  \node[vertex] (r2) at (8.9,0) {};
  \node[above left=3pt] at (r2) {$R_2$};
  \node[vertex] (d) at (8.9,1.05) {};
  \draw (a1)--(0.35,0) (1.05,0)--(a6)--(a7)--(r1)--(c1)
    (c1)--(5.25,0) (6.15,0)--(c235)--(c236)--(r2)
    (r1)--(b) (r2)--(d);
  \node at (0.7,0) {$\cdots$};
  \node at (5.7,0) {$\cdots$};
  \node[above=5pt] at (0.7,0) {$6$ vertices};
  \node[above=5pt] at (5.7,0) {$235$ vertices};
  \draw[->] (r2)--(10.15,0.55) node[right] {$F_1$};
  \draw[->] (r2)--(10.15,-0.55) node[right] {$F_2$};
  \draw[->,dashed] (a1)--(-0.7,-0.65) node[left] {$C_0$};
\end{tikzpicture}

\medskip
\begin{tikzpicture}[x=1cm,y=1cm,>=stealth,
  every node/.style={font=\small},
  vertex/.style={circle,fill,inner sep=1.3pt}]
  \node[anchor=west] at (-0.7,2.2)
    {$(g_1,g_2)$: $28$ exceptional components};
  \node[vertex] (a1) at (0,0) {};
  \node[vertex] (a2) at (1.2,0) {};
  \node[vertex] (a3) at (2.5,0) {};
  \node[vertex] (r1) at (3.7,0) {};
  \node[above right=3pt] at (r1) {$R_1$};
  \node[vertex] (b1) at (3.7,0.9) {};
  \node[vertex] (b2) at (3.7,1.8) {};
  \node[vertex] (c1) at (4.8,0) {};
  \node[vertex] (c21) at (6.6,0) {};
  \node[vertex] (r2) at (8.9,0) {};
  \node[above left=3pt] at (r2) {$R_2$};
  \draw (a1)--(a2)--(a3)--(r1)--(c1)
    (c1)--(5.25,0) (6.15,0)--(c21)--(r2)
    (r1)--(b1)--(b2);
  \node at (5.7,0) {$\cdots$};
  \node[above=5pt] at (5.7,0) {$21$ vertices};
  \draw[->] (r2)--(10.15,0.55) node[right] {$G_1$};
  \draw[->] (r2)--(10.15,-0.55) node[right] {$G_2$};
  \draw[->,dashed] (a1)--(-0.7,-0.65) node[left] {$C_0$};
\end{tikzpicture}
\end{minipage}
\end{center}

Since $(1066,1066,507)=13(82,82,39)$, Theorem~\ref{thm_curves}
shows that the two complete splittings have the same local remarkable
hyperplanes:
\[
s_1+1=0,\qquad s_2+1=0,\qquad
60(s_1+s_2)+17=0,\qquad 82(s_1+s_2)+39=0.
\]
Nevertheless, their embedded topological types differ, since each
$f_c$ has multiplicity $4$, whereas each $g_c$ has multiplicity $5$.
Adjoining the transverse line $\{x=0\}$ distinguishes the two cases:
the first rupture hyperplanes become, respectively,
\[
2s_0+60(s_1+s_2)+17=0,
\qquad
5s_0+60(s_1+s_2)+17=0.
\]
\end{example}

\begin{example}\label{ex_rmkhyp_no_nonnegative_specialization}
Consider the germ
\[
F=(xy^2,xy^3):\mathbb C^2 \to \mathbb C^2.
\]
The identity is a log resolution, and
\[
\mathrm{lct}(\underline D)
=\{(s_1,s_2)\in\mathbb R^2\mid
s_1+s_2+1>0,\quad 2s_1+3s_2+1>0\},
\]
where $\underline D=(\divi(xy^2),\divi(xy^3))$. Thus
$H=\{s_1+s_2+1=0\}$ is a remarkable hyperplane of $\underline D$. For any $(a,b)\in\mathbb Z_{\geq0}^2\setminus\{(0,0)\}$, put
\[
h=(xy^2)^a(xy^3)^b=x^uy^v,\qquad
u=a+b,\quad v=2a+3b>u.
\]
The only remarkable number of $h$ is hence $1/v$.
Thus no nonzero nonnegative integral specialization of $H$ gives the
negative of a remarkable number.
\end{example}

	\clearpage
	\appendix
	\patchcmd{\subsubsection}{\itshape}{\bfseries\boldmath}{}{\PackageError{appendix}{Could not set the subsubsection title font}{Check the amsart heading definition.}}
	\section{Strong monodromy conjecture for remarkable hyperplanes}\label{app:smc}
	\begin{center}
		\scshape by S. Kurokboev, Q. Shi,  H. Zuo
	\end{center}
	
	\medskip

	\subsection{Main results} The following were proved with AI, and verified and arranged by the authors of this Appendix.

	\begin{theorem}\label{MC_rmk} (SMC for remarkable hyperplanes)
			With notation as in \ref{Rmk_subs}, suppose in addition that $D_i=\operatorname{div}(f_i)$ for nonzero regular functions $f_i$, and put $F=(f_1,\ldots,f_r)$. If $H$ is a remarkable hyperplane, its complexification $H_{\mathbb C}$ is contained in the zero locus of $B_F$, the Bernstein--Sato ideal of $F$, and $\mathrm{Exp}(H_{\mathbb C})$ is contained in the monodromy support locus of $F$. Here $\mathrm{Exp}(s_1,\ldots,s_r):=(e^{2\pi\sqrt{-1}s_1},\ldots,e^{2\pi\sqrt{-1}s_r})$.
	\end{theorem}
	By~\cite{BVWZ}, the Bernstein--Sato conclusion implies the monodromy support conclusion. For $r=1$, the theorem recovers the strong monodromy conjecture for remarkable numbers \cite[Theorem 1.1]{BSZ_SMC}. 
	
	The proof has two steps. Firstly, every remarkable hyperplane is witnessed by a boundary divisor of the log canonical model of $\bigl(X,(D_1+\cdots+D_r)_{\red}\bigr)$. This is a straightforward generalization of the similar assertion for remarkable numbers \cite[Theorem 1.2]{BSZ_SMC}. 
	
		Secondly, each such divisor defines a hyperplane in the zero locus of $B_F$. 
		The proof of this step for $r=1$ differs from the proof of the similar assertion for remarkable numbers \cite[Theorem 1.3]{BSZ_SMC}; it uses recent results of Davis and Yang \cite{DY} on specializations of multivariate $V$-filtrations together with an appropriate choice of an auxiliary function. %This proof, together with Theorem \ref{thm_topological_type},  justify partly an older hope that guided since \cite{Bud15} the analysis of singularity invariants associated to a tuple $F$, that by using appropriate auxiliary functions one can address the monodromy conjecture for the case of a single polynomial $f$.
		The same argument also applies to the monomial specializations
		in Theorem~\ref{ThmC} that preserve the reduced zero locus.
		\begin{theorem}\label{thm:specialized-BS}
			Under the assumptions of Theorem~\ref{MC_rmk}, let
			$M=(m_{ij})\in\mathbb Z_{\geq 0}^{p\times r}$ and put
			\begin{displaymath}
				F^M=(g_1,\ldots,g_p),\qquad g_i=\prod_{j=1}^r f_j^{m_{ij}}.
			\end{displaymath}
			Assume that $V(g_1\cdots g_p)=V(f_1\cdots f_r)$, where $V(-)$ denotes the zero locus.
			For every $H\in\mathrm{Rmk}(\underline D)$, the complexification of
			$\sigma_M^{-1}(H)$ is a hyperplane contained in $V(B_{F^M})$.
			%For every closed point $x\in X$ and every
			%$H\in\mathrm{Rmk}_x(\underline D)$, the complexification of
			%$\sigma_M^{-1}(H)$ is an irreducible component of $V(B_{F^M,x})$,
			%where $B_{F^M,x}$ denotes the local Bernstein--Sato ideal of $F^M$ at $x$ \cite{BVWZ}.
		\end{theorem}
		\begin{remark}\label{rmk:diagonal-specialization}
			Let $f=f_1\cdots f_r$. If a remarkable hyperplane of $F$ contains
			$(s_0,\ldots,s_0)$, then Theorem~\ref{thm:specialized-BS} with $M=(1,\ldots,1)$
			gives $b_{f}(s_0)=0$. In the multibranch plane curve examples of
			\cite[Example~5.8]{BPSZ}, summing the branch multiplicities destroys
			divisibility relations, so the negatives of some poles contributed by rupture
			divisors, cf. Section \ref{Sec_curves}, are not remarkable numbers of $f$. By Theorem~\ref{ThmD},
			every rupture divisor gives a remarkable hyperplane of the complete splitting,
			so Theorem~\ref{thm:specialized-BS} recovers these poles as roots of $b_{f}$.
		\end{remark}

	\subsection{Auxiliary theorems}\label{subs_APA}
		By Remark~\ref{rmk_local} and \cite[Theorem~1.5.1(i)]{BVWZ}, to prove Theorem \ref{MC_rmk}, it suffices
		to show that every remarkable hyperplane at a point $x$ is contained in
		the zero locus of the local Bernstein--Sato ideal at $x$.
		This will follow from Theorems~\ref{thm:witness} and~\ref{thm:lcBF}.

	Let $X$ be a smooth irreducible complex variety, let $F=(f_1,\ldots,f_r)$
	be a tuple of regular functions, and assume that $f=f_1\cdots f_r$ is
	nonzero and non-invertible. Set $D=(\Div_X f)_{\red}$.
	For a prime divisor $E$ over $X$, put
	\begin{displaymath}
		\bm N_E=(\ord_E f_1,\ldots,\ord_E f_r),\qquad
		\nu_E=A_X(E),\qquad \ell_E(\bs)=\bm N_E\cdot\bs+\nu_E,
	\end{displaymath}
	where $\bs=(s_1,\ldots,s_r)$ and $A_X(E)$ is the log discrepancy of $E$
	over $X$. For the divisorial valuation $v=\ord_E$ we also write $N_{vj}=v(f_j)$,
	$\bm N_v=\bm N_E$, $\nu_v=\nu_E$, and $\ell_v=\ell_E$.
	
	Fix a closed point $x\in X$. Let $\mu:Y\to X$ be a log resolution
	which is an isomorphism off $D$. Write $(\Div_Y f)_{\red}=\bigcup_{i\in J}E_i$
	and $E_I=\bigcap_{i\in I}E_i$. For a full-rank sublattice $L\subseteq\Z^r$
	and $\bm a\in\Z^r$, put
	\begin{align*}
		J_x&=\bigcup\left\{I\subseteq J \,\middle|\,
		I\ne\varnothing,\ x\in\mu(E_I),\
		\bm a\in L+\sum_{i\in I}\Z\bm N_{E_i}\right\},\\
		P_x&=\bigcap_{i\in J_x}\{\bs\in\R^r\mid\ell_{E_i}(\bs)\ge0\}.
	\end{align*}
	An empty intersection of half-spaces is $\R^r$. For fixed $(\bm a,L)$ the
	polyhedron $P_x$ is independent of the choice of log resolution
	\textcolor{blue}{by} Theorem~\ref{intrinsic_cl} and Remark~\ref{rmk_local}. A \emph{remarkable
		hyperplane at $x$} is the complex affine span of a facet of $P_x$ for
	some $(\bm a,L)$.

	Let $c:(Z,B_Z)\to(X,D)$ be the log canonical model over $X$, where
	$B_Z=(\Div_Z f)_{\red}$ \cite[Definition~2.1 and Theorem~1.1]{OX}.
	To each prime component $Q$ of $B_Z$ we associate the hyperplane
	$H_Q=\{\bs\in\C^r\mid\ell_Q(\bs)=0\}$.

	\begin{theorem}\label{thm:witness}
	  For every remarkable hyperplane $H$ at $x$, there is a prime component
	  $Q$ of $B_Z$ such that $x\in c(Q)$ and $H=H_Q$.
	\end{theorem}
	
	For the Bernstein--Sato ideal, we work on $X^{\mathrm{an}}$. Set
	$S=\C[\bs] = \mathbb C[s_1,\ldots,s_r]$ and let $\mathscr D_{X,x}$ be the ring of germs of holomorphic differential
	operators at $x$. Write $F^{\bs}=\prod_{j=1}^r f_j^{s_j}$ for the formal
	generator of $\mathcal O_{X,x}[1/f][\bs]F^{\bs}$. Put
	\begin{displaymath}
		\mathcal N_{F,x}
		=\frac{\mathscr D_{X,x}[\bs]F^{\bs}}{\mathscr D_{X,x}[\bs]\,fF^{\bs}},
		\qquad \gamma_x=[F^{\bs}].
	\end{displaymath}
	The local simultaneous Bernstein--Sato ideal is
	$B_{F,x}=\operatorname{Ann}_S(\gamma_x)$. It is nonzero by \cite{Sab87}. For the regular tuple $F$ the algebraic
	and analytic germs give the same ideal.
	Write $V(B_{F,x})\subseteq\C^r$ for the zero locus of $B_{F,x}$.
	
	\begin{theorem}\label{thm:lcBF}
		For every prime component $Q$ of $B_Z$ and every closed point
		$x\in c(Q)$, the hyperplane $H_Q$ is an irreducible component of
		$V(B_{F,x})$. In particular, for every closed point $x\in X$, each
		remarkable hyperplane at $x$ is an irreducible component of $V(B_{F,x})$.
	\end{theorem}

	In proving Theorem~\ref{thm:lcBF}, we obtain a specialization result for the
	enlarged lct polytope. Let $g$ be a nonconstant regular function on $X$
	such that $V(f)\subseteq V(g)$. For the tuple $(F,g)$, choose a log resolution of $\Div_X(fg)$
	dominating the one fixed above, and denote it again by $\mu:Y\to X$. Following Remark~\ref{rmk:enlarged-lct}, put
	\begin{displaymath}
		\mathrm{lct}_x(F,g)=
		\bigcap_{\substack{E\subset(\Div_Y(fg))_{\red}\\ x\in\mu(E)}}
		\{(\bs,t)\in\R^{r+1}\mid
		    \bm N_E\cdot\bs+(\ord_E g)t+\nu_E>0\},
	\end{displaymath}
	where $E$ runs over the prime components indicated. The negative of
	$\mathrm{lct}_x(F,g)$ is called the enlarged lct polytope at $x$.
	By a facet of $\mathrm{lct}_x(F,g)$ we mean a facet of its closure.
	Write $\mathrm{lctf}_x(F,g)$ for the real affine hyperplanes spanned by these facets.

	\begin{theorem}\label{thm:facet-specialization}
		Let $H=\{\bm a\cdot\bs+bt+c=0\}\in\mathrm{lctf}_x(F,g)$.
		If the facet of $\mathrm{lct}_x(F,g)$ supported by $H$ meets $\{t>0\}$, then
		$\{\bs\in\C^r\mid\bm a\cdot\bs+c=0\}$ is an irreducible component
		of $V(B_{F,x})$. If $V(f)=V(g)$ in a neighbourhood of $x$, the
		conclusion holds for every $H\in\mathrm{lctf}_x(F,g)$.
	\end{theorem}

	\subsection{Witness of remarkable hyperplanes on the lc model}
	Fix $(\bm a,L)$. It suffices to prove Theorem~\ref{thm:witness} for the complex affine span of each facet of $P_x$.
	
	For a normal pair $(W,\Gamma)$ with $K_W+\Gamma$ $\Q$-Cartier, the log discrepancy
	of a prime divisor $G$ on a birational model $q:W'\to W$ is
	$A_{W,\Gamma}(G)=1+\ord_G(K_{W'}-q^*(K_W+\Gamma))$,
	with compatible canonical divisors. The pair is lc if all its log discrepancies are nonnegative.
	
	Write $\Theta=\Z^r/L$ and $\sigma=\bm a+L$. The inclusion
		$L\subseteq\Z^r$ of cocharacter lattices defines an isogeny of complex
		tori $T_L\to(\C^*)^r$. We identify its deck group with $\Theta$ by
		monodromy of the positive coordinate loops. Pull this cover back by
		$F$ over $X\setminus D$. For each normal birational model $M\to X$,
		let $\rho_M:\widetilde M\to M$ be the normalization in the full generic
		algebra of this pullback, including all field factors. This is a finite
		map with a $\Theta$-action. A proper birational morphism $M'\to M$
		induces a $\Theta$-equivariant morphism $\widetilde M'\to\widetilde M$
		by normalization. It is proper, since $\widetilde M'$ is proper over
		$M$ and $\widetilde M$ is separated over $M$, and is birational on
		corresponding components.
	
	For a nonempty $I\subseteq J$ with $E_I\ne\varnothing$, we call $I$
	and the corresponding intersection $E_I$ \emph{$(\bm a,L)$-admissible} if $
	\bm a\in L+\sum_{i\in I}\Z\bm N_{E_i}$.
	A closed point of $(\Div_Y f)_{\red}$ is \emph{$(\bm a,L)$-admissible}
	if the indices of all boundary components through it form an admissible set.
	Thus $i\in J_x$ precisely when $E_i$ contains an admissible point over $x$.
	
	\begin{lemma}[The stabilizer over a crossing]\label{lem:stabilizer}
		Let $y\in Y$ be a closed point and put $I(y)=\{i\in J\mid y\in E_i\}$.
		The fiber $\rho_Y^{-1}(y)$ is a single $\Theta$-orbit.
		The stabilizer of every point of this fiber is the subgroup
		\begin{displaymath}
			H_y=\bigl(L+\textstyle\sum_{i\in I(y)}\Z\bm N_{E_i}\bigr)/L.
		\end{displaymath}
		Consequently, $i\in J_x$ if and only if some $y\in E_i\cap\mu^{-1}(x)$
		has a $\sigma$-fixed lift to $\widetilde Y$.
	\end{lemma}
	\begin{proof}
		We work in the analytic topology on $Y$ and $\widetilde Y$. Take a
		coordinate polydisc $U$ centered at $y$ in which the only boundary
		components are the coordinate hyperplanes $E_i\cap U$, $i\in I(y)$,
		and put $U^\circ=U\setminus\bigcup_{i\in I(y)}E_i$. Choose pairwise
		disjoint neighborhoods of the points of $\rho_Y^{-1}(y)$. By properness
		we may shrink $U$ so that $\rho_Y^{-1}(U)$ lies in their union.
		
		Each connected component $W$ of $\rho_Y^{-1}(U)$ thus contains at most
		one point over $y$. By our construction of $\widetilde Y$, the inverse image of each $E_i\cap U$ is a closed
		$\Theta$-invariant analytic subset. Since $W$ is normal, removing these
		divisors leaves $W^\circ=W\cap\rho_Y^{-1}(U^\circ)$ connected. Hence
		$W^\circ$ is a connected component of the covering over $U^\circ$ and
		maps onto $U^\circ$. As $\rho_Y(W)$ is closed in $U$ by properness, we have $y\in \rho_Y(W)$. So $W$ contains exactly one point $\tilde y$ over $y$.
		
		The action of $\Theta$ on each fiber over $U^\circ$ is simply transitive,
		so it permutes the components $W^\circ$ transitively, and hence also the
		points over $y$. The stabilizer of $\tilde y$ is thus the subgroup of $\Theta$ preserving
		$W^\circ$. Its restriction to $W^\circ$ is the full deck group of
		$W^\circ\to U^\circ$, since it acts transitively on each fiber of $\rho_Y|_{W^\circ}$ and a
		deck transformation of a connected covering is determined by the image
		of one point.
		
		Two points of a fiber over $U^\circ$ lie in the same component exactly
		when a lift of a loop joins them. Thus this subgroup consists of the
		translations obtained by monodromy. The fundamental group of
		$U^\circ$ is generated by positive meridians around the $E_i$, $i\in I(y)$.
		On a meridian around $E_i$ the winding number of $f_j\circ\mu$ is $\ord_{E_i} f_j$,
		so its action is translation by $\bm N_{E_i}+L$. The stabilizer is therefore
		generated by these classes, which gives $H_y$ as claimed. For
		$y\in E_i\cap\mu^{-1}(x)$, a lift is $\sigma$-fixed exactly when
		$\sigma\in H_y$, equivalently when $I(y)$ is $(\bm a,L)$-admissible.
		The characterization of $J_x$ now gives the last assertion.
	\end{proof}

	After shrinking $X$ to an affine neighbourhood of $x$, choose, as in
	\cite[proof of Theorem~1.2]{BSZ_SMC}, a small projective
	$\Q$-factorialization $r:T\to Z$ such that the composite $p=c\circ r$
	is an isomorphism over $X\setminus D$.
	The components $B_j$ of $B_T=(\Div_T f)_{\red}$ are the strict transforms of
	the components of $B_Z$. Moreover, since $K_T+B_T=r^*(K_Z+B_Z)$, the pair $(T,B_T)$ is lc.
	Recall that $\rho_T:\widetilde T\to T$ is the normalized $\Theta$-cover
		induced by the torus isogeny $T_L\to(\C^*)^r$ and the tuple $F\circ p$.
	
	For a prime divisor $E$ over $T$, pulling back the divisors $\Div_T f_k$ and
	$K_T+B_T-p^*K_X=\sum_j\nu_{B_j}B_j$ gives
	\begin{equation}\label{eq:row}
		\ell_E(\bm s)=A_{T,B_T}(E)+\sum_j\ord_E(B_j)\ell_{B_j}(\bm s).
	\end{equation}
	Here $\ord_E(B_j)\ge0$ since each $B_j$ is effective and $\Q$-Cartier.
	Let $\mathcal A_x(T)$ consist of those $B_j$ for which there exists a
	$\sigma$-fixed point $\tilde z\in\widetilde T$ satisfying
	$\rho_T(\tilde z)\in B_j$ and $p(\rho_T(\tilde z))=x$.
	
	\begin{lemma}\label{lem:push}
		The inequalities $\ell_{B_j}(\bm s)\ge0$ for $B_j\in\mathcal A_x(T)$
		imply every defining inequality of $P_x$.
	\end{lemma}
	
	\begin{proof}
		By Theorem~\ref{intrinsic_cl} and Remark~\ref{rmk_local}, we can compute $P_x$ on any log resolution over an affine neighbourhood of $x$.
		Choose a log resolution $Y\to X$ factoring through $T$ that is an isomorphism over $X\setminus D$. For each $i\in J_x$, Lemma~\ref{lem:stabilizer} gives a point
		$y\in E_i \subset Y$ over $x$ with a $\sigma$-fixed lift to $\widetilde Y$. Its image $z\in T$ has a fixed
		lift to $\widetilde T$ by equivariance. If $\ord_{E_i}(B_j)>0$, the center of $E_i$ on $T$
		lies in $B_j$, so $z\in B_j$ and $B_j\in\mathcal A_x(T)$.
		Now we conclude the proof by applying \eqref{eq:row} and $A_{T,B_T}(E_i)\ge0$.
	\end{proof}
	
	The following lemma follows from the argument for the reverse inequality in
	\cite[proof of Theorem~1.2]{BSZ_SMC}, replacing the $\zeta^d$-fixed points on the
	normalized cyclic cover by the $\sigma$-fixed points on the normalized torus cover.
	For $B\in\mathcal A_x(T)$, the divisor $G$ obtained in that argument satisfies
	$A_{T,B_T}(G)=0$, $\ord_G(B)>0$, and $\ord_G(B_j)=0$ for $B_j\ne B$.
	Hence \eqref{eq:row} gives $\ell_G=\ord_G(B)\ell_B$.
	
	\begin{lemma}\label{lem:pointed}
		For every $B\in\mathcal A_x(T)$ there is a log resolution
		$Y_B\to X$ and a component $G=E_i$ on $Y_B$ with $i\in J_x$
		such that $\ell_G=d_G\ell_B$ for some $d_G\in\Q_{>0}$.
	\end{lemma}
	
	\begin{proposition}\label{prop:facets}
		With the notation above,
		\begin{equation}\label{eq:lcwit}
			P_x=\bigcap_{B\in\mathcal A_x(T)}
			\{\bm s\in\R^r\mid\ell_B(\bm s)\ge0\}.
		\end{equation}
	\end{proposition}
	
	\begin{proof}
		Lemma~\ref{lem:push} gives the inclusion from right to left. For
		$B\in\mathcal A_x(T)$ compute $P_x$ on the resolution $Y_B$ of
		Lemma~\ref{lem:pointed}. There $\ell_G\ge0$ is a defining inequality
		of $P_x$, and $\ell_G=d_G\ell_B$ with $d_G>0$. This gives the reverse
		inclusion.
	\end{proof}
	
	\begin{proof}[Proof of Theorem~\ref{thm:witness}]
		By Proposition~\ref{prop:facets}, each facet of $P_x$ is again a facet of the right-hand side of \eqref{eq:lcwit}, which is always given by $\{\ell_Q(\bm s) = 0\}$ for some prime divisor $Q$ of $B_Z$. This proves our assertion. 
	\end{proof}

	\subsection{Boundary divisors of the lc model and the Bernstein--Sato ideal}
	\label{sec:bern} In this subsection, we will address Theorem \ref{thm:lcBF}.

	\subsubsection{Multivariate \texorpdfstring{$V$}{V}-filtration}
		We use the conventions of Davis and Yang \cite{DY} for the tuple $q = (q_1,\ldots,q_d) := (f_1,\ldots,f_r,g)$, where $g$ is any nonzero regular
		function on an affine neighbourhood of $x$. Put $d=r+1$,
		$D_q=\operatorname{div}(fg)$, $\bm z=(\bs,t)$, $R=S[t]$, and
		$e=F^{\bs} g^t$. Consider the localized graph module
		\begin{displaymath}
			\mathcal G=\mathcal O_X[1/(fg)][s_1,\ldots,s_r,t]\,e .
		\end{displaymath}
	In this presentation $s_j$ and $t$ act by polynomial multiplication.
	For a local derivation $\partial$ on $X$ and a local section
	$a\in\mathcal O_X[1/(fg)][\bm z]$, the differential action is
	\begin{equation}\label{eq:lcbs-action}
		\partial(ae)=
		\left(\partial a+
		a\sum_{j=1}^r s_j\frac{\partial f_j}{f_j}
		+at\frac{\partial g}{g}\right)e .
	\end{equation}
	Let $\iota_q:X\to X\times\mathbb C^d$ be the graph embedding, with target
	coordinates $\bm\tau=(\tau_1,\ldots,\tau_d)$. At $(x,q(x))$, sections
	of the graph direct image are finite sums of normal derivatives
	$\partial_{\bm\tau}^{\bm k}\delta_q$ of the graph generator $\delta_q$,
	with coefficients in $\mathcal O_{X,x}[1/(fg)]$.
	Each $q_j$ is invertible in this ring, so the Malgrange--Mellin transform
	\cite[Proposition 2.52]{DY} identifies the following $\mathscr D_X[\bm z]$-modules on $X$ via the graph projection
	\begin{displaymath}
		\iota_{q,+}\bigl(\mathcal O_X[1/(fg)]\bigr)\ \simeq\
		\mathcal O_X[1/(fg)][\bm z]e,\qquad
		\delta_q\longmapsto e,\qquad z_j=-\partial_{\tau_j}\tau_j .
	\end{displaymath}
	In particular $t=z_d=-\partial_{\tau_d}\tau_d$ acts as a polynomial
	variable. This identifies $\iota_{q,+}\mathcal O_X$ with a submodule
	of $\mathcal G$ stable under $\mathscr D_X[\bm z]$ and the $\tau_j$.
	Writing $\bm e_j$ for the $j$-th coordinate vector, we have
	$\tau_j(a(\bm z)e)=q_j a(\bm z+\bm e_j)e$, so $\tau_d e=ge$.

	The \textit{localized $V$-filtration} $V_*^\bullet\mathcal G$ is an exhaustive, decreasing,
	$\mathbb R^d$-indexed, \textit{wall and chamber filtration}, see below, by $\mathscr D_X[\bm z]$-submodules of $\mathcal G$
	\cite[Definition 2.34 and Corollary 2.38]{DY}. Thus
	$V_*^{\bm b}\mathcal G\subseteq V_*^{\bm a}\mathcal G$ whenever
	$a_j\le b_j$ for all $j$, which we write as $\bm a\le\bm b$.
	\begin{definition}[{\cite[Definitions 2.2--2.5]{DY}}]
		A \emph{wall and chamber filtration} is a filtration for which, locally
		on $X$, one can choose a collection $\mathcal W$ of affine hyperplanes
		in $\mathbb R^d$, finite modulo translation by $\mathbb Z^d$.
		Each wall $W\in\mathcal W$ has an equation $\bm L\cdot\bm c=\gamma$
		with $0\ne\bm L\in\mathbb Q_{\ge0}^d$ and $\gamma\in\mathbb R$.
		Its two sides are
		\begin{displaymath}
			W^-=\{\bm c \mid \bm L\cdot\bm c\le\gamma\},\qquad
			W^+=\{\bm c \mid \bm L\cdot\bm c>\gamma\}.
		\end{displaymath}
		The half-spaces $\{W^-,W^+ \mid W \in\mathcal  
		W\}$ divide $\mathbb R^d$ into \emph{chambers}. The filtration is
		constant on each chamber.
	\end{definition}
	
	Two indices $\bm a\le\bm b$ are \emph{separated by a single wall} $W = \{\bm c \mid \bm L \cdot \bm c = \gamma\}$
	if $\bm L\cdot\bm a\le\gamma<\bm L\cdot\bm b$ for $W$ and the two
	indices lie on the same side of every other wall. For such a pair $(\bm a,\bm b)$,  \cite[Definition 2.34(2), Proposition 3.12 and Corollary 5.8]{DY}
	give an integer $m\ge1$ locally such that
	\begin{equation}\label{eq:nilpotent}
		(\bm L\cdot\bm z+\gamma)^m
		\bigl(V_*^{\bm a}\mathcal G/V_*^{\bm b}\mathcal G\bigr)=0 .
	\end{equation}
	
	Let $V^\bullet\iota_{q,+}\mathcal O_X$ denote the ordinary multivariate
	$V$-filtration on the graph module along the coordinate hyperplanes
	$\tau_j=0$ \cite[Definition 2.16 and Theorem 3.6]{DY}.
	It is also an exhaustive, decreasing, $\mathbb R^d$-indexed wall and chamber
	filtration by $\mathscr D_X[\bm z]$-submodules.
	By \cite[Corollary 2.38]{DY}, under the inclusion
	$\iota_{q,+}\mathcal O_X\subseteq\mathcal G$ one has $
	V^{\bm c}\iota_{q,+}\mathcal O_X=V_*^{\bm c}\mathcal G$ {for every }$\bm c\in\mathbb R_{>0}^d$.

	Taking stalks at $x$, consider $L_{\bm c} := V_*^{\bm c}\mathcal G_x$ for $\bm c=(c_1,\ldots,c_d)\in\mathbb Q^d$.
	\begin{corollary}[{\cite[Corollary 2.38 and its proof]{DY}}]
		For $\bm n\in\mathbb Z_{\ge0}^d$ such that $\bm c+\bm n>\bm 0$ in every coordinate, one has
		\begin{equation}\label{eq:lcbs-shift}
			L_{\bm c}=\bm{\tau}^{-\bm n}\bigl(V^{\bm c+\bm n}\iota_{q,+}\mathcal O_X\bigr)_x,\qquad
			\bm{\tau}^{-\bm n}=\prod_{j=1}^d\tau_j^{-n_j}.
		\end{equation}
	\end{corollary}
	
	\subsubsection{Hodge filtration}
	On $\mathcal G$ the localized Hodge filtration is given by total degree
	in $\bm z$. With $e$ in degree zero, its $p$-th piece is
	\begin{displaymath}
		H_p\mathcal G=\mathcal O_X[1/(fg)][\bm z]_{\le p}\,e
		\qquad (H_p\mathcal G=0 \text{ for }p<0).
	\end{displaymath}
	The ordinary Hodge filtration on $\iota_{q,+}\mathcal O_X$ is its
	restriction, $H_p(\iota_{q,+}\mathcal O_X)=(\iota_{q,+}\mathcal O_X)\cap H_p\mathcal G$. In the notation of Davis and Yang,
	$H_p\mathcal G=F_p^*\mathcal G$. Under the Malgrange--Mellin transform this piece
	corresponds to $F^*_{p+d}\bigl(\iota_{q,+}(\mathcal O_X[1/(fg)])\bigr)$
	on the localized graph module. On the ordinary graph module,
	$H_p(\iota_{q,+}\mathcal O_X)=F_{p+d}(\iota_{q,+}\mathcal O_X)$
	\cite[\S5.1 and Proposition 5.28]{DY}. At $x$ we put $H_pL_{\bm c}=L_{\bm c}\cap(H_p\mathcal G)_x$ and $\mathrm{gr}_p^H L_{\bm c} = H_{p}L_{\bm c}/H_{p-1} L_{\bm c}$. In particular, $\operatorname{gr}^H_0L_{\bm c}=H_0L_{\bm c}$.

	Since by definition $z_j H_pL_{\bm c} \subset H_{p+1}L_{\bm c}$, it induces an action on $\operatorname{gr}^H L_{\bm c} = \bigoplus_{k\geq 0} \operatorname{gr}^H_k L_{\bm c}$ by
	\begin{equation}\label{eq:symboll}
		z_j : \operatorname{gr}^H_k L_{\bm c} \to \operatorname{gr}^H_{k+1} L_{\bm c},\ [u] \mapsto [z_ju],\ \text{for all } k \in \mathbb N.
	\end{equation}
	For $\bm a\le\bm b$, the inclusion $L_{\bm b}\hookrightarrow L_{\bm a}$
	is strict for the Hodge filtration. Indeed, both filtrations are induced from
	$H_\bullet\mathcal G_x$, so
	$L_{\bm b}\cap H_pL_{\bm a}=L_{\bm b}\cap H_p\mathcal G_x=H_pL_{\bm b}$. The quotient filtration is given by
	\begin{displaymath}
		H_p(L_{\bm a}/L_{\bm b})
		=\operatorname{im}\!\left(H_pL_{\bm a}\to L_{\bm a}/L_{\bm b}\right)
		\simeq H_pL_{\bm a}/H_pL_{\bm b}.
	\end{displaymath}
	Its associated graded is a $\mathbb C[\bm z]$-module with the action
	described in \eqref{eq:symboll}.

	Choose a log resolution $\pi:Y\to X$ of $D_q$ that is an isomorphism
	off its support. For each component $E_i$ of $(\pi^*D_q)_{\rm red}$,
	abbreviate $\bm N_{E_i}$ to $\bm N_i=(N_{ij})_{j=1}^r$, where
	$N_{ij}=\ord_{E_i}f_j$, and write $\nu_i=\nu_{E_i}=1+\ord_{E_i}K_{Y/X}$,
	where $K_{Y/X}=K_Y-\pi^*K_X$. Put $G_i=\ord_{E_i}g$ and
	$\bm M_i=(\bm N_i,G_i)$.
	
	\begin{corollary}
		For every $\bm c\in\mathbb Q^d$, the lowest Hodge piece of $L_{\bm c}$ is
		\begin{equation}\label{eq:lcbs-lowest}
			H_0L_{\bm c}=
			\left(\pi_*\mathcal O_Y\!\left(
			\sum_i(\nu_i-\lceil\bm M_i\cdot\bm c\rceil)E_i\right)\right)_x e .
		\end{equation}
	\end{corollary}
	
	\begin{proof}
		Choose $\bm n\in\mathbb Z_{\ge0}^d$ so that $\bm c+\bm n>\bm0$ in
		every coordinate, and put
		$\Delta=\sum_{j=1}^d(c_j+n_j-\delta)\operatorname{div}(q_j)$ for
		$0<\delta\ll1$. By \cite[Corollary 5.40]{DY},
		\begin{equation}\label{eq:multiplier}
			H_0\bigl(V^{\bm c+\bm n}\iota_{q,+}\mathcal O_X\bigr)_x
			=\mathcal J(X,\Delta)_x e =\left(\pi_*\mathcal O_Y\!\left(
			\sum_i\bigl(\nu_i-\lceil\bm M_i\cdot(\bm c+\bm n)\rceil\bigr)E_i
			\right)\right)_x e,
		\end{equation}
		where $\mathcal J(X,\Delta)$ denotes the multiplier ideal of $(X,\Delta)$.
		On $\mathcal J(X,\Delta)_x e$, the shift $\bm \tau^{-\bm n}$ acts by
		$\bm\tau^{-\bm n}(ae)=(\prod_{j=1}^d q_j^{-n_j})ae$.
		This adds $\bm M_i\cdot\bm n$ to the coefficient of $E_i$ on the right-hand side of \eqref{eq:multiplier}, giving \eqref{eq:lcbs-lowest}.
	\end{proof}

	Let $\bm a\le\bm b$ and let $\mathcal W$ be a set of walls for
	$V_*^\bullet$ near $x$. A wall $W=\{\bm c \mid \bm L\cdot\bm c=\gamma\}$
	in $\mathcal W$ separates $\bm a$ and $\bm b$ precisely when
	\begin{equation}\label{eq:lcbs-levels}
		\bm L\cdot\bm a\le\gamma<\bm L\cdot\bm b .
	\end{equation}
	Let $W_i=\{\bm c \mid \bm L_i\cdot\bm c=c_i\}$, $1\le i\le k$, be all
	the walls in $\mathcal W$ separating $\bm a$ and $\bm b$. By
	\cite[Remark 2.42]{DY}, one can pass from the chamber of $\bm a$ to
	that of $\bm b$ by crossing these walls one at a time in increasing
	order of the filtration indices. Applying \eqref{eq:nilpotent} to the
	successive quotients gives, for some integer $m\ge1$,
	\begin{displaymath}
		\left(\prod_{i=1}^k(\bm L_i\cdot\bm z+c_i)\right)^m
		\bigl(L_{\bm a}/L_{\bm b}\bigr)=0 .
	\end{displaymath}
	
	\begin{definition}[{\cite[Definitions 4.18 and 5.34]{DY}}]\label{def:lcbs-relative}
		Let $W=\{\bm c \mid \bm L\cdot\bm c=\gamma\}$. An $R$-module is
		\emph{flat} (resp.\ \emph{free}) \emph{relative to $W$} if it is flat (resp.\ free) over $\mathbb C[U]$ for
		every subspace $U\subseteq\operatorname{span}_{\mathbb C}\{z_1,\ldots,z_d\}$
		with $\bm L\cdot\bm z\notin U$. Here $\mathbb C[U]\subseteq R$ is
		the polynomial algebra generated by $U$.
	\end{definition}
	
	\begin{theorem}[{\cite[Theorems 4.17 and 5.33]{DY}}]\label{thm:lcbs-free}
		\textup{(1)} For every $\bm c \in \mathbb Q^d$,  both $L_{\bm c}$ and
		$\operatorname{gr}^H L_{\bm c} $ are free over $R$, where $z_j$ acts on $\operatorname{gr}^H L_{\bm c}$ as in \eqref{eq:symboll}.
		
		\textup{(2)} If $\bm a\le\bm b$ are separated by a single wall $W$, then
		$L_{\bm a}/L_{\bm b}$ (resp.\ $\operatorname{gr}^H(L_{\bm a}/L_{\bm b})$)
		is flat (resp.\ free) relative to $W$, where the quotient carries the
		quotient Hodge filtration.
	\end{theorem}

	\subsection{Proof of Theorem \ref{thm:lcBF}}
	Fix a prime component $Q$ of $B_Z$ and $x\in c(Q)$. We write $\bm{N}_Q,\nu_Q$ for the data of
	$\operatorname{ord}_Q$. Then $\bm{N}_Q\ne\bm 0$, since $\operatorname{ord}_Q f>0$.  For the local argument we may replace $X$
	by an affine neighbourhood of $x$ trivializing $K_X$. Write
	$\bm 1_k \in \mathbb N^k$ for the vector of $k$ ones.
	
	\begin{lemma}\label{lem:lcbs-wall}
		There exist a nonzero $g\in f\mathcal O_X(X)$ and an index
		$\bm{\alpha}^0\in\mathbb Q^{r+1}$ with $\alpha^0_{r+1}<0$ such that
		\begin{displaymath}
			(\bm{N}_v,v(g))\cdot\bm{\alpha}^0\le\nu_v
		\end{displaymath}
		for every normalized divisorial valuation $v$ over $X$, with equality
		precisely at $v=\operatorname{ord}_Q$.
	\end{lemma}
	
	\begin{proof}
		Choose compatible canonical divisors and write
		\begin{displaymath}
			P=K_Z+B_Z-c^*K_X=\sum_{D\subset B_Z}\nu_D D,\qquad H_Z=\operatorname{div}_Z f,
		\end{displaymath}
		where the sum runs over the prime components of $B_Z$. Since $P$ is
		$c$-ample, it separates $Q$ as follows. Choose a smooth open set
		$U\subset Z$ meeting $Q$ on which $B_Z|_U=Q|_U$ is smooth.
		For a sufficiently divisible integer $a\gg0$, the ideal of $Z\setminus U$
		twisted by $\mathcal O_Z(aP)$ is generated by global sections over the affine
		base $X$. A section nonzero at the generic point of $Q$ gives a rational function
		$\phi \in \mathbb C(Z) = \mathbb C(X)$ such that $C=\operatorname{div}_Z\phi+aP$ is effective Cartier,
		its support contains $Z\setminus U$, and $Q$ is not a component of $C$.
		Take $b>0$ integral with $(b-1)H_Z\ge aP$. Then $g=f^b\phi$ belongs
		to $f\mathcal O_X(X)$, since $g/f$ is regular on $Z$ and
		$c_*\mathcal O_Z=\mathcal O_X$.
		
		For each divisorial valuation $v$, put $\bm{M}_v=(\bm{N}_v,v(g))$. At the index
		$\bm{\alpha}^0=((b/a)\bm 1_r,-1/a)$, the defining equation for $C$ gives
		\begin{equation}\label{eq:lcbs-discrepancy}
			\nu_v-\bm{M}_v\cdot\bm{\alpha}^0
			=A_{(Z,B_Z)}(v)+v(C)/a\ \ge 0 .
		\end{equation}
		Here $A_{(Z,B_Z)}(v)=\nu_v-v(P)$ is the log discrepancy of the pair and
		$v(C)$ is the order of a local equation of $C$. Equality holds at $Q$.
		Conversely, equality forces $v(C)=0$, so the generic point of the centre
		of $v$ lies in $U$. Since a smooth pair with one smooth reduced prime boundary is plt, that boundary is the unique divisor with log discrepancy zero. This proves
		uniqueness among normalized valuations.
	\end{proof}
	
	\vspace{0.5em}

	For the function $g$ and the index $\bm\alpha^0$ of
	Lemma~\ref{lem:lcbs-wall}, choose a log resolution of $\Div_X(fg)$
	dominating $Z$ and denote the strict transform of $Q$ again by $Q$.
	At $-\bm\alpha^0$ the inequality defining
	$\overline{\mathrm{lct}_x(F,g)}$ is an equality for $Q$ and is strict
	for every other component over $x$. Hence $-\bm\alpha^0$ lies in the
	relative interior of the facet supported by
	$\bm N_Q\cdot\bs+(\ord_Q g)t+\nu_Q=0$. By Lemma~\ref{lem:lcbs-wall}, the last coordinate of $-\bm\alpha^0$ is positive.\par

	\begin{proof}[Proof of Theorem~\ref{thm:lcBF}]
		Assuming Theorem~\ref{thm:facet-specialization}, we specialize the
		preceding facet hyperplane at $t=0$. This gives
		$H_Q=\{\bs\in\C^r\mid\bm N_Q\cdot\bs+\nu_Q=0\}$ as an irreducible
		component of $V(B_{F,x})$. The assertion for remarkable hyperplanes
		follows from Theorem~\ref{thm:witness}.
	\end{proof}

	\subsubsection{Proof of Theorem~\ref{thm:facet-specialization}}
	Now take $g$ and $H$ as in the statement. Work near $x$ on the log
	resolution in \eqref{eq:lcbs-lowest}, discarding components whose
	images miss $x$. Choose a component $Q$ whose inequality defines the
	facet of $\mathrm{lct}_x(F,g)$ supported by $H$, and put $G_Q=\ord_Q g$, $\bm M_Q=(\bm N_Q,G_Q)$ and
	$\ell(\bs,t)=\bm M_Q\cdot\bm z+\nu_Q$. Choose $-\bm\alpha$ in the
	relative interior of this facet. If the facet meets $\{t>0\}$, choose
	$\alpha_d<0$. In this case $\bm N_Q\ne\bm0$, since $G_Qt+\nu_Q>0$ for $t>0$.
	If $V(f)=V(g)$ near $x$, the reduced total transforms of $fg$ and $f$
	coincide there. Thus $\ord_{E_i}f=\sum_j N_{ij}>0$ for every $E_i$
	under consideration. In particular every $\bm N_i$ is nonzero, and
	$\bm\alpha$ may be chosen without a restriction on $\alpha_d$.
	In either case we can take $\bm\alpha$ rational and away from every
	hyperplane $\bm M_i\cdot\bm c=k$, $k\in\Z$, distinct from
	$\bm M_Q\cdot\bm c=\nu_Q$, by local finiteness.
	Set $\bm\beta=\bm\alpha+\epsilon\bm 1_d$ for a sufficiently small
	rational $\epsilon>0$. The only wall separating $\bm\alpha$ and $\bm\beta$
	is $\{\bm c\mid\bm M_Q\cdot\bm c=\nu_Q\}$. In the coordinates
	$(\bs,t)=-\bm c$, this wall is
	$H=\{(\bs,t)\mid\bm N_Q\cdot\bs+G_Qt+\nu_Q=0\}$,
	even if several components define it. Take $\epsilon$ small enough that
	$\lceil\bm M_i\cdot\bm\beta\rceil\le\nu_i+1$ for every $i$ and that
	$\beta_d<0$ whenever $\alpha_d<0$. Write $A_-=L_{\bm\alpha}$, $A_+=L_{\bm\beta}$,
	$\overline A=A_-/A_+$ and $\mathcal K=\overline A/t\overline A$.\par
	
	\begin{lemma}\label{lem:lcbs-specialize}
		The generator $e$ belongs to $H_0 A_-$, and its image $e_0$ in
		$\mathcal K$ is nonzero. Moreover $fe\in A_+$, and $\ell(\bs,t)^m \overline{A}=0$
		for some $m\ge1$.
	\end{lemma}
	
	\begin{proof}
		By \cite[Corollary 3.7 and Example 3.8]{DY}, we may take as a set of
		walls for $V_*^\bullet$ the hyperplanes $\bm M_i\cdot\bm c=k$,
		$k\in\mathbb Z$, associated to the divisors $E_i$ on our log
		resolution. By \eqref{eq:lcbs-levels} and the choice of $\bm\alpha,\bm\beta$,
		the only separating wall is $\bm M_Q\cdot\bm c=\nu_Q$. Thus
		\eqref{eq:nilpotent} gives $\ell(\bs,t)^m\overline{A}=0$ for some $m\ge1$. 
		
		Coefficients of all $E_i$ in \eqref{eq:lcbs-lowest} are nonnegative at $\bm c = \bm{\alpha}$. At $Q$ the coefficient drops from $0$ to $-1$ at $\bm{\beta}$. Thus a
		section at $\bm{\beta}$ must vanish along $Q$, whereas the coefficient of $e$ is allowed to be $1$ at $\bm{\alpha}$. Since $\pi(Q)$ contains
		$x$, this gives $e\in H_0A_-$ and $e\notin H_0A_+$ at the chosen
		stalk. Thus $\bar e\ne0$ in $H_0\overline{A}$. 
		
		For $fe$, recall that $\bm\beta$ was chosen so that
		$\lceil\bm M_i\cdot\bm\beta\rceil\le\nu_i+1$ for every $i$.
		Thus, if $\bm N_i\ne\bm 0$, then
		$\operatorname{ord}_{E_i}(f)=\sum_{j=1}^r N_{ij}
		\ge1\ge\lceil\bm M_i\cdot\bm\beta\rceil-\nu_i$. For $E_i$ with $\bm{N}_i=\bm 0$ we have
		$G_i>0$ and $\bm{M}_i\cdot\bm{\beta}=G_i\beta_d<0$.
		In both cases, we have
		$\operatorname{ord}_{E_i}(f)\ge\lceil\bm M_i\cdot\bm\beta\rceil-\nu_i$.
		Hence \eqref{eq:lcbs-lowest} gives $fe\in H_0 A_+$.
		
		Since $\bm{N}_Q\ne\bm 0$, $\bm{N}_Q\cdot\bs+G_Qt$ does not lie in the
		span of $t$. Apply Theorem~\ref{thm:lcbs-free}(2) with
		$U=\mathbb C t$. It gives freeness of $\operatorname{gr}^H \overline{A}$ over
		$\mathbb C[U]=\mathbb C[t]$. Thus, if $u\in H_p\overline{A}\setminus H_{p-1}\overline{A}$, the leading
		class of $tu$ is $t$ times the leading class of $u$, and is nonzero
		in degree $p+1$. Since $H_{-1}\overline{A}=0$, this gives $t\overline{A}\cap H_0\overline{A}=0$.
		Thus the nonzero class $\bar e\in H_0\overline{A}$ has nonzero image $e_0\in \mathcal K$.
	\end{proof}
	
	\begin{lemma}\label{lem:lcbs-original}
		There are natural maps
		\begin{displaymath}
			\mathcal N_{F,x}\ \longrightarrow\
			\mathcal G_x/(A_++t\mathcal G_x)\ \longleftarrow\ \mathcal K
		\end{displaymath}
		under which $\gamma_x$ and $e_0$ have the same image. The right arrow is
		injective and $\operatorname{Ann}_S(e_0)\subseteq(\ell_Q)$.
		In particular $B_{F,x}\subseteq(\ell_Q)$.
	\end{lemma}
	
	\begin{proof}
		For the right arrow, we show $A_-\cap t\mathcal G_x=tA_-$.
		Let $u\in\mathcal G_x$ be such that $tu\in A_-$. By exhaustivity,
		$u\in L_{\bm{\alpha}-j\bm 1_d}$ for some integer $j>0$.
		Using the same set of walls for $V_*^\bullet$, take a power
		$D(\bs,t)$ of the product of the factors $\bm M_i\cdot\bm z+k$ for
		the walls separating $\bm\alpha-j\bm 1_d$ and $\bm\alpha$.
		Then $D(\bs,t)u\in A_-$ by the annihilation statement following
		\eqref{eq:lcbs-levels}.

		For a factor of $D$ with $\bm{N}_i\ne\bm 0$, its specialization at $t = 0$  is a nonzero polynomial
		$\bm{N}_i\cdot\bs+k$. For a factor $\bm M_i\cdot\bm z + k$ with $\bm{N}_i=\bm 0$, its specialization at $t=0$ is the integer
		$k<G_i\alpha_d<0$  by \eqref{eq:lcbs-levels}. In particular $D(\bs,0)\ne0$.
		Now $D(\bs,t)\,tu=t(D(\bs,t)u)\in tA_-$, so $D(\bs,0)$ kills the class
		of $tu$ in $A_-/tA_-$. This module is free over $S$, since $A_-$
		is free over $S[t]$ by Theorem~\ref{thm:lcbs-free}(1).
		Thus $tu\in tA_-$, and cancelling $t$ in
		$\mathcal G_x$ gives $u\in A_-$. As $A_+\subset A_-$, we obtain
		\begin{displaymath}
			A_-\cap(A_++t\mathcal G_x)=A_++tA_- .
		\end{displaymath}
		Since $\mathcal K=A_-/(A_++tA_-)$, this proves injectivity.
		
		Modulo $t$, formula \eqref{eq:lcbs-action} gives the differential action
		on $F^{\bs}$. Thus $F^{\bs}\mapsto e\bmod t$ identifies
		$\mathcal O_{X,x}[1/f][\bs]F^{\bs}$ with a $\mathscr D_{X,x}[\bs]$-submodule
		of $\mathcal G_x/t\mathcal G_x$. Since $fe\in A_+$ and $A_+$ is
		stable under the $\mathscr D_{X,x}[\bm s]$-action, the map descends to $\mathcal N_{F,x}$.
		It sends $\gamma_x$ to the image of $e_0$ under the right arrow. By
		injectivity, every element of $B_{F,x}$ kills $e_0$.
		
		Choose $j$ with $N_{Qj}>0$ and put $S'=\mathbb C[s_i:i\ne j]$.
		Since $\sum_{i=1}^r N_{Qi}s_i+G_Qt$ has nonzero $s_j$-coefficient, $\overline A$ is flat
		over $S'[t]$ by Theorem~\ref{thm:lcbs-free}(2).
		Thus $\mathcal K=\overline A/t\overline A$ is flat, hence torsion-free, over $S'$.
		By Lemma~\ref{lem:lcbs-specialize}, $e_0\ne0$ and $\ell_Q^m e_0=0$ in $\mathcal K = \overline A / t\overline A$.
		Take $k\ge0$ maximal with $v=\ell_Q^k e_0\ne0$, so $\ell_Qv=0$.
		If $p\in\operatorname{Ann}_S(e_0)$, then $pv=0$, and hence the residue of $p$
		in $S/(\ell_Q)\simeq S'$ kills $v$. Torsion-freeness forces this residue
		to vanish, proving $\operatorname{Ann}_S(e_0)\subseteq(\ell_Q)$.
	\end{proof}

	By Lemma~\ref{lem:lcbs-original}, $B_{F,x}\subseteq(\ell_Q)$.
	Since $B_{F,x}\ne0$ and $\ell_Q(\bs)=\ell(\bs,0)$ is proportional to
	$\bm a\cdot\bs+c$, its zero hyperplane is an irreducible component
	of $V(B_{F,x})$. This proves Theorem~\ref{thm:facet-specialization}.
	\hfill$\Box$\par

	\begin{example}
	For $f=(x^2-y^3)^2-xy^5$ and $g=yf$, $\mathrm{lct}_0(f)=\tfrac5{12}$ and Theorem~\ref{thm:facet-specialization} gives $b_{f,0}(-\tfrac{11}{26})=0$ as follows.
	The dual graph of the total transform of $g$ is as follows.
	{\color{black}
	\begin{displaymath}
		\xymatrix@C=17pt@R=8pt{
			&& E_2 & \widetilde{V(f)} \\
			\widetilde{V(y)}
			& E_1 \ar[l] \ar@{-}[r]
			& E_3 \ar@{-}[u] \ar@{-}[r]
			& E_5 \ar[u] \ar@{-}[r]
			& E_4
		}
	\end{displaymath}
	}
	We depict the enlarged lct polytope in the following graph.
	\par\smallskip
	{\centering\color{black}
	\begin{tikzpicture}[font=\footnotesize,
		curve/.style={color=blue!65!black}, liney/.style={color=violet!90!black},
		zone/.style={color=orange!80!black}, ztwo/.style={color=black!60},
		zthree/.style={color=teal!85!black}, zfour/.style={color=purple!75!black},
		zfive/.style={color=red!75!black}, candidate/.style={densely dashed},
		facet/.style={line width=1.1pt}]
		% The first three blowups resolve x^2-y^3. Near the last exceptional
		% curve take x=p^3 q, y=p^2 q and w=(1-q)/q. Then
		% f=p^12 q^6 (w^2-p). Two more blowups resolve this tangency.
		% The triples below are (ord_E f, ord_E g, A_X(E)).
		% At (-41/26,1) only the E_5 form vanishes.
		\begin{scope}[x=1.05cm,y=0.82cm]
			\begin{scope}
				\clip (-5.62,-1.41) rectangle (1.59,4.66);
				\fill[blue!8] (-5.62,5.3) -- (-5.62,4.62) -- ({-19/4},{15/4})
					-- (-1,0.5) -- (0.75,-1) -- (1.59,-1) -- (1.59,5.3) -- cycle;
				\foreach \a/\b/\c/\ink in {1/1/1/curve,0/1/1/liney,
						4/5/2/zone,6/7/3/ztwo,12/14/5/zthree,
						13/15/6/zfour,26/30/11/zfive}{
					\draw[\ink,candidate] (-5.62,{(5.62*\a-\c)/\b}) --
						(1.59,{(-1.59*\a-\c)/\b});
				}
				\draw[curve,facet] (-5.62,4.62) -- ({-19/4},{15/4});
				\draw[zfive,facet] ({-19/4},{15/4}) -- (-1,0.5);
				\draw[zthree,facet] (-1,0.5) -- (0.75,-1);
				\draw[liney,facet] (0.75,-1) -- (1.59,-1);
			\end{scope}
			\draw[->,black!65] (-5.69,0) -- (1.78,0) node[right] {$s$};
			\draw[->,black!65] (0,-1.48) -- (0,4.87) node[above] {$t$};
			\fill[curve] ({-19/4},{15/4}) circle (1.4pt);
			\fill[zfive] (-1,0.5) circle (1.4pt);
			\fill[zthree] (0.75,-1) circle (1.4pt);
			\fill[zfive] ({-41/26},1) circle (1.8pt)
				node[above right,text=black,inner sep=2pt] {$(-\tfrac{41}{26},1)$};
			\draw[black!70,thin] ({-11/26+1/312},0) circle (2.5pt);
		\end{scope}
		\begin{scope}[shift={(2.20,3.75)},x=1cm,y=1cm]
			\fill[blue!8] (0,-0.09) rectangle (0.38,0.10);
			\node[anchor=west] at (0.43,0) {$\overline{\mathrm{lct}_0(f,g)}$};
			\foreach \pos/\ink/\entry in {
				-0.47/curve/{f=0:\ s+t+1=0},
				-0.84/liney/{y=0:\ t+1=0},
				-1.21/zone/{E_1:\ 4s+5t+2=0},
				-1.58/ztwo/{E_2:\ 6s+7t+3=0},
				-1.95/zthree/{E_3:\ 12s+14t+5=0},
				-2.32/zfour/{E_4:\ 13s+15t+6=0},
				-2.69/zfive/{E_5:\ 26s+30t+11=0}}{
				\draw[\ink,candidate] (0,\pos) -- (0.38,\pos);
				\node[anchor=west] at (0.43,\pos) {$\entry$};
			}
			\draw[black!65,facet] (0,-3.11) -- (0.38,-3.11);
			\node[anchor=west] at (0.43,-3.11) {Facet};
		\end{scope}
		% Inset coordinates are w=156(s+11/26), z=156t.
		\begin{scope}[shift={(3.53,-0.81)},x=1.38cm,y=0.82cm]
			\node[anchor=west,inner sep=0pt] at (-0.75,1.07) {Detail at $t=0$};
			\begin{scope}
				\clip (-0.75,-0.79) rectangle (1.75,0.90);
				\fill[blue!8] (-0.75,2) -- (-0.75,1.5) -- (1.75,{-4.5/7})
					-- (1.75,2) -- cycle;
				\draw[zfive,candidate] (-0.75,0.65) -- (1.75,{-13*1.75/15});
				\draw[zthree,facet] (-0.75,1.5) -- (1.75,{-4.5/7});
				\draw[black!65] (-0.75,0) -- (1.75,0);
				\draw[zfive,fill=white] (0,0) circle (1.65pt);
				\fill[zthree] (1,0) circle (1.65pt);
				\node[below left,inner sep=1.5pt] at (0,0) {$-\tfrac{11}{26}$};
				\node[below left,inner sep=1.5pt] at (1,0) {$-\tfrac5{12}$};
			\end{scope}
			\draw[black!45,thin] (-0.75,-0.79) rectangle (1.75,0.90);
		\end{scope}
	\end{tikzpicture}\par
	}
	The point $(-\tfrac{41}{26},1)$ lies in the relative interior of the $E_5$-facet.
	Its supporting line $26s+30t+11=0$ meets $t=0$ at $s=-\tfrac{11}{26}$,
	so Theorem~\ref{thm:facet-specialization} gives the asserted root.
	\end{example}

	\subsection{Proof of Theorem~\ref{thm:specialized-BS}}
	Suppose $H$ is remarkable at $x$. By Theorem~\ref{thm:witness}, there is a
	component $Q$ of $B_Z$ with $x\in c(Q)$ and $H_{\mathbb C}=H_Q$.
	The tuples $F$ and $F^M$ have the same reduced zero divisor, hence the
	same log canonical model. Since $\ord_Q g_i=\sum_j m_{ij}\ord_Q f_j$,
	the hyperplane of $Q$ for $F^M$ is $(\sigma_M^{-1}(H))_{\mathbb C}$.
	Theorem~\ref{thm:lcBF} applied to $F^M$ gives the local version of Theorem \ref{thm:specialized-BS}. The global assertion follows from the same reduction as at the beginning of Subsection \ref{subs_APA}.\hfill$\Box$\par


\begin{thebibliography}{BCKSZ26}

		\bibitem[AGLN25]{AGLN25}
		J.~\`Alvarez Montaner, M.~Gonz\'alez Villa, E.~Le\'on-Cardenal, and L.~N\'u\~nez-Betancourt,
		\emph{Bernstein--Sato polynomial and related invariants for meromorphic functions},
		Trans. Amer. Math. Soc. \textbf{378} (2025), no.~7, 4929--4954.

		\bibitem[Bla24]{Bla24}
		G.~Blanco,
		\emph{Topological roots of the Bernstein--Sato polynomial of plane curves},
		Amer. J. Math., to appear, \href{https://arxiv.org/abs/2406.09034}{arXiv:2406.09034} (2024).

		\bibitem[Bud15]{Bud15}
		N.~Budur,
		\emph{Bernstein--Sato ideals and local systems},
		Ann. Inst. Fourier \textbf{65} (2015), no.~2, 549--603.

	

		\bibitem[BPSZ26]{BPSZ}
		N.~Budur, E.~de Lorenzo Poza, Q.~Shi, and H.~Zuo,
		\emph{A remarkable subset of poles of the motivic zeta function},
		\href{https://arxiv.org/abs/2602.13508}{arXiv:2602.13508v2} (2026).

		\bibitem[BSZ26]{BSZ_SMC}
		N.~Budur, Q.~Shi, and H.~Zuo,
		\emph{The monodromy conjecture for remarkable poles},
		preprint (2026).

		\bibitem[BVVW24]{BVVW24}
		N.~Budur, R.~van der Veer, and A.~Van Werde,
		\emph{Estimates for zero loci of Bernstein--Sato ideals},
		Publ. Res. Inst. Math. Sci. \textbf{60} (2024), no.~3, 607--628.

		\bibitem[BVWZ21]{BVWZ}
		N.~Budur, R.~van der Veer, L.~Wu, and P.~Zhou,
		\emph{Zero loci of Bernstein--Sato ideals},
		Invent. Math. \textbf{225} (2021), no.~1, 45--72.

		\bibitem[CL08]{CL}
		R.~Cluckers and F.~Loeser,
		\emph{Constructible motivic functions and motivic integration},
		Invent. Math. \textbf{173} (2008), no.~1, 23--121.
		
			\bibitem[CKSZ26]{BCKSZ-data}
		M.~\v{C}mrlec, S.~Kurokboev, Q.~Shi, and H.~Zuo,
		\emph{Data and code for ``A remarkable subset of polar hyperplanes of the multivariable motivic zeta function''},
		GitHub repository (2026),
		\url{https://github.com/thusq20/counterexample_rational_function}.

		\bibitem[DGPS25]{Singular}
		W.~Decker, G.-M.~Greuel, G.~Pfister, and H.~Sch\"onemann,
		\emph{\textup{\textsc{Singular}} 4.4.1---A computer algebra system for polynomial computations},
		2025, \url{https://www.singular.uni-kl.de}.

		\bibitem[DL98]{DL98}
		J.~Denef and F.~Loeser,
		\emph{Motivic {I}gusa zeta functions},
		J. Algebraic Geom. \textbf{7} (1998), 505--537.

		\bibitem[DY26]{DY}
		D.~Davis and R.~Yang,
		\emph{Multivariate $V$-filtrations and the Strong Monodromy Conjecture for hyperplane arrangements},
		\href{https://arxiv.org/abs/2605.25314v4}{arXiv:2605.25314v4} (2026).

		\bibitem[ELM04]{ELM04}
		L.~Ein, R.~Lazarsfeld, and M.~Musta\c t\u a,
		\emph{Contact loci in arc spaces},
		Compos. Math. \textbf{140} (2004), no.~5, 1229--1244.

		\bibitem[GB00]{GB00}
		E.~R.~Garc\'ia Barroso,
		\emph{Sur les courbes polaires d'une courbe plane r\'eduite},
		Proc. London Math. Soc. (3) \textbf{81} (2000), no.~1, 1--28.

		\bibitem[GVL14]{GVL14}
		M.~Gonz\'alez Villa and A.~Lemahieu,
		\emph{The monodromy conjecture for plane meromorphic germs},
		Bull. Lond. Math. Soc. \textbf{46} (2014), no.~3, 441--453.

		\bibitem[GLM98]{GLM98}
		S.~M.~Gusein-Zade, I.~Luengo, and A.~Melle-Hern\'andez,
		\emph{Zeta functions for germs of meromorphic functions, and Newton diagrams},
		Funct. Anal. Appl. \textbf{32} (1998), no.~2, 93--99.

		\bibitem[Kol97]{Kol97}
		J.~Koll\'ar,
		\emph{Singularities of pairs},
		in \emph{Algebraic geometry---Santa Cruz 1995},
		Proc. Sympos. Pure Math. \textbf{62}, Part~1 (1997), 221--287.

		\bibitem[KM98]{KM98}
		J.~Koll\'ar and S.~Mori,
		\emph{Birational geometry of algebraic varieties},
		Cambridge Tracts in Mathematics, vol.~134, Cambridge University Press, 1998.

		\bibitem[Lic89]{Lic89}
		B.~Lichtin,
		\emph{Poles of $|f(z,w)|^{2s}$ and roots of the $B$-function},
		Ark. Mat. \textbf{27} (1989), 283--304.

		\bibitem[Mum61]{Mum61}
		D.~Mumford,
		\emph{The topology of normal singularities of an algebraic surface and a criterion for simplicity},
		Publ. Math. Inst. Hautes \'Etudes Sci. \textbf{9} (1961), 5--22.

		\bibitem[OX12]{OX}
		Y.~Odaka and C.~Xu,
		\emph{Log-canonical models of singular pairs and its applications},
		Math. Res. Lett. \textbf{19} (2012), no.~2, 325--334.

		\bibitem[Rai12]{Rai12}
		M.~Raibaut,
		\emph{Motivic Milnor fibers of a rational function},
		C. R. Math. Acad. Sci. Paris \textbf{350} (2012), no.~9--10, 519--524.

		\bibitem[Rai13]{Rai13}
		M.~Raibaut,
		\emph{Motivic Milnor fibers of a rational function},
		Rev. Mat. Complut. \textbf{26} (2013), no.~2, 705--734.

		\bibitem[Sab87]{Sab87}
		C.~Sabbah,
		\emph{Proximit\'e \'evanescente. II. \'Equations fonctionnelles pour plusieurs fonctions analytiques},
		Compositio Math. \textbf{64} (1987), no.~2, 213--241.

		\bibitem[Ta23]{Ta23}
		K.~Takeuchi,
		\emph{On a Bernstein--Sato polynomial of a meromorphic function},
		Nagoya Math. J. \textbf{251} (2023), 715--733.

		\bibitem[Vey90]{Vey90}
		W.~Veys,
		\emph{On the poles of Igusa's local zeta function for curves},
		J. Lond. Math. Soc. (2) \textbf{41} (1990), no.~1, 27--32.

		\bibitem[VZG17]{VZG17}
		W.~Veys and W.~A.~Z\'u\~niga-Galindo,
		\emph{Zeta functions and oscillatory integrals for meromorphic functions},
		Adv. Math. \textbf{311} (2017), 295--337.

		\bibitem[Yan83]{Yan83}
		T.~Yano,
		\emph{$b$-functions and exponents of hypersurface isolated singularities},
		in \emph{Singularities (Arcata, Calif., 1981)},
		Proc. Sympos. Pure Math. \textbf{40}, Part~2 (1983), 641--652.

		% References retained for earlier drafts.

		% \bibitem[AC75]{AC} N. A’Campo. \newblock{La fonction z\^eta d’une monodromie.} {\em Comment. Math. Helv.} 50: 233-248, 1975.

		% \bibitem[ACLM05]{Artal} E. Artal Bartolo, P. Cassou-Noguès, I. Luengo, A. Melle Hernández. Quasi-ordinary power series and their zeta functions. Mem. Amer. Math. Soc. 178, vi+85 pp. (2005).

		% \bibitem[BBNV25]{Birational zeta function}
		% T. Biesbrouck, N. Budur, J. Nicaise,  W. Veys.
		% \newblock Birational zeta functions.
		% \newblock{arxiv:2509.03352}.

		% \bibitem[Bla22]{Bla22} G. Blanco.	\newblock{An algorithm for Hodge ideals}. {\tt hodge.lib}.		\newblock{\em Math. Comp.}, 338: 2955--2967.

		% \bibitem[BK86]{BK} E. Brieskorn, H. Kn\"orrer.
		% {Plane algebraic curves.} Birkh\"auser/Springer,  x+721 pp. (1986).

		% \bibitem[B03]{Bu} N. Budur. On Hodge spectrum and multiplier ideals.  Math. Ann.  327,  257-270 (2003).

		% \bibitem[B12]{Bsurvey} N. Budur. Singularity invariants related to Milnor fibers: survey. Contemp. Math. 566, 161-187 (2012).

		% \bibitem[BGG12]{BGV}	N. Budur, P. González-Pérez, M. González Villa.
		% Log canonical thresholds of quasi-ordinary hypersurface singularities.
		% Proc. Amer. Math. Soc. 140, 4075-4083 (2012).

		% \bibitem[B+24]{BBPP24}
		% N. Budur, J. de la Bodega, E. de Lorenzo Poza, J. Fern\'andez de Bobadilla, T. Pe\l ka.
		% \newblock{On the embedded Nash problem}.
		% \newblock{Forum Math. Pi} 12, Paper No. e15, 28 pp. (2024).

		% \bibitem[BFLN22]{Cohomology_of_Contact_Loci}
		% N. Budur, J.~Fern\'andez de~Bobadilla, Q.T. L\^e,  H.D.
		% Nguyen.
		% \newblock Cohomology of contact loci.
		% \newblock {J. Differential Geom.} 120, 389-409 (2022).

		% \bibitem[BPSZ26]{BPSZ} N. Budur, E. de Lorenzo Poza, Q. Shi,  H. Zuo. \newblock{A remarkable subset of poles of the motivic zeta function}. \newblock{ arXiv:2411.01305.}

		% \bibitem[BSZ24]{BSZ} N. Budur, Q. Shi,  H. Zuo. \newblock{Motivic principal value integrals for hyperplane arrangements}. \newblock{ arXiv:2411.01305.}

		% \bibitem[BN20]{BN20}
		% E. Bultot, J. Nicaise.
		% \newblock{ Computing motivic zeta functions on log smooth models.}
		% \newblock{ Math. Zeit.} 295, 427-462 (2020).

		% \bibitem[BL23]{BL}
		% J. de la Bodega, E. de Lorenzo Poza.
		% \newblock{The arc-Floer conjecture for plane curves}.
		% \newblock{arXiv:2308.00051}. To appear in J. Differential Geom.

		% \bibitem[DL99]{DL99} J. ~Denef,  F. ~Loeser. \newblock {Germs of arcs on singular algebraic varieties and motivic integration}.	\newblock {\em Invent. Math.}, 135(1):201--232, 1999.

		% \bibitem[ELSV04]{ELSV04}
		% L.~Ein, R.~Lazarsfeld, K.~E.~Smith, and D.~Varolin,
		% \emph{Jumping coefficients of multiplier ideals},
		% Duke Math. J. \textbf{123} (2004), no.~3, 469--506.

		% \bibitem[FP24]{BP} J. Fern\'andez de Bobadilla, T. Pe\l ka. Symplectic monodromy at radius zero and equimultiplicity of $\mu$-constant families. { Ann. Math.} 200, 153-299 (2024).

		% \bibitem[Fi12]{Fic} G. Fichou. The motivic real Milnor fibres. Manuscripta Math. 139, 167–178 (2012).

		% \bibitem[GG14]{GV}  P. González Pérez, M. González Villa.
		% Motivic Milnor fiber of a quasi-ordinary hypersurface. J. Reine Angew. Math. 687, 159-205 (2014).

		% \bibitem[GR25]{GP} P. González Pérez, M. Robredo Buces. Quasi-ordinary hypersurfaces, multiplier ideals and local tropicalizations. arXiv:2510.21009.

		% \bibitem[I18]{Ish14}
		% S. Ishii.
		% \newblock{Introduction to singularities}.
		% \newblock Springer, 2018, x+236 pp.

		% \bibitem[KT19]{KT}
		% M. Kontsevich, Y. Tschinkel.
		% \newblock{Specialization of birational types.}
		% \newblock{ Invent. Math.} 217, 415-432 (2019).

		% \bibitem[Li89]{Li} B. Lichtin. Poles of $|f(z,w)|^{2s}$ and roots of the $b$-function. Ark. Mat. 27, 283-304 (1989).

		% \bibitem[Lo88]{L} F. Loeser. Fonctions d'Igusa $p$-adiques et polyn\^{o}mes de
		% Bernstein. Amer. J. Math. 110, 1-21 (1988).

		% \bibitem[Ma20]{Man} F. Mangolte. Real algebraic varieties.
		% Springer, 2020, xviii+444 pp.

		% \bibitem[McCP03]{McP} C. McCrory, A. Parusi\'nski. Virtual Betti numbers of real algebraic varieties. C. R. Math. Acad. Sci. Paris 336, 763–768 (2003).

		% \bibitem[McL19]{McL19}
		% M. McLean.
		% \newblock Floer cohomology, multiplicity and the log canonical threshold.
		% \newblock { Geom. Topol.} 23, 957-1056 (2019).

		% \bibitem[Mu02]{Mus02}
		% M. Musta\c t\u a.
		% \newblock{ Singularities of pairs via jet schemes.}
		% \newblock{ J. Amer. Math. Soc.} 15, 599-615 (2002).

		% \bibitem[NX16]{NX16}
		% J. Nicaise, C. Xu.
		% \newblock Poles of maximal order of motivic zeta functions.
		% Duke Math. J. 165, 217-243 (2016).

		% \bibitem[RV03]{RV} B. Rodrigues and W. Veys.
		% Poles of zeta functions on normal surfaces.
		% { Proc. London Math. Soc.} 87, 164-196, 2003.

		% \bibitem[S00]{Sexp} M. Saito. Exponents of an irreducible plane curve singularity. arXiv:math/0009133.

		% \bibitem[S07]{S-real} M. Saito.  On real log canonical thresholds. arXiv:0707.2308.

		% \bibitem[Va76]{Var76}
		% A. Varchenko.
		% \newblock Zeta-function of monodromy and Newton's diagram.
		% \newblock{ Invent. Math.} 37, 253-262 (1976).

		% \bibitem[Va82]{Va} A. Varchenko. The complex singularity index does not change along the stratum $\mu=$const. Funct. Anal. Appl. 16, 1-9  (1982).

		% \bibitem[Laz04]{Laz04}
		% R.~Lazarsfeld,
		% \emph{Positivity in algebraic geometry II: Positivity for vector bundles, and multiplier ideals},
		% Ergebnisse der Mathematik und ihrer Grenzgebiete, 3. Folge, vol.~49,
		% Springer-Verlag, Berlin, 2004.

		% \bibitem[Ve95]{Ve} W. Veys. Determination of the poles of the topological zeta function for curves. Manuscripta
		% Math 87, 435-448 (1995).

	\end{thebibliography}
\end{document}